\documentclass{authesis}

\usepackage{amsmath}
\usepackage{amsfonts}
\usepackage{amssymb}
\usepackage{amscd}
\usepackage[matrix,arrow,cmtip,rotate,curve]{xy}
\SelectTips{cm}{}
\input xy
\xyoption{all}
\usepackage{graphics}
\usepackage{graphicx}
\usepackage{fancyhdr}
\usepackage{ushort}

\renewcommand{\thesection}{\thechapter.\arabic{section}}

\DeclareSymbolFont{bbold}{U}{bbold}{m}{n}
\DeclareMathSymbol{\numberone}{\mathord}{bbold}{`1}

\newcommand{\bbR}{\mathbb R}
\newcommand{\bbZ}{\mathbb Z}

\newcommand{\mcA}{\mathcal A}
\newcommand{\mcB}{\mathcal B}
\newcommand{\mcC}{\mathcal C}

\newcommand{\mcF}{\mathcal F}
\newcommand{\mcI}{\mathcal I}
\newcommand{\mcJ}{\mathcal J}

\newcommand{\mcP}{\mathcal P}

\newcommand{\mcR}{\mathcal R}
\newcommand{\mcS}{\mathcal S}
\newcommand{\mcU}{\mathcal U}
\newcommand{\mcV}{\mathcal V}
\newcommand{\mcW}{\mathcal W}

\newcommand{\mfX}{\mathfrak X}

\newcommand{\mfm}{\mathfrak m}
\newcommand{\mfn}{\mathfrak n}

\newcommand{\scE}{\textsc{E}}
\newcommand{\scM}{\textsc{M}}
\newcommand{\scU}{\textsc{U}}
\newcommand{\scV}{\textsc{V}}
\newcommand{\scS}{\textsc{S}}

\newcommand{\Ab}{\mathsf{Ab}}
\newcommand{\AffRmod}{\mathsf{Aff}\textrm{-}\mcR\textrm{-}\mathsf{mod}}
\newcommand{\Alg}{\mathsf{Alg}}
\newcommand{\Amod}{\mcA\textrm{-}\mathsf{mod}}
\newcommand{\Bmod}{\mcB\textrm{-}\mathsf{mod}}

\newcommand{\Ch}{\mathsf{Ch}}
\newcommand{\ConvTop}{\mathsf{ConvTop}}

\newcommand{\End}{\mathsf{End}}

\newcommand{\LCTVS}{\mathsf{LCTVS}}
\newcommand{\Rmod}{\mcR\textrm{-}\mathsf{mod}}
\newcommand{\Set}{\mathsf{Set}}
\newcommand{\Top}{\mathsf{Top}}
\newcommand{\TVS}{\mathsf{TVS}}

\newcommand{\Bilin}{\mathop\mathsf{Bilin}}

\newcommand{\map}{\mathop\mathsf{map}}

\newcommand{\xra}[1]{\xrightarrow{#1}}

\newcommand{\ra}{\rightarrow}

\newcommand{\xlra}[1]{\xrightarrow{\ #1\ }}
\newcommand{\xlla}[1]{\xleftarrow{\ #1\ }}
\newcommand{\lra}{\longrightarrow}

\newdir{c}{{}*!/-5pt/@^{(}}
\newdir{d}{{}*!/-5pt/@_{(}}
\newdir{s}{{}*!/+10pt/@{}}

\newcommand{\pullback}{\ar@{}[rd]|(.3){\scalebox{1.5}{$\displaystyle\lrcorner$}}}
\newcommand{\pushout}{\ar@{}[lu]|(.3){\scalebox{1.5}{$\displaystyle\ulcorner$}}}

\newcommand{\Cech}{\v Cech }
\newcommand{\ch}{\mathop\mathrm{ch}}
\newcommand{\Charts}{\mathrm{Charts}}
\newcommand{\colim}{\mathop\mathrm{colim}}
\newcommand{\Cov}{\mathop\mathrm{Cov}}
\newcommand{\dirsp}{\mathrm{dir}\ }
\newcommand{\Diff}{\mathrm{Diff}}

\newcommand{\hofib}{\mathrm{hofib}}
\newcommand{\Hom}{\mathop\mathrm{Hom}\nolimits}
\newcommand{\SurHom}{\mathop\mathrm{SurHom}\nolimits}
\newcommand{\id}{\mathrm{id}}
\newcommand{\im}{\mathop\mathrm{im}}

\newcommand{\Lan}{\mathop\mathrm{Lan}\nolimits}
\newcommand{\Nbhd}{\mathop\mathrm{Nbhd}}
\newcommand{\Op}{\mathop\mathrm{Op}}
\newcommand{\Pre}{\mathsf{Pr}}
\newcommand{\PreSheaf}[2]{\mathop\mathrm{{#2}^{\Op({#1})^{op}}}}
\newcommand{\Sh}{\mathsf{Sh}}
\newcommand{\sign}{\mathrm{sign}}
\newcommand{\spec}{\mathrm{sp}}
\newcommand{\supp}{\mathop\mathrm{supp}}
\newcommand{\tensor}[1]{\otimes_{#1}}
\newcommand{\Tot}{\mathrm{Tot}}

\newcommand{\tr}{\mathrm{tr}}
\newcommand{\vol}{\mathop\mathrm{vol}}
\newcommand{\wspec}{\mathrm{wsp}}

\begin{document}



\Title{%
A generalization of\\
Vassiliev's h-principle\\
}


\Author{Luk\'a\v s Vok\v r\'inek}


\Year{September 2006}


\phdfrontpage

\begin{Abstract}%

This thesis consists of two parts which share only a slight
overlap.

The first part is concerned with the study of ideals in the ring $C^\infty(M,\bbR)$ of smooth functions on a compact smooth manifold $M$ or more generally submodules of a finitely generated $C^\infty(M,\bbR)$-module $\mcV$. We define 
a topology on the space $\scM_d(\mcV)$ of all submodules of $\mcV$ of a fixed finite codimension $d$. Its main property is that it is compact Hausdorff and, when $\mcV=C^\infty(M,\bbR)$, it contains as a subspace the configuration space of $d$ distinct unordered points in $M$ and therefore gives a ``compactification'' of this configuration space. We present a concrete description of this space for low codimensions.

The main focus is then put on the second part which is concerned
with a generalization of Vassiliev's h-principle. This principle in
its simplest form asserts that the jet prolongation map
$j^r:C^\infty(M,V)\ra\Gamma(J^r(M,V))$, defined on the space of
smooth maps from a compact manifold $M$ to a Euclidean space $V$ and
with target the space of smooth sections of the jet bundle
$J^r(M,V)$, is a cohomology isomorphism when restricted to certain
``nonsingular'' subsets (these are defined in terms of a certain
subset $R\subseteq J^r(M,V)$). Our generalization then puts this
theorem in a more general setting of topological
$C^\infty(M,\bbR)$-modules. As a reward we get a strengthening of
this result asserting that all the homotopy fibres have zero
homology.
\end{Abstract}


\StartThesis

\pagestyle{fancy} \fancyhf{}
\renewcommand{\chaptermark}[1]{\markboth{}{\thechapter.\ #1}}
\renewcommand{\sectionmark}[1]{\markboth{}{\thesection.\ #1}}
\rhead{\thepage} \lhead{\scshape\rightmark}
\renewcommand{\headrulewidth}{0.4pt}





\setcounter{chapter}{-1}
\chapter{Introduction}

The history of h-principles (homotopy principles) goes back to
Smale who classified in \cite{Smale} immersions of spheres into a
Euclidean space up to regular homotopy and to Hirsch who
generalized Smale's work to immersions $M\ra N$ between any
manifolds (in \cite{Hirsch1} for $\dim M<\dim N$ and in
\cite{Hirsch2} for $\dim M=\dim N$ with $M$ open). More elaborated
h-principles can be stated as follows. Let $R\subseteq J^r(M,N)$
(corresponding to jets of immersions in the preceding) and let us
consider the set
 \[\Gamma R:=\{s\in \Gamma(J^r(M,N))\ |\ \im s\subseteq R\}\]
of those sections of the jet bundle $J^r(M,N)\ra M$ whose image
lies in $R$. The subset
 \[\Gamma_{\mathrm{hol}}R:=\{j^rf\in \Gamma R\ |\ f\in C^\infty(M,N)\}\]
of holonomic sections can be clearly identified with
 \[\{f\in C^\infty(M,N)\ |\ j^r_xf\in R\ \forall x\in M\}\subseteq C^\infty(M,N)\]
and this identification is in fact a homeomorphism. A (parametric)
h-principle for $R$ generally asserts that the inclusion
$\Gamma_{\mathrm{hol}}R\hookrightarrow\Gamma R$ is a weak homotopy
equivalence. The statements about the set of immersions up to
regular homotopy are then translated into ones about $\pi_0$. In
\cite{Gromov} Gromov proved that the h-principle holds for any open
$\Diff M$-invariant $R$ provided that $M$ is open.

For the case of compact $M$ the situation is more complicated and in
fact the \mbox{h-principle} for immersions between manifolds of the
same dimension does not hold. A partial result for maps $M\ra V$
into a Euclidean space $V$ is given in \cite{Vassiliev}: Vassiliev
proves that if $R$ is open and if its complement is a semialgebraic
$\Diff M$-invariant subset of codimension at least $\dim M+2$ then
the inclusion $\Gamma_{\mathrm{hol}}R\ra\Gamma R$ is a cohomology
isomorphism. Moreover he constructs a spectral sequence converging
to the cohomology of these spaces.

The purpose of this thesis is to generalize this theorem in few
ways. Our proof is based upon interpolation theory and
transversality theory. Both these work in a more general setting
than $C^\infty(M,V)$ and allow us to construct Vassiliev's spectral
sequence for the homotopy fibres of the inclusion
$\Gamma_{\mathrm{hol}}R\ra\Gamma R$ proving that they are acyclic
(have zero homology). This is again under the assumption that the
codimension of the complement of $R$ in $J^r(M,V)$ is at least $\dim
M+2$.

Our setting is that of topologically finitely generated affine
$C^\infty(M,\bbR)$-modules (topological quotients of free
$C^\infty(M,\bbR)$-module of finite rank). A ``representation'' of
such a module $\mcV$ on an affine bundle $E\ra M$ is a special map
$\varphi:\mcV\ra\Gamma E$ which enables us to determine whether an
element $v\in\mcV$ ``lies in'' an open subset $R\subseteq E$: let
us call $v$ regular if $\varphi(v)\subseteq R$. Then we can
consider the subset of all regular $v$ and ask what its homotopy
type is (or homology groups in our case). The space
$\Gamma_{\mathrm{hol}}R$ is the space of regular elements for
 \[j^r:\mcW_{\mathrm{hol}}=C^\infty(M,V)\ra\Gamma(J^r(M,V))\]
and $\Gamma R$ the space of regular elements for
 \[\id:\mcW=\Gamma(J^r(M,V))\ra\Gamma(J^r(M,V))\]
Interesting examples of topologically finitely generated affine
$C^\infty(M,\bbR)$-modules arise as affine submodules of
$\mcW_{\mathrm{hol}}$ of functions whose $r$-jets along a
submanifold (for example along the boundary or at a point) are fixed
and the corresponding affine submodules of $\mcW$.

Returning to the general case we denote for simplicity by $A$ the
complement of $R$ in $E$, by $\varphi_x$ the composition
$\mcV\xra{\varphi}\Gamma E\xra{ev_x}E_x$ for $x\in M$ and by
$\mcV_A$ the set of regular elements. Our main theorem (Theorem
\ref{thm_main_theorem} in the main text) can be stated as follows:

\begin{stheorem}
Let $M$ be a compact smooth manifold, let $\alpha:\mcU\ra\mcV$ be an
affine $C^\infty(M,\bbR)$-homomorphism between two topologically
finitely generated affine \linebreak $C^\infty(M,\bbR)$-modules, let
$\varphi:\mcV\rightarrow\Gamma E$ be a ``representation'' on an
affine bundle $E\ra M$ and $A\subseteq E$ a manifold stratified
subset of codimension at least $\dim M+2$ such that outside the
set\footnote{This (a bit technical) condition is justified by the
example of maps with a fixed $r$-jet along a submanifold where one
could allow only jets in $R=E-A$.}
\[\hat{M}=\{x\in M\ |\ (\varphi\alpha)_x\textrm{ is surjective}\}\]
we have $\im\varphi_x\cap A=\emptyset$. Then each homotopy fibre
$\hofib_v\alpha_A$ of the restriction $\alpha_A:\mcU_A\ra\mcV_A$
of $\alpha$ to the sets of regular elemets is acyclic, i.e.
$\tilde{H}_*(\hofib_v\alpha_A)=0$.
\end{stheorem}

Now for the actual contents of the thesis. The first chapter is
concerned with the interpolation on smooth manifolds. All that is
required for further chapters is the first section which
guarantees an existence of a finite dimensional linear subspace in
any finitely generated $C^\infty(M,\bbR)$-module $\mcV$ that is
transverse to all submodules of a fixed finite codimension. In the
case of the ring $C^\infty(M,\bbR)$ itself we have special
submodules/ideals for points $x_1,\ldots,x_n\in M$:
 \[\{f:M\ra\bbR\ |\ j^r_{x_1}f=\cdots=j^r_{x_n}f=0\}\]
For a subspace $D\subseteq C^\infty(M,\bbR)$ to be transverse to
this ideal is equivalent to the multi-jet evaluation map
 \[(j^r_{x_1},\ldots,j^r_{x_n}):C^\infty(M,\bbR)\ra J^r_{x_1}(M,\bbR)\times\cdots\times J^r_{x_n}(M,\bbR)\]
being surjective when restricted to $D$, i.e. given any collection
of $r$-jets at points $x_1,\ldots,x_n$ there is a function lying in
$D$ that realizes them. This is the way the interpolation property
is used in the main proof.

The remainder of Chapter \ref{interpolation} is devoted to studying
the space of submodules of $\mcV$ of a fixed finite codimension.
Generalizing Glaeser's article \cite{Glaeser} it is given a
canonical topology making it into a compact Hausdorff space. We
present a few examples showing what this topology looks like. We
also include a criterion for an $C^\infty(M,\bbR)$-module to be
topologically finitely generated justifying the above example of
submodules with fixed $r$-jets along a submanifold.

The second and the third chapters are preparing the ground for the
construction of the main spectral sequence.

In the fourth chapter, assuming that the codimension of $A$ in $E$
is at least $\dim M+1$, we derive the main spectral sequence for
the homology of $\mcV_A$ and show that it converges if this
codimension is at least $\dim M+2$. As promised we use the
interpolation theory of (the first section of) Chapter
\ref{interpolation} and the transversality theory of Chapter
\ref{transversality}.

In Chapter \ref{homotopy_fibre} we construct the homotopy fibre of
the map $\alpha_A$ as the space $\mcW_{A\times I}$ of regular
elements in some topologically finitely generated affine
$C^\infty(M\times I,\bbR)$-module $\mcW$ equipped with a
representation $\mcW\ra\Gamma(E\times I)$ on the bundle $E\times
I\ra M\times I$. This allows us to apply our spectral sequence on
it and finally in Chapter \ref{main_theorem} prove our main
result.


\chapter{Interpolation on smooth manifolds} \label{interpolation}

Let us first state clearly that this chapter is based on the article
\cite{Glaeser} of Glaeser. In a sense it is just a (nontrivial)
generalization of the ideas of this article from parallelepiped in
$\bbR^m$ to general compact manifolds. The main structure of the
proof of the compactness of the space of ideals (or more generally
submodules) remains the same although at some point one has to come
up with a new approach as Glaeser's proof is very ``affine''. In the
first section we gather well-known facts about the ideals in the
algebra of smooth functions. The results in subsequent sections are
original.

In this chapter let $M$ be a compact smooth manifold,
$\mcR=C^\infty(M,\bbR)$ the ring of smooth functions on $M$. We
will be considering finitely generated $\mcR$-modules $\mcV$ and
their submodules of a fixed finite codimension $d$ over $\bbR$.
Our ultimate goal is to endow the set $\scM_d$ of all such
submodules with a topology. It will be compact Hausdorff. Together
with $\scM_d$ we will also topologize
\[\scE_d=\{(B,w)\ |\ B\in\scM_d,\ w\in\mcV/B\}\]
in such a way that the canonical projection map $\scE_d\ra\scM_d$
will become a vector bundle. If $\mcV$ is a topological
$\mcR$-module and it is topologically finitely generated then the
map
\[\mcV\times\scM_d\lra\scE_d\]
sending $(v,B)$ to $(B,v+B)$ is a continuous homomorphism of
vector bundles. This is certainly a property one would require
from any such topology.

To demonstrate other properties let us now specialize to the case
$\mcV=\mcR$ so that $\scM_d$ is the set of all ideals of $\mcR$ of
codimension $d$. For any set $Y\subseteq M$ consisting of $d$
points one has an ideal
\[\mfm_Y:=\{f\in\mcR\ |\ f(y)=0\textrm{ for all }y\in Y\}\]
and easily $\mfm_Y\in\scM_d$. In this way one can embed into
$\scM_d$ the configuration space $M^{[d]}$ of $d$ distinct
unordered points in $M$
\[M^{[d]}\hookrightarrow\scM_d\]
For our topology on $\scM_d$ this inclusion map is a topological
embedding and therefore one can think of $\scM_d$ as a
``compactification'' of $M^{[d]}$.

There is an inverse procedure of associating to each $I\in\scM_d$
an unordered $d$-tuple of points in $M$. However for ideals not of
the form $\mfm_Y$ these points do not have to be distinct and
therefore this procedure yields a map
\[\wspec:\scM_d\ra M^d/\Sigma_d\]
We call $\wspec(B)$ the weighted spectrum of $B$. For any reasonable
topology on $\scM_d$ this map should be continuous as well. This is
indeed the case for our topology.

Let us try to indicate now how we construct the topology on
$\scM_d$. If $\mcV$ was finite dimensional we could simply think
of $\scM_d$ as a subset of the Grassmannian manifold $G_d(\mcV)$
of linear subspaces of $\mcV$ of codimension $d$. This is of
course almost never the case. On the other hand suppose that there
is a finite dimensional linear subspace $D\subseteq\mcV$ that is
transverse to each $B\in\scM_d$. Then the intersection with $D$
produces a map
\[\scM_d\lra G_d(D)\]
from $\scM_d$ to the Grassmannian manifold of codimension $d$
linear subspaces of $D$. If this map is moreover injective one can
give $\scM_d$ the subspace topology. The key step is therefore to
show that there is a transversal $D$ for which the map $\scM_d\lra
G_d(D)$ is injective. This is the first section and it is rather
elementary. It starts by describing the structure of submodules of
a given $\mcR$-module $\mcV$ so that one is able to understand
this transversality condition. At the end of the section we
produce a transversal.

In the second section we derive an ``interpolation formula''. The
idea is that if $D$ is transverse to some submodule $B\in\scM_d$
then for each $v\in\mcV$ there is some $d\in D$ such that $v=d$
modulo $B$. Therefore one can ``interpolate'' elements of $\mcV$
``at $B$'' by elements of $D$. Moreover if we assume that the
dimension of $D$ is precisely $d$ there is exactly one such
interpolation. Intuitively $D$ should then be also transverse to
all submodules that are close to $B$ and one should get a
continuous interpolation map
\[\mcV\times\Nbhd(B)\rightarrow D\]
A problem is that $\scM_d$ does not posses any topology so far and
so we cannot talk about a neighbourhood of $B$. However we have a
topology on $M^d/\Sigma_d$. Using locally the affine structure of
$M$ we construct a continuous interpolation map\footnote{One has to
specify what an interpolation property with respect to elements of
$M^d/\Sigma_d$ means. It forces taking a slightly bigger $D$ and as
a result there is no preferred/unique interpolation map.}
\[\mcV\times\Nbhd(Y)\rightarrow D\]
where $Y\in M^d/\Sigma_d$ is the weighted spectrum of $B$. It serves
as a tool both in the proof of the compactness of $\scM_d$ and in
the proofs of its main properties.

The third section is dedicated to describing the topology of
$\scM_d(\mcR)$ in some special cases (codimensions $1$, $2$ and $3$ and
the case of one-dimensional manifolds). Along with these examples
two general theorems are proved, the first of which states that
the inclusion of the configuration space into $\scM_d$ is a smooth
embedding (as was already mentioned above). The second theorem is
concerned with the set of all ideals $I$ with a fixed isomorphism
type of the quotient algebra $\mcR/I$. These ideals possess a
structure of a smooth manifold (as was proved in \cite{Alonso})
and the content of our second theorem is that they form an
immersed submanifold in $\scM_d$. As a set $\scM_d$ is then a
disjoint union of these submanifolds.

One of the important notions that turn up in the course of this
chapter is the following. An $\mcR$-module $\mcV$ is topologically
finitely generated if $\mcV$ is a topological quotient of a free
$\mcR$-module of finite rank (with its canonical topology). In the
fourth section we give a simple criterion for an $\mcR$-module to
be topologically finitely generated. It suffices for it to be
locally topologically finitely generated. As an application we
give some examples of topologically finitely generated modules.

Finally in the last section we briefly discuss the results for
non-compact manifolds. This is still a work in progress and as a
result almost no proofs are included.

\section{The structure of submodules of $\mcV$}

For the first few results we do not need to restrict to finitely
generated modules, hence we assume now that $\mcV$ is any
$\mcR$-module. We investigate the structure of submodules of
$\mcV$, first looking at the special case of ideals of the ring
$\mcR$ itself.

\begin{lemma}
The maximal ideals in $\mcR$ are identified with points in $M$.
For a point $x\in M$ the corresponding ideal is
\[\mfm_x:=\{f:M\rightarrow\bbR\ |\ f(x)=0\}\]
\end{lemma}

\begin{proof}
To prove that $\mfm_x$ is maximal it is enough to observe that it
is the kernel of the evaluation homomorphism
$ev_x:\mcR\rightarrow\bbR$ and that $\bbR$ is a field.

On the other hand suppose that $I$ is a proper ideal and we will
prove that there is $x\in M$ such that $I\subseteq \mfm_x$. Assume
on the contrary that no such point exists. Then for every point
$x$ there is a function $f_x\in I$ such that $f_x\neq 0$ on a
neighbourhood $U_x$ of $x$. Replacing $f_x$ by $f_x^2$ if needed
we can assume that each $f_x$ is nonnegative. By compactness there
is a finite set $\{x_1,\ldots,x_k\}\subseteq M$ such that
$U_{x_1},\ldots,U_{x_k}$ cover $M$. Therefore
$f=f_{x_1}+\cdots+f_{x_k}\in I$ is nonzero on $M$ and consequently
a unit in $\mcR$, a contradiction to the properness of $I$.
\end{proof}

For any ideal $I\subseteq\mcR$ the \textit{spectrum} of $I$ is
defined to be the closed subset
\[\spec(I):=\bigcap_{f\in I}f^{-1}(0)\]
The previous lemma then says that it can also be described as the
subset of those points $x\in M$ for which $I\subseteq \mfm_x$. For
any closed subset $A\subseteq M$ we denote by $\mfn_A$ the ideal
of functions which vanish on a neighborhood of $A$. The quotient
$\mcR/\mfn_A$ is then the ring of germs at $A$ of smooth functions
on $M$.

\begin{lemma} \label{int_adjointness_for_ideals}
Let $I\subseteq\mcR$ be an ideal and $A\subseteq M$ a closed
subset. Then
\[\spec(I)\subseteq A\quad\Leftrightarrow\quad \mfn_A\subseteq I\]
\end{lemma}

\begin{remark}
If we define $\mfm_A$ to be the ideal of functions vanishing on
$A$ then also $A\subseteq\spec(I)$ if and only if $I\subseteq
\mfm_A$. Both these statements can be phrased as adjointness of
the respective functors.
\end{remark}

\begin{proof}
Let $f\in \mfn_A$ be any function $M\rightarrow\bbR$ vanishing on
a neighborhood of $A$. We set $C=\supp(f)$ and we have
\[C\subseteq M-A\]
We can cover $C$ by open subsets $U_i$ for which there is $g_i\in
I$ that is positive on $U_i$. By the means of a partition of unity
we can glue them to get a function $g\in I$ which is positive on a
neighborhood of $C$. Then $f=\frac{f}{g}g$ expresses $f$ as an
element of $I$.
\end{proof}

\begin{lemma} \label{int_properties_of_JA}
Let $A$ and $A'$ be closed subsets of $M$. Then the following
holds
\begin{itemize}
\item[$(i)$]{$\mfn_A\cap \mfn_{A'}=\mfn_{A\cup A'}$}
\item[$(ii)$]{$\mfn_A+\mfn_{A'}=\mfn_{A\cap A'}$}
\end{itemize}
\end{lemma}

\begin{proof}
The part $(i)$ is obvious and so is the inclusion $\subseteq$ in
the part $(ii)$.
Hence let $f\in \mfn_{A\cap A'}$, $\hat{A}=A\cap\supp(f)$ and
$\hat{A}'=A'\cap\supp(f)$. As $\hat{A}\cap\hat{A}'=\emptyset$ we can
find a function $\lambda\in\mcR$ such that $\lambda=0$ on a
neighborhood
of $\hat{A}$ and $\lambda=1$ on a neighborhood of $\hat{A}'$%
. Then clearly $\lambda f\in \mfn_A$ and $(1-\lambda)f\in
\mfn_{A'}$ and so
\[f=\lambda f+(1-\lambda)f\in \mfn_A+\mfn_{A'}\]
\end{proof}

\begin{lemma} \label{int_intersection_versus_product}
Let $I$ and $I'$ be ideals with disjoint spectra and
$B\subseteq\mcV$ any submodule. Then
\begin{itemize}
\item[$(i)$]{$II'=I\cap I'$}
\item[$(ii)$]{$II'\mcV=I\mcV\cap I'\mcV$}
\item[$(iii)$]{$B+II'\mcV=(B+I\mcV)\cap(B+I'\mcV)$}
\end{itemize}
\end{lemma}

\begin{proof}
Obviously $(iii)$ is the most general case. We give here only the
proof of $(ii)$ for the simplicity sake. Let us denote the spectra
of $I$ and $I'$ by $A$ and $A'$ respectively. The inclusion
$\subseteq$ is obvious and
 {\setlength\arraycolsep{2pt}
 \begin{eqnarray*}
  I\mcV\cap I'\mcV & = & (\mfn_A+\mfn_{A'})(I\mcV\cap
  I'\mcV)=\mfn_A(I\mcV\cap I'\mcV)+\mfn_{A'}(I\mcV\cap I'\mcV)\\
  & \subseteq & \mfn_AI'\mcV+\mfn_{A'}I\mcV\subseteq II'\mcV
  \end{eqnarray*}}
(the first equality follows from Lemma \ref{int_properties_of_JA}
and the last inclusion from Lemma \ref{int_adjointness_for_ideals}).
The proof for $B\neq 0$ follows the same idea.
\end{proof}

Let $B\subseteq\mcV$ be any submodule. One has an
ideal
$(B:\mcV)$ of $\mcR$ defined by
\[(B:\mcV):=\{f\in\mcR\ |\ f\mcV\subseteq B\}\]
Alternatively, it is the kernel of the action map
\[\mcR\rightarrow\End(\mcV/B)\]
We define a \textit{spectrum} of $B$ (inside $\mcV$) to be the
spectrum of the ideal $(B:\mcV)$ and denote it by $\spec_\mcV(B)$
or simply by $\spec(B)$.

\begin{lemma} \label{int_adjointness_for_submodules}
Let $B\subseteq\mcV$ be a submodule and $A\subseteq M$ a closed
subset. Then
 \[\spec(B)\subseteq A\quad\Leftrightarrow\quad \mfn_A\mcV\subseteq B\]
\end{lemma}

\begin{proof}
This is clear as $\mfn_A\mcV\subseteq B$ if and only if
$\mfn_A\subseteq(B:\mcV)$ if and only if
$\spec(B)=\spec(B:\mcV)\subseteq A$ according to Lemma
\ref{int_adjointness_for_ideals}.
\end{proof}

Our next goal is to decompose and thus simplify the submodules of
$\mcV$. First we need yet another lemma.

\begin{lemma} \label{int_properties_of_spectra}
Let $B$ and $B'$ be submodules of $\mcV$. Then the following
holds.
\begin{itemize}
\item[$(i)$]{If $B\subseteq B'$ then $\spec(B')\subseteq\spec(B)$}
\item[$(ii)$]{$\spec(B\cap B')=\spec(B)\cup\spec(B')$}
\item[$(iii)$]{$\spec(B+B')\subseteq\spec(B)\cap\spec(B')$}
\end{itemize}
\end{lemma}

\begin{proof}
To prove $(i)$ one observes that $(B:\mcV)\subseteq(B':\mcV)$ and
so
\[\spec(B')=\spec(B':\mcV)\subseteq\spec(B:\mcV)=\spec(B)\]
For the inclusion $\subseteq$ of $(ii)$ let us denote $A=\spec(B)$
and $A'=\spec(B')$ so that (using Lemma
\ref{int_adjointness_for_submodules} in both directions)
 \[\mfn_{A\cup A'}\mcV\subseteq \mfn_A\mcV\cap \mfn_{A'}\mcV\subseteq B\cap
 B'\quad\Rightarrow\quad \spec(B\cap B')\subseteq A\cup A'\]
The remaining claims are trivial consequences of $(i)$.
\end{proof}


\begin{corollary}
If $B_1,\ldots,B_n$ are submodules of $\mcV$ with pairwise disjoint
spectra, then they are in general position, i.e.
\[B_i+(B_1\cap\cdots\cap\widehat{B_i}\cap\cdots\cap B_n)=\mcV\]
\qed\end{corollary}

Let $B\subseteq\mcV$ be any submodule whose spectrum is a disjoint
union $A=A_1\cup\cdots\cup A_n$ of closed subsets $A_i\subseteq
M$. Then according to Lemma \ref{int_intersection_versus_product},
Lemma \ref{int_properties_of_JA} and Lemma
\ref{int_adjointness_for_submodules} one has
\[(B+\mfn_{A_1}\mcV)\cap\cdots\cap(B+\mfn_{A_n}\mcV)=B+\mfn_A\mcV=B\]
We observe that (by Lemma \ref{int_properties_of_spectra})
$\spec(B+\mfn_{A_i}\mcV)\subseteq\spec(B)\cap A_i=A_i$ while on
the other hand $A=\spec(B)=\bigcup\spec(B+\mfn_{A_i}\mcV)$ so that
$\spec(B+\mfn_{A_i}\mcV)=A_i$ and the submodules
$B+\mfn_{A_i}\mcV$ are in general position. We will now prove the
uniqueness of such a decomposition.

\begin{lemma}
Let $B_1,\ldots,B_n$ be submodules of $\mcV$ with disjoint spectra
$A_1,\ldots,A_n$, let $B=B_1\cap\cdots\cap B_n$. Then
$B_i=B+\mfn_{A_i}\mcV$.
\end{lemma}

\begin{proof}
We have
\[B^i:=B+\mfn_{A_i}\mcV\subseteq B_i+\mfn_{A_i}\mcV=B_i\]
Suppose now that $B^j\subsetneqq B_j$ for some $j$. Because $B^i$
are in general position we get
\[B^j+(B^1\cap\cdots\cap\widehat{B^j}\cap\cdots\cap B^n)=\mcV\]
and therefore
\[B=B^1\cap\cdots\cap B^n\subsetneqq
B_j\cap(B^1\cap\cdots\cap\widehat{B^j}\cap\cdots\cap B^n)\subseteq
B_1\cap\cdots\cap B_n=B\] a contradiction.
\end{proof}

We will now apply these ideas to submodules of finite codimension
over $\bbR$. We fix an integer $d$ and denote the collection of
all the submodules of codimension $d$ by $\scM_d=\scM_d(\mcV)$.
First observe that the spectrum of any such submodule is finite.
This is obvious for ideals $I\subseteq\mcR$ as for any collection
$y_1,\ldots,y_n\in\spec(I)$ of distinct points one can find
functions $f_1,\ldots,f_n$ such that $f_i(y_i)=1$ and $f_i(y_j)=0$
if $i\neq j$. The $f_i$'s are obviously linearly independent and
hence their span is an $n$-dimensional linear subspace that
clearly intersect $I$ trivially. If $B\in\scM_d(\mcV)$ then
\[(B:\mcV)=\ker(\mcR\rightarrow\End(\mcV/B))\]
is an ideal of finite codimension and hence
$\spec(B)=\spec(B:\mcV)$ is finite.

Denoting the spectrum of $B$ by $Y=\{y_1,\ldots,y_n\}$ we get the
decomposition
\[B=(B+\mfn_{y_1}\mcV)\cap\cdots\cap(B+\mfn_{y_n}\mcV)\]
where the submodules $B+\mfn_{y_1}\mcV,\ldots,B+\mfn_{y_n}\mcV$
are in general position and hence
\[\dim\mcV/B=\dim\mcV/(B+\mfn_{y_1}\mcV)+\cdots+\dim\mcV/(B+\mfn_{y_n}\mcV)\]
We set $k_i=\dim\mcV/(B+\mfn_{y_i}\mcV)$ and we see that the
spectrum of $B$ has more structure if $B$ has finite codimension:
each point $y_i$ in the spectrum has associated a weight $k_i$
with it.

We define a space $\scS_d(M)$ (where these more structured spectra
will be defined) by
\[\scS_d(M):=M^d/\Sigma_d\]
and give it the quotient topology. Very often we abbreviate
$\scS_d(M)$ to $\scS_d$. The space $\scS_d$ consists of unordered
collections of $d$ not necessarily distinct points in $M$. If
$y_1,\ldots,y_n$ are all the points in such a collection $Y$ and if
each $y_i$ appears in it exactly $k_i$-times then we say that $k_i$
is the weight of $y_i$ and use an alternative notation
\[Y=\{(y_1,k_1),\ldots,(y_n,k_n)\}\]
We call $|Y|:=\{y_1,\ldots,y_n\}$ the \textit{support} of $Y$. There
is a weight function $|Y|\rightarrow\bbZ^+$ associating to each
point $y_i$ its weight $k_i$ and by the definition the total weight
$\sum k_i$ is $d$. Therefore we also call $Y$ a set of points with
weights.

Hence with every submodule $B\subseteq\mcV$ of finite codimension
$d$ there is associated a canonical set of points with weights
$Y\in\scS_d$ whose support is the spectrum of $B$ (and whose total
weight is $d$). It is called the \textit{weighted spectrum} of $B$.

On the other hand if $Y=\{(y_1,k_1),\ldots,(y_n,k_n)\}\in\scS_d$
we define an ideal $\mfm_Y\subseteq\mcR$ by
 \[\mfm_Y=(\mfm_{y_1})^{k_1}\cdots(\mfm_{y_n})^{k_n}=
 (\mfm_{y_1})^{k_1}\cap\cdots\cap(\mfm_{y_n})^{k_n}\]
(with the equality guaranteed by Lemma
\ref{int_intersection_versus_product}) and hence also get
submodules $\mfm_Y\mcV\subseteq\mcV$. They need not be of finite
codimension unless $\mcV$ is finitely generated. We will see later
that for a finitely generated $\mcV$ they do have a finite
codimension.

\begin{lemma}
Let $B\subseteq\mcV$ be a submodule of finite codimension with
weighted spectrum $Y$. Then $\mfm_Y\mcV\subseteq B$.
\end{lemma}

\begin{proof}
Let us first prove the lemma in the case $Y=\{(y,k)\}$. We set
$I=(B:\mcV)$ and note that by our assumptions $I$ is contained
only in one maximal ideal, namely in $\mfm_y$. We have a
$k$-dimensional $\mcR/I$-module $\mcV/B$ and hence
\[(\mfm_y)^{k+1}(\mcV/B)=(\mfm_y)^k(\mcV/B)\]
for dimensional reasons. Now we apply Nakayama's lemma (see for
example \cite{Lang}) to conclude that $(\mfm_y)^k(\mcV/B)=0$, i.e.
$(\mfm_y)^k\mcV\subseteq B$.

In general when $|Y|=\{y_1,\ldots,y_n\}$ we use the decomposition
\[B=(B+\mfn_{y_1}\mcV)\cap\cdots\cap(B+\mfn_{y_n}\mcV)\]
and Lemma \ref{int_intersection_versus_product} to reduce the
proof to the case $n=1$.
\end{proof}

\begin{lemma}
For $(\mfm_y)^k$ the following holds
\[(\mfm_y)^k=\{f\in\mcR\ |\ j^{k-1}_yf=0\}\]
and it is a finitely generated ideal.
\end{lemma}

\begin{proof}
The inclusion $\subseteq$ is obvious from the product formula for
the derivative. Hence let $f\in\mcR$ be such that $j^{k-1}_yf=0$.
Let $\lambda:M\rightarrow\bbR$ be a function which is supported in
a coordinate neighborhood around $y$ and is identically $1$ on a
neighborhood of $y$. Then
\[f=\lambda\cdot f+(1-\lambda)\cdot f\]
where the second summand is in $\mfn_y\subseteq(\mfm_y)^k$ so it
remains to show that the first summand $g=\lambda\cdot f$ lies in
the same. Note that $j^{k-1}_yg=0$ as $f$ and $g$ agree near $y$
and thus we can write in the coordinate chart
\[g(x)=\frac{1}{k!}\sum\limits_{i_1,\ldots i_k=1}^{\dim M}
a_{i_1\ldots i_k}(x)\cdot x_{i_1}\cdots x_{i_k}\] with
$a_{i_1\ldots i_k}$ smooth. If $\rho:M\rightarrow\bbR$ is any
function that is identically $1$ on a neighborhood of
$\supp(\lambda)$ and is supported in the above coordinate chart we
get smooth functions $\hat{x}_i=\rho\cdot x_i$ extended by $0$ to
$M$ and similarly $\hat{a}_{i_1\ldots i_k}$. Clearly on the whole
whole of $M$ we have
\[g(x)=\frac{1}{k!}\sum\limits_{i_1,\ldots i_k=1}^{\dim M}
\hat{a}_{i_1\ldots
i_k}(x)\cdot\hat{x}_{i_1}\cdots\hat{x}_{i_k}\in(\mfm_y)^k\] and
moreover we see that $(\mfm_y)^k$ is generated as an ideal by
$1-\lambda$ and the functions $\hat{x}_{i_1}\cdots\hat{x}_{i_k}$.
\end{proof}

\begin{note}
Last lemma is not true in the $C^r$ case. As an example, let
$M=\bbR$ and $r=1$, then any element of $\mfm_0$ can be written in
the form $f(x)=a(x)x$ with $a(x)$ continuous. Hence any element of
$(\mfm_0)^2$ can be written in the form $f(x)=a(x)x^2$ with $a(x)$
continuous again. In particular for any such function the limit
\[\lim_{x\rightarrow 0} \frac{f(x)}{x^2}=\lim_{x\rightarrow 0}a(x)=a(0)\]
exists. Certainly, the function $g(x)=x|x|$ is $C^1$ and
\[g(0)=g'(0)=0\]
but
\[\lim\limits_{x\rightarrow 0}\frac{g(x)}{x^2}\] does not
exist.
\end{note}

\begin{lemma} \label{int_spectral_ideals}
For any set of points with weights
$Y=\{(y_1,k_1),\ldots,(y_n,k_n)\}$ the following holds
\[\mfm_Y=\{f\in\mcR\ |\ j^{k_1-1}_{y_1}f=\cdots=j^{k_n-1}_{y_n}f=0\}\]
In other words $\mfm_Y$ is the kernel of
 \[(j^{k_1-1}_{y_1},\ldots,j^{k_n-1}_{y_n})^T:\mcR\lra
 J^{k_1-1}_{y_1}(M,\bbR)\times\cdots\times J^{k_n-1}_{y_n}(M,\bbR)\]
Moreover it has a finite codimension in $\mcR$ and it is finitely
generated as an ideal. \qed
\end{lemma}

To proceed further we need to specialize to the case of finitely
generated $\mcV$. If $E$ is a finite dimensional real vector
space, it induces a free $\mcR$-module $\mcR\otimes E\cong
C^\infty(M,E)$ which we will abbreviate to $\mcR E$. We will call
$\mcV$ \emph{topologically finitely generated} if there is a
surjective $\mcR$-module map $\mcR E\rightarrow\mcV$ that is a
continuous quotient map, i.e. such that the topology on $\mcV$ is
the quotient topology induced by this map. This notion will be
important later
.

%

\begin{corollary}
If $\mcV$ is finitely generated then so is every $B\in\scM_d$.
\end{corollary}

\begin{proof}
Let $Y$ be the weighted spectrum of $B$. It is clear that when
$\mcV$ is finitely generated then $\mfm_Y\mcV$ has a finite
codimension (an easy consequence of the case $\mcV=\mcR$ from
Lemma \ref{int_spectral_ideals}). The statement then follows from
the fact that $\mfm_Y\mcV$ is finitely generated and its
codimension in $B$ is finite.
\end{proof}

Let $D$ be a linear subspace of an $\mcR$-module $\mcV$. We write
$D\pitchfork\scM_d(\mcV)$ or just $D\pitchfork\scM_d$ if $D$ is
transverse to every element of $\scM_d$, i.e. to every submodule
of $\mcV$ of codimension $d$, and $D\pitchfork\scS_d(\mcV)$ or
just $D\pitchfork\scS_d$ if $D$ is transverse to all submodules
$\mfm_Y\mcV$, whenever $Y\in\scS_d$.

If $\mcV$ is finitely generated we construct a finite dimensional
linear subspace $D\pitchfork\scS_d$ (and hence also
$D\pitchfork\scM_d$). First we find this subspace locally in the
case $\mcV=\mcR$.

\begin{lemma} \label{int_local_existence_of_transversal_for_R}
Let $M=\bbR^m$. There is a finite dimensional linear subspace
$D\subseteq\mcR$ satisfying $D\pitchfork\scS_d$.
\end{lemma}

\begin{proof}
Let $D$ be the linear subspace of $\mcR$ of all polynomials
$\bbR^m\rightarrow\bbR$ of degree at most $d-1$ and let
$Y\in\scS_d$ be any set of points with weights. As the projection
$\mcR\rightarrow\mcR/\mfm_Y$ can be identified with the jet
evaluation (see Lemma \ref{int_spectral_ideals})
\[\mcR\lra J^{k_1-1}_{y_1}(\bbR^m,\bbR)\times\cdots\times
J^{k_n-1}_{y_n}(\bbR^m,\bbR)=:J_Y(\bbR^m,\bbR)\] we need to show
that the composition
\[D\subseteq\mcR E\rightarrow J_Y(\bbR^m,\bbR)\]
is surjective. We identify the jet spaces
$J^{k_i-1}_{y_i}(\bbR^m,\bbR)$ with the truncated polynomial
algebra of polynomials $\bbR^m\rightarrow\bbR$ of degree at most
$k_i-1$ in such a way that the polynomial $p$ corresponds to the
jet $j^{k_i-1}_{y_i}p(x-y_i)$. This translation ensures that the
identification is an isomorphism of algebras. To prove
surjectivity let
\[(0,\ldots,0,p_i,0,\ldots,0)\in J_Y(\bbR^m,\bbR)\]
We choose a linear form $\alpha:\bbR^m\rightarrow\bbR$ which is
injective on the support $|Y|$ and consider the following
polynomial
\[q(x)=\prod_{j\neq i}(\alpha(x)-\alpha(y_j))^{k_j}\]
By construction $q(y_i)\neq 0$ and thus $q(x+y_i)$ is invertible
in the truncated polynomial algebra. Hence we can find a
polynomial $r(x)$ of degree at most $k_i-1$ such that
$q(x+y_i)r(x)=p_i(x)$. This product corresponds to the jet
$j^{k_i-1}_{y_i}(q(x)r(x-y_i))$ and so $q(x)r(x-y_i)\in D$ is a
preimage of $(0,\ldots,0,p_i,0,\ldots,0)$.
\end{proof}

\begin{theorem} \label{int_global_existence_of_transversal_for_R}
Let $M$ be any compact manifold. Then there is a finite
dimensional linear subspace $D\subseteq\mcR$ for which
$D\pitchfork\scS_d$.
\end{theorem}

\begin{proof}
Let us choose an embedding $\iota:M\hookrightarrow V$ of $M$ into
a Euclidean space $V$ and note that the canonical map
 \[\iota^*:C^\infty(V,\bbR)\lra C^\infty(M,\bbR)\]
is surjective (i.e. any smooth map $M\ra\bbR$ is extensible to
$V$). By Lemma \ref{int_local_existence_of_transversal_for_R}
there exists a finite dimensional linear subspace $D\subseteq
C^\infty(V,\bbR)$ such that $D\pitchfork\scS_d(V)$. Then we claim
that $\iota^*D=\{f\iota\ |\ f\in D\}$ satisfies
$\iota^*D\pitchfork\scS_d(M)$. This is verified by the commutative
diagram
 \[\xymatrixnocompile{
  D \ar@{c->}[r] \ar[d]^{\iota^*} &
  C^\infty(V,\bbR) \ar@{->>}[r] \ar@{->>}[d]^{\iota^*} &
  J_{\iota(Y)}(V,\bbR) \ar@{->>}[d]^{\iota^*}
 \\
  \iota^*D \ar@{c->}[r] &
  C^\infty(M,\bbR) \ar@{->>}[r] &
  J_Y(M,\bbR)
 }\]
The composition across the top row is surjective and therefore so
must be the composition across the bottom one.
\end{proof}

\begin{theorem} \label{int_existence_of_transversal}
Let $M$ be a compact manifold and let $\mcV$ be a finitely
generated $\mcR$-module. Then there is a finite dimensional linear
subspace $D\subseteq\mcV$ for which $D\pitchfork\scS_d$.
\end{theorem}

\begin{proof}
First we give a proof in the special case $\mcV=\mcR E$ of a free
$\mcR$-module of finite rank. By Theorem
\ref{int_global_existence_of_transversal_for_R} there exists a
finite dimensional linear subspace $D'\subseteq\mcR$ such that
$D'\pitchfork\scS_d$. Clearly
 \[\hat{D}:=D'\otimes E\subseteq\mcR\otimes E=\mcR E\]
does the job as $\mfm_Y\mcR E=\mfm_Y\otimes E$. In the general
case let $\varphi:\mcR E\rightarrow\mcV$ be any surjective map of
$\mcR$-modules. As we just saw there is a $\hat{D}\subseteq\mcR E$
such that $\hat{D}\pitchfork\scS_d(\mcR E)$. We claim now that
$\varphi(\hat{D})\pitchfork\scS_d(\mcV)$. Hence let $Y\in\scS_d$.
Then
\[B=\varphi^{-1}(\mfm_Y\mcV)=\mfm_Y\mcR E+\ker\varphi\]
and so $\hat{D}+B=\mcR E$. Therefore
\[\varphi(\hat{D})+\mfm_Y\mcV=\mcV\]
and indeed $D=\varphi(\hat{D})\pitchfork\scS_d(\mcV)$.
\end{proof}

Let us fix a finite dimensional linear subspace $D\subseteq\mcV$.
Then there is the following adjoint pair of functors between
posets
\[\xymatrix@C+10pt{
\{\textrm{linear subspaces of $D$}\} \ar[r]<3pt>^-{F_D} &
\{\textrm{$\mcR$-submodules of $\mcV$}\} \ar[l]<3pt>^-{U_D}}\]
where $F_D$ maps $L\subseteq D$ to the $\mcR$-submodule $\mcR L$
generated by $L$ while $U_D$ maps $B\subseteq\mcV$ to $B\cap D$.

%

\begin{proposition} \label{int_injectivity_proposition}
If $D\pitchfork\scM_{d+1}$ then $F_DU_D|_{\scM_d}=\id$. In
particular $U_D|_{\scM_d}$ is injective and if $B'$ and $B''$ are
submodules of codimension at most $d$ then
\[B'\subseteq B''\quad\Leftrightarrow\quad D\cap B'\subseteq D\cap B''\]
\end{proposition}

\begin{proof}
Let $B\in\scM_d$. For brevity we denote
\[\hat{B}=\mcR(D\cap B)=F_DU_D(B)\]
and let us assume that $\hat{B}\subsetneqq B$. We choose a maximal
submodule $\tilde{B}$ with
\[\hat{B}\subseteq\tilde{B}\subsetneqq B\]
which is possible by the finite generation of $B$. By maximality
$B/\tilde{B}\cong\mcR/\mfm_y$ for some maximal ideal $\mfm_y$ and in
particular $\tilde{B}\in\scM_{d+1}$, hence $D\pitchfork\tilde{B}$.
By transversality then $D\cap\hat{B}\subseteq
D\cap\tilde{B}\subsetneqq D\cap B$. This is a contradiction with the
triangle identity $U_DF_DU_D=U_D$ for the adjoint pair $(F_D,U_D)$.

The second claim follows for in that case $B'=F_DU_D(B')$,
$B''=F_DU_D(B'')$ and both $F_D$ and $U_D$ preserve inclusions.
\end{proof}

\begin{remark}
Mere $D\pitchfork\scM_d$ would not be enough for injectivity of
$U_D$ as the example of polynomials $\bbR^m\ra\bbR$ of degree at
most $d-1$ shows.
\end{remark}

\begin{proposition} \label{int_naturality}
If $D\subseteq D'$ are finite dimensional linear subspaces of
$\mcV$ such that $D\pitchfork\scM_d$ then the map
\[U_{D'}\scM_d\longrightarrow U_{D}\scM_d\]
sending a subspace $L\subseteq D'$ to its intersection with $D$ is
continuous where both spaces are given the subspace topology from
the respective Grassmannian manifolds of subspaces of codimension
$d$. Also the map of the canonical $d$-dimensional bundles (with
fibres $D'/L$ over $L$ and corresponding fibre $D/(D\cap L)$ over
the image $D\cap L$) is continuous.
\end{proposition}

\begin{proof}
Let $V\subseteq G_d(D')$ be the subset of subspaces of $D'$ of
codimension $d$ transverse to $D$. Then the map $V\rightarrow
G_d(D)$ sending $L$ to $D\cap L$ is continuous (even smooth
) and the map from the statement is just its restriction to the
subset $U_{D'}\scM_d$. The same works for the canonical bundles.
\end{proof}


%

\section{The topology on $\scM_d$}

It turns out that one of the most important tools in the proofs in
this section is an existence of a continuous interpolation map. It
is a map $A:\mcR\times\scS_d\ra D$ (for some $D\pitchfork\scS_d$)
with the property $A(f,Y)=f$ modulo $\mfm_Y$. By definition for
each $f\in\mcR$ and $Y\in\scS_d$ such an $A(f,Y)$ exists but there
is no obvious choice and in particular it is not obvious that
there is a continuous way of choosing it. For example if $D$ is
the space of polynomials of degree at most $d-1$, we constructed
an interpolation in the proof of Lemma
\ref{int_local_existence_of_transversal_for_R}. However as the
following example shows it fails to be continuous.

\begin{example}
Let $M=\bbR^2$, $Y_n=[(0,0),(-1/n,1/n),(1/n,1/n)]$. Let us
consider the function $f(x,y)=y$ and let us construct the
interpolations in the space $D$ of truncated polynomials of degree
at most $2$ by the method from the proof of Lemma
\ref{int_local_existence_of_transversal_for_R} with respect to
$\alpha(x,y)=x$ which clearly is injective on each $Y_n$. Easily
 \[p_n(x,y)=0\cdot (x+1/n)(x-1/n)+n/2\cdot x(x-1/n)+n/2\cdot x(x+1/n)=nx^2\]
and therefore $p_n$ does not converge in $D$. In fact there is no
choice of $\alpha$ which would produce a convergent sequence. \qed
\end{example}

In the last example it is very easy to produce a convergent
interpolation sequence (after all $f\in D$, so we can as well take
$p_n=f$). The following construction is one way how to construct
canonical interpolations in the local case $M=\bbR^m$.

Let $f:\bbR^m\rightarrow E$ be a smooth function. If
$(x_0,\ldots,x_r)\in (\bbR^m)^{r+1}$ we denote by
$[x_0,\ldots,x_r]:\Delta^r\rightarrow\bbR^m$ the unique affine map
sending the vertices of the standard $r$-simplex $\Delta^r$ to
$x_0,\ldots,x_r$. It is obvious that this gives a bijective
correspondence between $(\bbR^m)^{r+1}$ and affine maps
$\Delta^r\rightarrow\bbR^m$. We will denote a general affine map
$\Delta^r\rightarrow\bbR^m$ by $\sigma$, if we do not want to
emphasize the values at vertices. By embedding it linearly into
$\bbR^r$ we give $\Delta^r$ the Lebesgue measure in which the volume
is $1$. Then we are able to define unambiguously
\[\mcI(f,\sigma)\in\Hom(S^r\bbR^m,E)\]
a symmetric $r$-form on $\bbR^m$ with values in $E$ to be
\[\mcI(f,\sigma)=\int_{\Delta^r}f^{(r)}\sigma\]
where $f^{(r)}:\bbR^m\ra\Hom(S^r\bbR^m,E)$ denotes the $r$-fold
derivative of $f$. In an obvious way  by omitting $x_i$ we get
$\partial_i\sigma:\Delta^{r-1}\rightarrow\bbR^m$ and thus forms
\[\mcI(f,\partial_i\sigma)\in\Hom(S^{r-1}\bbR^m,E)\]

\begin{lemma}\label{int_integration_formula}
For any smooth function $g:\bbR^m\rightarrow E$ and for any
$\sigma=[x_0,\ldots,x_r]$, the following holds for $0\leq i,j\leq
r$
\[\int_{\Delta^r}g'_{x_j-x_i}\sigma=
r\left(\int_{\Delta^{r-1}}g(\partial_i\sigma)-
\int_{\Delta^{r-1}}g(\partial_j\sigma)\right)\] where $g'_{x_j-x_i}$
denotes the derivative of $g$ in the direction $x_j-x_i$. In
particular by taking $g=f^{(r-1)}$ we have
\[\mcI(f,\sigma)(v_1,\ldots,v_{r-1},x_j-x_i)=
r(\mcI(f,\partial_i\sigma)-
\mcI(f,\partial_j\sigma))(v_1,\ldots,v_{r-1})\]
\end{lemma}

\begin{proof}
Define a map $\delta:\Delta^{r-2}\rightarrow\Delta^r$ to be
\[[e_0,\ldots,\hat{e}_i,\ldots,\hat{e}_j,\ldots,e_r]\]
where $e_n$ are the vertices of $\Delta^r$. We think of
$\Delta^{r-1}$ as a convex span of $\Delta^{r-2}$ (with vertices
$e_1,\ldots,e_{r-1}$) and an additional point $e_0$. Then we can
define a homotopy
\[h:\Delta^{r-1}\times I\rightarrow\Delta^r\]
by a formula (with $x$ running over $\Delta^{r-2}$)
\[h(t_0e_0+(1-t_0)x,t)=t_0((1-t)e_i+te_j)+(1-t_0)\delta(x)\]
One sees easily that $h(-,0)$ is the inclusion of the $j$-th face of
$\Delta^r$ and $h(-,1)$ the inclusion of the $i$-th one. To compute
the determinant of $h'$ we choose a basis
$(e_1-e_0,\ldots,e_{r-1}-e_0,e)$ of $T(\Delta^{r-1}\times I)$ where
$e$ is the unit tangent vector of $I$. Then
\[h'(t_0e_0+(1-t_0)x,t)(e_n-e_0)=\delta(e_n)-e_i-t(e_j-e_i)\]
and
\[h'(t_0e_0+(1-t_0)x,t)(e)=t_0(e_j-e_i)\]
Hence we easily get a formula
\[|\det h'(t_0e_0+(1-t_0)x,t)|=ct_0\]
for some constant $c$ and it is not difficult to see that $c=r$.
Then
\[\int_{\Delta^r}g'_{x_j-x_i}\sigma=
\int_{\Delta^{r-1}\times I}rt_0g'_{x_j-x_i}\sigma h=
r\int_{\Delta^{r-1}}\int_0^1t_0g'_{x_j-x_i}\sigma h(-,t)\
\textrm{d}t\] Now note that
\[\frac{\partial}{\partial t}(g\sigma h)=t_0g'_{x_j-x_i}\sigma h\]
and so
\begin{eqnarray*}
\int_{\Delta^{r-1}}\int_0^1t_0g'_{x_j-x_i}\sigma h(-,t)
\textrm{d}t & = & \int_{\Delta^{r-1}}g\sigma
h(-,1)-\int_{\Delta^{r-1}}g\sigma h(-,0)\\
& = & \int_{\Delta^{r-1}}g(\partial_i\sigma)-
\int_{\Delta^{r-1}}g(\partial_j\sigma)
\end{eqnarray*}
\end{proof}

\begin{corollary} \label{int_taylor_expansion}
The following formula holds
 \begin{equation} \label{int_taylor_expansion_expression}
  \begin{split}
   f(x_r)=\mcI(f,[x_0])
   +\cdots & +\frac{1}{i!}\cdot\mcI(f,[x_0,\ldots,x_i])(x_r-x_0,\ldots,x_r-x_{i-1})+\cdots
  \\
   & +\frac{1}{r!}\cdot\mcI(f,[x_0,\ldots,x_r])(x_r-x_0,\ldots,x_r-x_{r-1})
  \end{split}
 \end{equation}
\end{corollary}

\begin{proof}
We use induction with respect to $r$. According to the previous
lemma
 \begin{equation} \label{int_highest_term_in_taylor_expansion}
  \frac{1}{r!}\cdot\mcI(f,[x_0,\ldots,x_r])(x_r-x_0,\ldots,x_r-x_{r-1})
 \end{equation}
is equal to
 \[\frac{1}{(r-1)!}\cdot(\mcI(f,[x_0,\ldots,\widehat{x_{r-1}},x_r])-\mcI(f,[x_0,\ldots,x_{r-1}]))(x_r-x_0,\ldots,x_r-x_{r-2})\]
Adding the remaining terms of
(\ref{int_taylor_expansion_expression}) to
(\ref{int_highest_term_in_taylor_expansion}) and using the inductive
hypothesis on $[x_0,\ldots,\widehat{x_{r-1}},x_r]$ we prove the
inductive step.
\end{proof}

\begin{corollary} \label{int_equivalence}
The following conditions are equivalent
\begin{itemize}
\item[$(i)$]{$\mcI(f,\tau)=0$ for all the faces $\tau$ of $\sigma$.}
\item[$(ii)$]{$\mcI(f,\partial_1\cdots\partial_r\sigma)=\cdots=
\mcI(f,\partial_i\cdots\partial_r\sigma)=\cdots=\mcI(f,\partial_r\sigma)=\mcI(f,\sigma)=0$.}
\end{itemize}
\end{corollary}

\begin{proof}
We assume $(ii)$. By induction we can also assume that for all the
faces $\tau$ of $\partial_r\sigma$, we have $\mcI(f,\tau)=0$. By
the previous lemma $\mcI(f,\partial_i\sigma)=0$ for all $i$. As
there is a common face of $\partial_i\sigma$ and
$\partial_r\sigma$ we see that up to a renumbering of vertices the
condition $(ii)$ is satisfied for $\partial_i\sigma$ and by
induction again we get $(i)$ for all the faces of
$\partial_i\sigma$.
\end{proof}

\begin{corollary} \label{int_vanishing_of_integrals_implies_vanishing_of_jets}
Let $\sigma=[x_0,\ldots,x_r]$. If
\[\mcI(f,\partial_1\cdots\partial_r\sigma)=\cdots=
\mcI(f,\partial_i\cdots\partial_r\sigma)=\cdots=\mcI(f,\partial_r\sigma)=\mcI(f,\sigma)=0\]
and if $y$ appears $k$-times in $x_0,\ldots,x_r$ then
$j^{k-1}_yf=0$.
\qed
\end{corollary}

\begin{remark}
The converse is not true in general (with the exception
$x_0=\cdots=x_r$) for dimensional reasons unless $m=1$: if
$\{(y_1,k_1),\ldots,(y_n,k_n)\}$ denotes the class of
$(x_0,\ldots,x_r)$ in $\scS_{r+1}(\bbR)=\bbR^{r+1}/\Sigma_{r+1}$
then $j^{k_i-1}_{y_i}f=0$ for all $i=1,\ldots,n$ implies
$\mcI(f,[x_0,\ldots,x_j])=0$ for all $j=0,\ldots,r$.
\end{remark}

Now we will explain how this leads to an interpolation map. First we
restrict ourselves to interpolation at points close to a single
point, later generalizing to a number of points. This is only to
lighten the notation a bit. Going back from the local situation to
the case of a compact manifold $M$ we consider $\mcV=\mcR
E=C^\infty(M,E)$, a free $\mcR$-module of a finite rank. We fix
$y\in M$ and $k\geq 1$ and identify a neighborhood of $y$ with
$\bbR^m$. We also fix a complementary linear subspace $F$ to the
submodule $(\mfm_y)^k\mcV$ and define the following map
\[\begin{split}
G:\mcV\times(\bbR^m)^k\lra &
\Hom(S^0\bbR^m,E)\times\cdots\times\Hom(S^{k-1}\bbR^m,E)\times(\bbR^m)^k
\\
\cong\  & J^{k-1}_*(\bbR^m,E)\times(\bbR^m)^k
\end{split}\]
(with $J^{k-1}_*(\bbR^m,E)$ being any $J^{k-1}_x(\bbR^m,E)$ -- they
are all identified via translations) by the formula
\[G(f,(y_1,\ldots,y_k))=(\mcI(f,[y_1]),\ldots,\mcI(f,[y_1,\ldots,y_k]),(y_1,\ldots,y_k))\]
We denote by $G_F$ its restriction
\[G_F:F\times(\bbR^m)^k\rightarrow J^{k-1}_*(\bbR^m,E)\times(\bbR^m)^k\]
Note that for each $f\in\mcV$ the map $G(f,-)$ is continuous (in
fact smooth) and therefore so is $G_F$. Our transversality
condition on $F$ implies that on the fibres over $(y,\ldots,y)$
\[(G_F)_{(y,\ldots,y)}:F\rightarrow J^{k-1}_*(\bbR^m,E)\]
is a linear isomorphism. Hence we find a neighborhood of
$(y,\ldots,y)$ in $M^k$ of the form $U^k$ with $U$ compact convex,
such that the restriction
\[G_F:F\times U^k\rightarrow J^{k-1}_*(\bbR^m,E)\times U^k\]
is an isomorphism of vector bundles over $U^k$. Hence we can
define a map
\[\hat{A}:\mcV\times U^k\xlra{G}J^{k-1}_*(\bbR^m,E)\times U^k\xlra{G_F^{-1}}F\times U^k\rightarrow F\]
Now note that according to Corollary \ref{int_equivalence} the
value of $\hat{A}$ does not depend on the ordering of the points
and hence we get
\[A:\mcV\times U^k/\Sigma_k\rightarrow F\]
with the property that $A(f,Y)$ is an interpolation of $f$ at $Y$,
i.e. such that $f=A(f,Y)$ modulo $\mfm_Y\mcV$.

\begin{theorem}
The interpolation map $A$ is continuous and $A(f,\{(y,k)\})=0$
whenever $f\in(\mfm_y)^k\mcV$.
\end{theorem}

\begin{proof}
By the construction of $A$ it is enough to show the continuity of
each component
\[G_i:\mcV\times U^k\rightarrow\Hom(S^i\bbR^m,E)\]
of $G$, $i=0,\ldots,k-1$. Hence let us fix
$(f,(y_1,\ldots,y_k))\in\mcV\times U^k$ and denote
$Y=(y_1,\ldots,y_k)$ for a short. We choose a norm on
$\Hom(S^i\bbR^m,E)$ and a neighborhood
\[\{\alpha\in\Hom(S^i\bbR^m,E)\ |\ ||\alpha-G_i(f,Y)||<\varepsilon\}\]
Because of the continuity of $G_i(f,-)$ there is a neighborhood
$V$ of $Y$ on which
\[||G_i(f,-)-G_i(f,Y)||<\varepsilon/2\]
Hence if $g$ is such that $||g^{(i)}-f^{(i)}||<\varepsilon/2$ on
$U$ then also
\[||G_i(g,-)-G_i(f,-)||<\varepsilon/2\]
on $U$ and finally on $V$ we have
\[||G_i(g,-)-G_i(f,Y)||\leq||G_i(g,-)-G_i(f,-)||+||G_i(f,-)-G_i(f,Y)||<\varepsilon\]
Because the condition on $g^{(i)}$ defines a neighborhood of $f$ in
$\mcV$, this finishes the proof.
\end{proof}


\begin{example}
In the one dimensional local case $M=\bbR$ the space $P_{k-1}$ of
polynomials $\bbR\ra\bbR$ of degree at most $k-1$ is clearly a
complementary subspace of $(\mfm_y)^k$. For any
$(y_1,\ldots,y_k)\in\bbR^k$ consider the following basis of
$P_{k-1}$:
 \[\{1,(x-y_1),\ldots,1/i!\cdot(x-y_1)\cdots(x-y_i),\ldots,1/k!\cdot(x-y_1)\cdots(x-y_k)\}\]
Then one has
 \[\mcI(1/i!\cdot(x-y_1)\cdots(x-y_i),[y_1,\ldots,y_j])=\delta_{ij}\]
This is clear for $j\geq i$ by a direct computation and for $i>j$
this follows from the remark after Corollary
\ref{int_vanishing_of_integrals_implies_vanishing_of_jets}.
Consequently one obtains a formula
 \[A(f,\Sigma_k(y_1,\ldots,y_k))=\mcI(f,[y_1])+\cdots+1/(k-1)!\cdot\mcI(f,[y_1,\ldots,y_k])(x-y_1)\ldots(x-y_{k-1})\]
This is the so-called Lagrange's interpolation formula. It can be
found for example in Section~4.2 of \cite{Schatzman}.
\end{example}

Now if we have an arbitrary
$Y=\{(y_1,k_1,),\ldots,(y_n,k_n)\}\in\scS_d$ we identify a
neighbourhood of each $y_i$ with $\bbR^m$. Writing
\[\bbR^{md}\cong(\bbR^m)^{k_1}\times\cdots\times(\bbR^m)^{k_n}\]
we replace $G$ by the corresponding map
\[G:\mcV\times\bbR^{md}\rightarrow
J^{k_1-1}_*(\bbR^m,E)\times\cdots\times
J^{k_n-1}_*(\bbR^m,E)\times\bbR^{md}\] Again if $F$ is a
complementary linear subspace to $\mfm_Y\mcV$ then on the fibres
\[(G_F)_{((y_1,\ldots,y_1),\ldots,(y_n,\ldots,y_n))}:F\rightarrow
J^{k_1-1}_*(\bbR^m,E)\times\cdots\times J^{k_n-1}_*(\bbR^m,E)\] is
an isomorphism and we get a neighborhood of the form
$U_1^{k_1}\times\cdots\times U_n^{k_n}$ over which $G_F$ is an
isomorphism. Denoting the induced neighborhood of $Y$ in $\scS_d$ by
$W$ one gets an interpolation map
\[A:\mcV\times W\rightarrow F\]

\begin{theorem}
The interpolation map $A$ is continuous and $A(f,Y)=0$ whenever
$f\in \mfm_Y\mcV$. \qed
\end{theorem}

Next corollary is crucial for the proof of compactness of $\scM_d$.

\begin{corollary} \label{int_convergence_in_Sd}
Let $M$ be a compact manifold, let $\mcV=\mcR E$, let
$D\pitchfork\scS_d$. Let $Y_p\in\scS_d$ ($p=1,2,\ldots$) be a
sequence converging to $Y$. Then every $v\in D\cap \mfm_Y\mcV$ is
a limit of some sequence $v_p\in D\cap \mfm_{Y_p}\mcV$.
\end{corollary}

\begin{proof}
One chooses a complementary linear subspace $F\subseteq D$ to
$\mfm_Y\mcV$ and gets a sequence $v_p=v-A(v,Y_p)\rightarrow v$.
\end{proof}

\begin{theorem} \label{int_existence_of_continuous_interpolation}
Let $M$ be a compact manifold, let $\mcV=\mcR E$, let
$D\pitchfork\scS_d$. Then there is a continuous fibrewise linear
interpolation map
\[A:\mcV\times\scS_d\rightarrow D\]
i.e. a map such that $A(f,Y)=f$ modulo $\mfm_Y\mcV$.
\end{theorem}

\begin{proof}
One glues the local interpolation maps using a partition of unity
on $\scS_d$. Here one uses the fact that the interpolation maps
form an affine space.
\end{proof}

With Corollary \ref{int_convergence_in_Sd} at hand we can prove
the main theorem in the same way Glaeser did in \cite{Glaeser} for
a parallelepiped in $\bbR^m$.

\begin{theorem} \label{int_compactness_for_RE}
Let $M$ be a compact manifold, let $\mcV=\mcR E$ be a free
$\mcR$-module of finite rank and let $D\pitchfork\scS_{d+1}$. Then
$U_D\scM_d\subseteq G_d(D)$ is a closed subset hence compact. Also
the map $\pi:U_D\scM_d\rightarrow\scS_d$ sending $D\cap B$ to the
weighted spectrum of $B$ is continuous.
\end{theorem}

\begin{proof}
Let $B_p$ be a sequence of submodules of codimension $d$ and let
$L_p=D\cap B_p$. Let us assume that $L_p$ converges to $L$ in
$G_d(D)$. We will construct a submodule $B$ such that $L=D\cap B$.
By taking a subsequence we can assume that the sequence of
weighted spectra corresponding to $B_p$ converges to $Y$. We set
$B=\mfm_Y\mcV+L$. According to Corollary
\ref{int_convergence_in_Sd}, $D\cap \mfm_Y\mcV\subseteq L$ and
this easily implies that $D\cap B=L$ and that the codimension of
$B$ is $d$. It remains to show that $B$ is indeed a submodule.

We choose a finite dimensional linear subspace $P\subseteq\mcR$
complementary to $\mfm_Y$. We can make sure that $1\in P$. Then
$\mfm_Y+P=\mcR$ and to prove that $B$ is a submodule one needs to
prove the inclusion labeled by $?$ in
\[\mcR B=(\mfm_Y+P)(\mfm_Y\mcV+L)\subseteq \mfm_Y\mcV+PL\overset{?}{\subseteq}\mfm_Y\mcV+L=B\]
In order to do so one constructs an analogous sequence $PD\cap
B_p$ in $PD$ and by taking a further subsequence one can assume
that this sequence converges to $L'$. As $1\in P$ clearly
$L\subseteq L'$ and so
\[\mfm_Y\mcV+L\subseteq \mfm_Y\mcV+L'\]
But the codimensions are the same (equal to $d$) and therefore
these two subspaces must be equal. As $P(D\cap B_p)\subseteq
PD\cap B_p$ taking the limit we get $PL\subseteq L'$ and so
\[\mfm_Y\mcV+PL\subseteq \mfm_Y\mcV+L'=\mfm_Y\mcV+L\]
finishing the proof that $B$ is a submodule.

To prove continuity of the map $U_D\scM_d\rightarrow\scS_d$ it is
enough (by compactness of $\scS_d$) to show that the weighted
spectrum of the submodule $B$ that we just constructed is indeed
$Y$. The idea is to use the primary decomposition and compare the
dimensions. We choose some disjoint closed neighborhoods $C_i$ of
$y_i$ and decompose
\[B_p=(B_p+\mfn_{C_1}\mcV)\cap\cdots\cap(B_p+\mfn_{C_n}\mcV)\]
By the convergence $\spec(B_p)\rightarrow Y$ we easily see that
the codimension of each $B_p^i:=B_p+\mfn_{C_i}\mcV$ is $k_i$, at
least for all big enough $p$. Assuming that each sequence $D\cap
B_p^i$ converges in $G_{k_i}(D)$ to $D\cap B^i$ we certainly have
 {\setlength\arraycolsep{2pt}
 \begin{eqnarray*}
   D\cap B & = & \lim(D\cap B_p^1)\cap\cdots\cap(D\cap B_p^n)\subseteq
  \\
   & \subseteq & (D\cap B^1)\cap\cdots\cap(D\cap B^n)= D\cap(B^1\cap\cdots\cap B^n)
 \end{eqnarray*}}
and by the transversality assumptions and Proposition
\ref{int_injectivity_proposition} we conclude that $B\subseteq
B^1\cap\cdots\cap B^n$. By the proof above the spectrum of each
$B^i$ consists just of $y_i$ hence the $B^i$ are in general
position and therefore the codimension of $B^1\cap\cdots\cap B^n$
is also $d$. Thus $B=B^1\cap\cdots\cap B^n$ is the primary
decomposition of $B$ proving the claim about its weighted
spectrum.
\end{proof}

Hence by Proposition \ref{int_naturality} we know that the
topology on $\scM_d$ does not depend on the choice of $D$ as long
as $D\pitchfork\scM_d$ and $U_D$ is injective.

\begin{remark}
Let us consider the case $\mcV=\mcR$ and for simplicity denote
$\scM_d(M)=\scM_d(\mcR)$. For an open subset $U\subseteq M$ we
define $\scM_d(U)\subseteq\scM_d(M)$ by the pullback square
 \[\xymatrix{
  \scM_d(U) \ar@{c->}[r] \ar[d] & \scM_d(M) \ar[d]^{\pi}
 \\
  \scS_d(U) \ar@{c->}[r] & \scS_d(M)
 }\]
In other words $\scM_d(U)$ is the space of ideals in $\mcR$ with
spectrum in $U$. Theorem \ref{int_compactness_for_RE} implies that
the topology of $\scM_d(U)$ is independent of $M$ and $\scM_d(M)$ is
a union of $\scM_d(U)$ as $U$ varies over coordinate charts. In
other words $\scM_d(M)$ is local. On the other hand the glueing maps
are not polynomial and so \cite{Glaeser} could not be used.
\end{remark}

Let $E_d(D)$ denote the canonical $d$-dimensional vector bundle
over $G_d(D)$, whose fibre over $F\subseteq D$ is $D/F$. Then
again by the same proposition the restriction
$E_d(D)|_{U_D\scM_d}$ is independent of $D$ and will be denoted by
$\scE_d$.

\begin{corollary}
Let $M$ be a compact manifold and let $\mcV=\mcR E$ be a free
$\mcR$-module of finite rank. Then the natural map
\[\mcV\times\scM_d\rightarrow\scE_d\]
defined on the fibre over $B$ by
\[\mcV\rightarrow\mcV/B\cong D/(D\cap B)\]
is a continuous quotient map (in the topological sense) of vector
bundles over $\scM_d$.
\end{corollary}

\begin{proof}
The map in question can be defined alternatively as
\[\mcV\times\scM_d\xlra{\id_\mcV\times(\pi,\id_{\scM_d})}\mcV\times\scS_d\times\scM_d
\xlra{A\times\id_{\scM_d}}D\times\scM_d\rightarrow\scE_d\] using
the interpolation map $A$ constructed above and hence is by
Theorem \ref{int_existence_of_continuous_interpolation}
continuous. Using an inner product on $D$ one can easily find a
section proving that it is a quotient map.
\end{proof}

\begin{corollary}
\label{int_submodules_containing_fixed_submodule}
Let $M$ be a compact manifold, let $\mcV=\mcR E$ be a free
$\mcR$-module of finite rank and let $K\subseteq\mcR E$ be a
submodule. Then the subset
 \[\{B\in\scM_d\ |\ K\subseteq B\}\subseteq\scM_d\]
is closed hence compact.
\end{corollary}

\begin{proof}
For each $v\in K$ the restriction of the canonical map
 \[\scM_d\cong\{v\}\times\scM_d\subseteq\mcV\times\scM_d\ra\scE_d\]
is continuous by the last corollary and the set $\{B\in\scM_d\ |\
v\in B\}$ is precisely the preimage of the zero section and
therefore closed. So is then their intersection over all $v\in K$,
the set from the statement.
\end{proof}

%

\begin{theorem}
Let $M$ be a compact manifold, let $\mcV$ be any finitely generated
$\mcR$-module and let $D\pitchfork\scM_d$. Then $U_D\scM_d\subseteq
G_d(D)$ is a closed subset hence compact. If $\mcV$ is topologically
finitely generated then the natural map
\[\mcV\times\scM_d\rightarrow\scE_d\]
is a continuous quotient map.
\end{theorem}

\begin{proof}
Let $\varphi:\mcR E\rightarrow\mcV$ be any surjective map of
$\mcR$-modules with $E$ a finite dimensional real vector space,
let $K=\ker\varphi$. As we can replace $D$ by any bigger subspace
to prove compactness we can assume that $D=\varphi(\hat{D})$ for
some $\hat{D}\pitchfork\scM_{d+1}(\mcR E)$. The submodules of
$\mcV\cong\mcR E/K$ are identified with those submodules of $\mcR
E$ that contain $K$. Now we reconstruct this relation inside
$\hat{D}$ and $D$: setting $\hat{K}:=\hat{D}\cap K$ we see that
$\varphi|_{\hat{D}}:\hat{D}\ra D$ can be identified with the
projection $\hat{D}\ra\hat{D}/\hat{K}$. Consequently the map
$L\mapsto\varphi(L)$ is continuous when restricted to subspaces
containing $\hat{K}$ and so the dashed arrow in the diagram
 \[\entrymodifiers={!!<0pt,\the\fontdimen22\textfont2>+}
  \xymatrixnocompile{
   U_{\hat{D}}\{B\ |\ K\subseteq B\} \ar@{c->}[r] \ar@{-->}[d] &
   \{L\ |\ \hat{K}\subseteq L\} \ar@{c->}[r] \ar[d] &
   G_d\hat{D}
  \\
   U_D\scM_d(\mcV) \ar@{c->}[r] &
   G_dD
 }\]
is then a continuous bijection. As $U_{\hat{D}}\{B\ |\ K\subseteq
B\}$ is compact by Corollary
\ref{int_submodules_containing_fixed_submodule} this dashed arrow
is a homeomorphism identifying $U_D\scM_d(\mcV)$ with a subspace
of $U_{\hat{D}}\scM_d(\mcR E)$. In the same way $\scE_d(\mcV)$ can
be thought of as the restriction of $\scE_d(\mcR E)$ to
$\scM_d(\mcV)$ and in the diagram
 \[\xymatrix{
   {}\mcV\times\scM_d(\mcV) \ar@{-->}[rd] &
   {}\mcR E\times\scM_d(\mcV) \ar@{c->}[r] \ar@{->>}[l] \ar[d] &
   {}\mcR E\times\scM_d(\mcR E) \ar[d]
  \\
   &
   {}\scE_d(\mcV) \ar@{c->}[r] &
   {}\scE_d(\mcR E)
 }\]
the dashed arrow exists and is continuous by the properties of
quotients.
\end{proof}

We consider the space $G^r(\mcV)$ of $r$-dimensional linear
subspaces of $\mcV$ and give it the quotient topology from
$V^r(\mcV)/GL_r(\bbR)$ where
\[V^r(\mcV)\subseteq\mcV\times\cdots\times\mcV\]
is the subset of linearly independent $r$-tuples of vectors in
$\mcV$ and the action of $GL_r(\bbR)$ is the usual one.

\begin{corollary}
The set of all $r$-dimensional linear subspaces $D$ of $\mcV$ for
which $D\pitchfork\scM_d$ is open in $G^r(\mcV)$.
\end{corollary}

\begin{proof}
As the canonical projection $V^r(\mcV)\rightarrow G^r(\mcV)$ is
(defined to be) a quotient map we only need to show that the
preimage $S\subseteq V^r(\mcV)$ of the set from the statement is
open. We consider the map
\[\alpha:\mcV\times\cdots\times\mcV\times\scM_d\lra
\scE_d\times_{\scM_d}\cdots\times_{\scM_d}\scE_d
\cong\Hom(\bbR^r,\scE_d)\] sending $((v_1,\ldots,v_r),B)$ to the
homomorphism $\bbR^r\rightarrow\mcV/B$ which is determined by
sending the basis vectors to the classes of the $v_i$'s in $\mcV/B$.
According to the previous theorem $\alpha$ is continuous. The subset
$\textrm{SurHom}(\bbR^r,\scE_d)$ of surjective homomorphisms is open
in the target and hence so is
\[T:=(V^r(\mcV)\times\scM_d)\cap\alpha^{-1}(\SurHom(\bbR^r,\scE_d))\]
We can then describe $S$ as the set of those $x\in V^r(\mcV)$ for
which $\{x\}\times\scM_d\subseteq T$ and consequently $S$ is also
open by compactness of $\scM_d$.
\end{proof}

\section{Examples and the structure of $\scM_d$}

We will now give few examples to explain what the spaces $\scM_d$
look like. We specialize to the case $\mcV=\mcR=C^\infty(M,\bbR)$
and write $\scM_d(M):=\scM_d(C^\infty(M,\bbR))$. First we prove a
useful lemma which makes defining smooth maps (and in particular
immersions) into $\scM_d(M)$ a bit easier.

\begin{lemma} \label{int_smooth_maps_into_Md}
Let $N$ be a smooth manifold, let $D\pitchfork\scM_d(M)$ and let
 \[\varphi:\mcR\times N\ra\bbR^d\]
be a map such that for each $f\in\mcR$ the partial map
$\varphi(f,-):N\ra\bbR^d$ is smooth and such that for each $x\in
N$ the partial map $\varphi(-,x):\mcR\ra\bbR^d$ is surjective
linear whose kernel $I_x=\ker\varphi(-,x)$ is an ideal in $\mcR$.
Then the map
 \[\psi:N\ra U_D\scM_d(M)\subseteq G_d(D)\]
sending $x\in N$ to $U_D(I_x)$ is smooth. If moreover
$D\pitchfork\scM_{d+1}(M)$ then $X\in T_xN$ lies in the kernel of
the differential $\psi_*:TN\ra TG_d(D)$ if and only if
 \begin{equation} \label{int_condition_on_ideals}
  I_x\subseteq\{f\in\mcR\ |\ d(\varphi(f,-))(X)=0\}
 \end{equation}
\end{lemma}

\begin{proof}
Let $x\in N$ and denote $L=\psi(x)=D\cap\ker\varphi(-,x)$. We
choose a complementary subspace $F$ to $L$ inside $D$ and get a
chart on $G_d(D)$
 \[\Hom(L,F)\lra G_d(D)\]
given by sending a map $\alpha$ to its graph inside $L\times
F\cong D$. Let us consider the restriction $\varphi_D$ of
$\varphi$ to $D\times N$ and write it in the form
 \[\varphi_D:L\times F\times N\lra\bbR^d\]
Observe that the differential $d\varphi_D|_F$ is an isomorphism of
$F$ on $\bbR^d$ near $L\times F\times\{x\}$. In particular there
is a unique solution to the equation
 \[\varphi_D(v,\alpha(y)(v),y)=0\]
and it is automatically smooth. Clearly $\alpha(y):L\ra F$ is the
expression of $\psi(y)$ in the above coordinate chart (with
$\alpha(x)=0$). Moreover we have a formula for the derivative
 \[d\alpha(X)(v)=(d(\varphi_D(v,-,x)))^{-1}(d(\varphi_D(v,0,-)))(X)\]
In particular $X\in\ker\psi_*$ if and only if for each $v\in L$ it
lies in the kernel of $d(\varphi_D(v,0,-))$. To explain this
condition we introduce
 \[I_X:=\{f\in\mcR\ |\ \varphi(f,x)=0,\ d(\varphi(f,-))(X)=0\}\]
As the name suggests it is an ideal and to prove this one observes
that for each $y\in N$ we have a multiplication on $\bbR^d$
arising from the identification $\mcR/I_y\cong\bbR^d$. This family
is smooth in the sense of the map
 \[\mu:N\lra\Hom(\bbR^d\otimes\bbR^d,\bbR^d)\]
being smooth. If we temporarily denote
 \[f(x):=\varphi(f,x)\quad\textrm{and}\quad
 df(X):=d(\varphi(f,-))(X)\]
then for $f,g\in\mcR$ we get
 \[d(fg)(X)=\mu(x)\bigl(f(x)\otimes dg(X)+df(X)\otimes g(x)\bigr)+d\mu(X)(f(x)\otimes g(x))\]
Therefore if one of $f$, $g$ lies in $I_X$ then so does their
product.

The condition (\ref{int_condition_on_ideals}) from the statement
is then equivalent to $I_X=I_x$. Assuming that this equality
holds, every $v\in L\subseteq I_x$ lies in $I_X$ implying that
$d(\varphi_D(v,0,-))(X)=0$. Therefore in this case $\psi_*(X)=0$.
If on the other hand $I_X\subsetneqq I_x$ then there is an ideal
$J$ which is maximal among those for which $I_X\subseteq
J\subsetneqq I_x$. Necessarily $J\in\scM_{d+1}(M)$ and by our
assumption $D\cap J\subsetneqq D\cap I_x=L$ so that there is $v\in
L$ for which $v\not\in I_X$ implying that
$d(\varphi_D(v,0,-))(X)\neq0$ and $\psi_*(X)\neq0$.
\end{proof}

\begin{remark}
The first part of the proof applies even when $\varphi(f,-)$ is merely continuous proving that $\psi$ is also continuous in this situation.
\end{remark}


The space $\scM_d(M)$ contains as a subspace the configuration
space
 \[M^{[d]}=M^{(d)}/\Sigma_d\subseteq\scS_d\]
We will show now that it is in general an embedded submanifold.
I do not know 
whether $\scM_d(M)$ is itself a manifold with corners (or a
manifold stratified subset) with the configuration space $M^{[d]}$
being its interior. It is certainly an interesting question to
pursue.

\begin{proposition} \label{int_configuration_space_inclusion}
Let $D$ be a finite dimensional linear subspace of $\mcR$ such
that $D\pitchfork\scM_{d+1}(M)$. Then the inclusion
 {\setlength\arraycolsep{2pt}
 \begin{eqnarray*}
  \psi:M^{[d]}\subseteq\scS_d & \ra & \scM_d(M)\subseteq G_d(D)
 \\
  Y & \mapsto & \mfm_Y
 \end{eqnarray*}}
is a smooth embedding.
\end{proposition}

\begin{proof}
As the map in question has a continuous inverse (restriction of
$\pi$ from Theorem \ref{int_compactness_for_RE}) we only need to
show that it is an immersion. First we express $\psi$ locally via
a map $\varphi:\mcR\times M^{[d]}\ra\bbR^d$ as in Lemma
\ref{int_smooth_maps_into_Md} and compute the kernel of $\psi_*$
using the same lemma. Therefore let $(x_1,\ldots,x_d)\in M^{(d)}$
and identify a neighbourhood of $[(x_1,\ldots,x_d)]\in M^{[d]}$
with a product $U_1\times\cdots\times U_d$ of disjoint
neighbourhoods $U_i$ of $x_i$. Then we can define
$\varphi:\mcR\times U_1\times\cdots\times U_d\ra\bbR^d$ by
 \[(f,y_1,\ldots,y_d)\mapsto(f(y_1),\ldots,f(y_d))\]
Clearly all the assumptions of Lemma \ref{int_smooth_maps_into_Md}
are satisfied and so $\psi$ is a smooth map. Also for
 \[(X_1,\ldots,X_d)\in T_{x_1}U_1\times\cdots\times T_{x_d}U_d\]
we have $d(\varphi(f,-))(X_1,\ldots,X_d)=(df(X_1),\ldots,df(X_d))$
and this can be zero on $\mfm_{\{x_1,\ldots,x_d\}}$ only if
$X_1=\cdots=X_d=0$.
\end{proof}

Now we investigate a (rather trivial) example where Proposition
\ref{int_configuration_space_inclusion} completely describes the
topology of $\scM_d(M)$.

\begin{example}[a description of $\scM_1(M)$]
As ideals of codimension $1$ are exactly the maximal ones the
canonical map $M\ra\scM_1(M)$ is a bijection. As it was shown to
be a smooth embedding in Proposition
\ref{int_configuration_space_inclusion} it is a diffeomorphism.
\end{example}

Another easy example is that of the space of ideals for a
one-dimensional manifold.

\begin{example}[one-dimensional manifolds]
It is clear that the only ideals in $J^r_0(\bbR,\bbR)$ are the
powers of the maximal ideal $\mfm_0$. In other words there is
exactly one ideal for each codimension and therefore the canonical
map
 \[\scM_d(M)\lra\scS_d(M)\]
is bijective. As it is also a continuous map between compact
Hausdorff spaces it is even a homeomorphism. In particular
 \[\scM_d(M)\cong\scS_d(M)=M^d/\Sigma_d\]
We will identify this space explicitly for $M=I$, the closed
interval. Namely, we claim that $\scM_d(I)$ is the $d$-dimensional
simplex $\Delta^d$. This is easily seen from the classical
subdivision of $I^d$ into simplices
 \[e_\sigma=\{(t_1,\ldots,t_d)\ |\ 0\leq t_{\sigma(1)}\leq\cdots\leq t_{\sigma(d)}\leq 1\}\]
indexed by permutations $\sigma\in\Sigma_d$. Note that the action
of $\Sigma_d$ just permutes them $\tau:e_\sigma\ra e_{\tau\sigma}$
(if we think of it as the left action) and therefore the
composition
 \[\Delta^d\xlra{e_\id}I^d\lra I^d/\Sigma_d\]
is a continuous bijection between compact Hausdorff spaces hence a
homeomorphism\footnote{Using Lemma \ref{int_smooth_maps_into_Md} on
the map $\varphi:\mcR\times\Delta^d\ra\bbR^d$ defined by
 \[(f,(t_1,\ldots,t_d))\mapsto(\mcI(f,[t_1]),\ldots,\mcI(f,[t_1,\ldots,t_d]))\]
one can even produce a diffeomorphism $\Delta^d\ra\scM_d(I)$.}.

As $S^1$ is a quotient of $I$ one can conclude that also
$\scS_d(S^1)$ is a quotient of $\scS_d(I)$. More precisely there is
a map $f:\partial_0\Delta^d\ra\partial_d\Delta^d$ from the $0$-th
face of $\Delta^d$ to its $d$-th face (an affine map sending the
$i$-th vertex $e_i$ of $\Delta^d$ to $e_{i-1}$) such that
$\scS_d(S^1)\cong\Delta^d/f$, the space obtained from $\Delta^d$ by
identifying the two faces via $f$ (e.g.~$\scS_2(S^1)$ is the
M\"obius strip). Also note that locally $\scM_d(S^1)$ looks like
$\scM_d(I)$ and so it is a manifold with corners. \qed
\end{example}

In order to describe $\scM_2(M)$ for a general manifold $M$ we
first consider the following construction, a slight modification
of a blow-up construction (cf. exercises 7-8 after Chapter 12 in
\cite{BrockerJanich}):

\begin{construction} \label{int_construction_of_blow_up}
Let $E\ra N$ be a smooth vector bundle and let us choose an inner
product on $E$. The construction below does not depend on this
choice and in fact can be given in an inner product free way.
Consider the unit sphere bundle $SE\ra N$ and a trivial ray bundle
$\eta E=\bbR_+\times SE\ra SE$ over $SE$. We have a canonical map
$\eta E\ra E$ given by sending $(t,v)$ to $tv$. Clearly it is a
diffeomorphism away from the zero sections: $\eta
E-0\xra{\cong}E-0$. One can also think of $\eta E$ as $E$ with an
open disk subbundle removed but then there is no preferred way of
defining a projection map $\eta E\ra E$ (the inclusion is not what
we are after).

Let $M$ be a manifold and $N$ a closed neat submanifold. We
consider the normal bundle $\nu$ of $N$ and give $M_N:=(M-N)\sqcup
S\nu$ a structure of a manifold (with corners in general) defined
in terms of an embedding $\iota:\nu\hookrightarrow M$ as
 \[(M-N)\cup_{\nu-0}\eta\nu\]
glued along $\nu-0$ via the embeddings $\iota:\nu-0\hookrightarrow
M-N$ and $\nu-0\cong\eta\nu-0\hookrightarrow\eta\nu$. It can be
shown that the resulting manifold does not depend on the choice of
the embedding $\iota$.

Now let us specialize to a manifold $M$ equipped with an
involution $\tau:M\ra M$ and denote the fixed point submanifold by
$N=M^\tau$. As $M-N$ is clearly dense in $M_N$ there may be only
one possible smooth extension of the involution $\tau:M\ra M$ to
$M_N$. We show that it exists and that on $S\nu$ it is just
multiplication by $-1$.

To do so we choose a $\tau$-invariant Riemannian metric on $M$.
The action of $\tau$ on $TM|_N$ is by an orthogonal involution and
therefore this bundle decomposes into a direct sum
 \[TM|_N=TN\oplus\nu\]
with $TN$ being the $(+1)$-eigenspace and $\nu$ the
$(-1)$-eigenspace. The canonical embedding
$\iota:\nu\hookrightarrow M$ given by
 \[(x,v)\mapsto\exp_x(v)\]
transforms the involution $\tau$ to multiplication by $-1$ (as
$\tau$ is an isometry and therefore preserves geodesics) and this
can be clearly extended to $\eta\nu$. As the extended involution
$\hat{\tau}:M_N\ra M_N$ has no fixed points
$M_{\hat{\tau}}:=M_N/\hat{\tau}$ is again a manifold with corners.
As a set it is a disjoint union of $(M-M^\tau)/\tau$ with the
projective bundle $S\nu/\{\pm 1\}$ of the fixed point submanifold
$M^\tau$ in $M$. In other words one could say that it is a ``blow
up'' of $M^\tau$ in $M/\tau$ (but as $M/\tau$ is not a manifold
this differs from the classical construction and in particular one
introduces a new boundary as was mentioned).
\end{construction}

\begin{example}[a description of $\scM_2(M)$]
Another easy example is that of codimension $2$ ideals. We claim
that for the involution $\tau:M^2\ra M^2$ switching the two
coordinates, the above constructed $(M^2)_{\hat{\tau}}$ is
homeomorphic to $\scM_2(M)$.

As $(M^2)_{\hat{\tau}}$ is easily seen to be compact, we only need
to produce a continuous map $(M^2)_{\hat{\tau}}\ra\scM_2(M)$ and
show that it is bijective. By the definition of
$(M^2)_{\hat{\tau}}$, such a map can be obtained from two
$\hat{\tau}$-invariant maps
 \[M^2-\Delta\lra\scM_2(M)\qquad\mathrm{and}\qquad\eta\nu\lra\scM_2(M)\]
where $\nu$ is the normal bundle of the diagonal $\Delta$ in $M^2$
and these maps have to agree on their ``common domain''. The
canonical map $(x,y)\mapsto \mfm_{\{x,y\}}$ will certainly do for
the first (see Proposition
\ref{int_configuration_space_inclusion}). Choose a Riemannian
metric on $M$. It is well-known that $\nu\cong TM$ and that the
inclusion $\nu\cong TM\hookrightarrow M^2$ can be chosen to be
$(x,v)\mapsto(x,\exp_xv)$ where $\exp$ is defined via
geodesics\footnote{More precisely $t\mapsto\exp_x(tv)$ is the
unique geodesic starting at $x$ with speed $v$.}. We define
 {\setlength\arraycolsep{2pt}
 \begin{eqnarray*}
  \varphi:\eta\nu\times\mcR=\bbR_+\times STM\times\mcR & \ra & \bbR^2
  \\
  (t,(x,v),f) & \mapsto & \biggl(f(x),\frac{f(\exp_x(tv))-f(x)}{t}\biggr)
 \end{eqnarray*}}
where for $t=0$ the second coordinate is to be interpreted as a
limit for $t\ra 0$, i.e. as $df(x)(v)$. As $\varphi$ satisfies all
the hypotheses of Lemma \ref{int_smooth_maps_into_Md} we conclude
that the corresponding map $\psi:\eta\nu\ra\scM_2(M)$ is smooth
(in fact an immersion by an easy computation). The two maps
clearly agree on their common domain and thus define a map
 \[(M^2)_\Delta\lra\scM_2(M)\]
Also it is clear that this map is invariant under the involution
$\hat{\tau}$ and finally induces a map
$(M^2)_{\hat{\tau}}\ra\scM_2(M)$. From the set theoretic
description of $(M^2)_{\hat{\tau}}$ as a disjoint union
$M^{[2]}\sqcup(STM/\{\pm1\})$ it is easily seen to be bijective:
the ideals $I\in\scM_2(M)$ have either 2 points in their spectrum
and then they lie in $M^{[2]}$, or only 1 point, say $x$, in which
case $(\mfm_x)^2\subseteq I\subseteq\mfm_x$ and therefore $I$ can
be thought of as a hyperplane in $\mfm_x/(\mfm_x)^2\cong(T_xM)^*$.
These correspond to one-dimensional subspaces of $T_xM$, i.e.
elements of $(STM)_x/\{\pm1\}$.

As both spaces are compact Hausdorff this map is a homeomorphism
(in fact a diffeomorphism). Again note that $\scM_2(M)$ is always
a manifold with corners. \qed
\end{example}


\begin{example}[a description of $\scM_3(S^m)$]
Let us first classify all ideals
 \[I\subseteq J^2_0(\bbR^m,\bbR)\]
of codimension $3$. We claim that there exist only the following two
types of such ideals. The ideals of the first type are determined by
an orbit of $1$-jets of immersions
$\sigma:(\bbR^2,0)\hookrightarrow(\bbR^m,0)$ under the action of the
origin preserving diffeomorphism group on $\bbR^2$. The
corresponding ideal is
 \[I_{[j^1_0\sigma]}=\{j^2_0f\ |\ j^1_0(f\sigma)=0\}\]
It can be easily seen that these ideals are all distinct. Similarly
the second type is determined by an orbit of $2$-jets of immersions
$\gamma:(\bbR,0)\hookrightarrow(\bbR^m,0)$ and the corresponding
ideal is
 \[I_{[j^2_0\gamma]}=\{j^2_0f\ |\ j^2_0(f\gamma)=0\}\]
Again all these ideals are distinct.

We prove the claim by reduction to ideals in $J^1_0(\bbR^m,\bbR)$.
Denoting the maximal ideal of $\mcJ:=J^2_0(\bbR^m,\bbR)$ by $\mfm$
we identify $J^1_0(\bbR^m,\bbR)$ with $\mcJ/\mfm^2$. Therefore if
$I$ is an ideal of $\mcJ$ of codimension $3$ its image in
$J^1_0(\bbR^m,\bbR)$ is $(I+\mfm^2)/\mfm^2$ and its codimension must
be $2$ or $3$ by an easy inspection. As we saw at the end of the
last example for codimension $2$ in suitable linear coordinates we
can write
 \[I+\mfm^2=\langle x_2,\ldots,x_m,x_1^2\rangle\]
Therefore $x_2+q_2,\ldots,x_m+q_m\in I$ where $q_2,\ldots,q_m$ are
some quadratic functions. Necessarily
 \[\langle x_2,\ldots,x_m\rangle\mfm\subseteq I\]
and therefore we can assume $q_i=c_i\cdot x_1^2$. By another linear
change of coordinates we can achieve
 \[\langle x_2+c\cdot x_1^2,x_3,\ldots,x_m\rangle\subseteq I\]
As the codimensions are equal we must have an equality and the left
hand side is the ideal $I_{[j^2_0\gamma]}$ for $\gamma(t)=(t,-c\cdot
t^2,0,\ldots,0)$.

If the codimension of $I+\mfm^2$ is $3$ then for dimensional reasons
$I=I+\mfm^2$ and analogously to the previous case we have
 \[I=\langle x_3,\ldots,x_m,x_1^2,x_1x_2,x_2^2\rangle\]
in suitable linear coordinates. For $\sigma(s,t)=(s,t,0,\ldots,0)$
we get $I=I_{[j^1_0\sigma]}$.

Now we will explain how $\scM_3(S^m)$ is identified with the space
 \[X:=\{(p,x_1,x_2,x_3)\ |\ p\textrm{ an affine plane in }\bbR^{m+1}; x_1,x_2,x_3\in p\cap S^m\}/\Sigma_3\]
where $\Sigma_3$ acts by permuting $x_1,x_2,x_3$. First we construct
a map $\psi:X\ra\scM_3(S^m)$ by specifying it on various subsets:
\begin{itemize}
\item{
When $p$ is tangent to $S^m$ then
necessarily $x_1=x_2=x_3$ and $p$, being a $2$-dimensional subspace
of $T_{x_1}S^m$, can be thought of as an orbit of $1$-jets of
immersions $(\bbR^2,0)\ra(S^m,x_1)$. We assign to the class of
$(p,x_1,x_2,x_3)$ the corresponding ideal $I_{[j^1_0\sigma]}$. 
}
\item{
If $x_1=x_2=x_3$ but $p$ is not tangent to $S^m$ then it intersects
$S^m$ in a circle which can be thought of as an orbit of $2$-jets of
immersions $(\bbR^1,0)\ra(S^m,x_1)$. Again the class is sent to the
corresponding $I_{[j^2_0\gamma]}$. 
}
\item{
When $x_1=x_2\neq x_3$ then $p$
prescribes a tangential line $l$ at $x_1$ and the class is sent to
$I_l\cdot\mfm_{x_3}$ where $I_l$ is the ideal of codimension $2$
corresponding to $l\in(STM)_{x_1}/\{\pm1\}$ as in the previous
example. 
}
\item{
If all $x_1,x_2,x_3$ are distinct then the image will be
taken to be $\mfm_{\{x_1,x_2,x_3\}}$.
}
\end{itemize}

Both $X$ and $\scM_3(S^m)$ are compact. We postpone the proof of
continuity of $\psi$ and show that it is bijective. If all
$x_1,x_2,x_3$ are distinct then the plane $p$ is determined by these
three points as well as it is by a point and a tangential direction
at a different point. This exhausts all the ideals of codimension
$3$ with spectrum consisting of more than one point. On them $\psi$
is a bijection. We classified all the ideals of codimension $3$ with
spectrum a singleton above. On ideals of type $I_{[j^1_0\sigma]}$ we
clearly get a bijection too. The same is true of ideals of type
$I_{[j^2_0\gamma]}$ when one observes that orbits of $2$-jets of
immersions $(\bbR,0)\hookrightarrow(S^m,x)$ are in bijection with
circles through $x$.\footnote{This is clear from the following two observations. Firstly a $2$-jet of a path $\gamma:(\bbR,0)\ra(\bbR^{m+1},x)$ is a $2$-jet of $(\bbR,0)\ra(S^m,x)$ iff $j^2_0(\langle\gamma,\gamma\rangle)$ is constant $1$. For any parametrization of a circle $c\subseteq\bbR^{m+1}$ this condition is equivalent to $c\subseteq S^m$ by an easy computation. Applying these considerations to the osculation circle of a path $\gamma:(\bbR,0)\ra S^m$ shows that it always lies on $S^m$ providing the bijection.}

It remains to prove continuity of $\psi$. We use a continuous
version of Lemma \ref{int_smooth_maps_into_Md}. First we need to
resolve $X$ by a bigger space
 \[\begin{split}
 Y:=\{ & (z,u,v,y_1,y_2,y_3,w_{12},w_{23},w_{31})\ |\ z\in D^{m+1}; u,v\in S^m; u\perp z,v\perp z,u\perp v; \\
 & y_i\in S^1; w_{ij}\in S^1; y_i-y_j\in\bbR_+ w_{ij}; y_i=y_j\Rightarrow w_{ij}\perp y_i;\textrm{not all }w_{ij}\textrm{ the same}\}
 \end{split}\]
Here $z,u,v$ prescribe an affine plane $p$ in $\bbR^{m+1}$ with
orthonormal basis (endowing $p$ with an orientation) that clearly
intersects $S^m$. We denote by $\alpha:\bbR^2\ra\bbR^{m+1}$ the
affine map (isometry) sending $0$ to $z$ and the standard basis
vectors to $u$ and $v$. Clearly $\alpha$ depends smoothly on
$z,u,v$. Next we set $\rho=\sqrt{1-|z|^2}$, the radius of the circle
$p\cap S^m$. If $p$ is not tangent to $S^m$ then $y_1,y_2,y_3$ can
be identified with points $x_i=\alpha(\rho\cdot y_i)$ in $p\cap
S^m$. Otherwise they are thought of as directions in the tangent
space (and clearly $x_i=z$). When $y_i$ are distinct then
$w_{ij}=\frac{y_i-y_j}{|y_i-y_j|}$ and if $\rho>0$ then also
$\alpha(w_{ij})=\frac{x_i-x_j}{|x_i-x_j|}$. If $y_i=y_j$ the vector
$w_{ij}$ prescribes an infinitesimal direction from $y_j$ to $y_i$.
The only condition is that when all three points collapse one of
these infinitesimal directions has to be opposite to the other two.

Clearly $X$ is a quotient space of $Y$ and it is enough to prove
continuity of the map $\tilde{\psi}:Y\ra\scM_3(S^m)$. On the open
subset
 \[Z:=\{(z,u,v,y_1,y_2,y_3,w_{12},w_{23},w_{31})\in Y\ |\ |z|<1; y_1,y_2,y_3\textrm{ all distinct}\}\]
we have a continuous map $\varphi:\mcR\times Z\ra\bbR^3$ defined by sending $(f,P)$ to
 \[\left(f(x_1),\frac{|f(x_1)-f(x_2)|}{|x_1-x_2|},\rho\cdot\frac{|f(x_1)(x_2-x_3)+f(x_2)(x_3-x_1)+f(x_3)(x_1-x_2)|}{|x_2-x_3|\cdot|x_3-x_1|\cdot|x_1-x_2|}\right)\]

To show how $\varphi$ continuously extends to $Y$ we use our version
of Taylor expansion (Corollary \ref{int_taylor_expansion}) on some
extension of $f$ to $\bbR^{m+1}$.
 \[\begin{split}
  f(x_2)=\ & f(x_1)+\mcI(f,[x_1,x_2])(x_2-x_1) \\
  f(x_3)=\ & f(x_1)+\mcI(f,[x_1,x_2])(x_3-x_1)+1/2\cdot\mcI(f,[x_1,x_2,x_3])(x_3-x_1,x_3-x_2)
 \end{split}\]
and for brevity we put
$a=\mcI(f,[x_1,x_2]),b=\mcI(f,[x_1,x_2,x_3])$. Then easily
 \[\frac{|f(x_1)-f(x_2)|}{|x_1-x_2|}=|a(\alpha(w_{12}))|\]

 \[\begin{split}
  & f(x_1)(x_2-x_3)+f(x_2)(x_3-x_1)+f(x_3)(x_1-x_2) \\
  =\ & a(x_1-x_2)(x_2-x_3)-a(x_2-x_3)(x_1-x_2)-1/2\cdot
  b(x_3-x_1,x_2-x_3)(x_1-x_2)
 \end{split}\]
Now if we denote $a^*=a(u)u+a(v)v$ the projection of the gradient of
$a$ into the plane spanned by $u$ and $v$ we can rewrite
 \[\begin{split}
  & a(x_1-x_2)(x_2-x_3)-a(x_2-x_3)(x_1-x_2) \\
  =\ & \langle a^*,x_1-x_2\rangle(x_2-x_3)-\langle
  a^*,x_2-x_3\rangle(x_1-x_2) \\
  =\ & \vol(x_1-x_2,x_2-x_3)Ra^*
 \end{split}\]
where $\vol$ denotes the oriented volume and $R$ rotation in $p$ by
$+\pi/2$ so that $Ra^*=a(u)v-a(v)u$. Easily $\vol(x_1-x_2,x_2-x_3)$
is twice the area of the triangle $x_1x_2x_3$. It is well known that
this can be expressed by the lengths of the sides and the radius of
the circumscribed circle as
 \[a(x_1-x_2)(x_2-x_3)-a(x_2-x_3)(x_1-x_2)=\sign\cdot 1/(2\rho)\cdot|x_1-x_2||x_2-x_3||x_3-x_1|Ra^*\]
Here $\sign$ depends on the orientation of the triangle $x_1x_2x_3$ inside $p$ which on the other hand depends continuously on the data\footnote{Easily $Y$ decomposes into a disjoint union of two subspaces according to the orientation of $y_1y_2y_3$ in $\bbR^2$ (this is clear for nondegenerate triangles and $Y$ is designed so as to remember via $w_{12},w_{23},w_{31}$ the orientation of degenerate triangles). This decomposition prescribes a continuous function $\sign:Y\ra\bbZ/2=\bbZ^*$}. Therefore the last
component of $\varphi$ is
 \[1/2\cdot|\sign\cdot Ra^*-\rho\cdot b(\alpha(w_{31}),\alpha(w_{23}))\alpha(w_{12})|\]
These formulas obviously extend continuously to $Y$ and Lemma \ref{int_smooth_maps_into_Md} then provides a continuous map
 \[\tilde{\psi}:Y\ra\scM_3(S^m)\]
It is easy to verify that $\psi$ factors through the quotient $X$ of $Y$ yielding the above described map $\psi$ and showing that it is indeed a homeomorphism.

\end{example}

\begin{example}[description of $\scM_1(\mcV)$ for a special $\mcV$]
Let $N\subseteq M$ be a closed neat\footnote{see \cite{Hirsch}}
submanifold and
 \[\mcV=\{f\in C^\infty(M,\bbR)\ |\ f(z)=0\ \forall z\in N\}\]
We show later in Example
\ref{int_example_of_topologically_finitely_generated_modules} that
$\mcV$ is topologically finitely generated. We can argue similarly
as in the last example to conclude that
 \[\scM_1(\mcV)\cong (M-N)\sqcup(S\nu/\{\pm1\})\]
where the right hand side is given a smooth structure in a way
similar to Construction \ref{int_construction_of_blow_up} but
using the projective bundle with its canonical line bundle over it
in place of the sphere bundle. This construction is called a
(real) blow-up (cf. exercises 7-8 after Chapter 12 in
\cite{BrockerJanich}).
\end{example}

Our next goal is to describe certain subsets of $\scM_d(M)$. They
are subsets of ideals of a ``fixed type'' and are injectively
immersed submanifolds. Also every ideal has some (unique) type and
so $\scM_d(M)$ is in fact a disjoint union (over all possible
types) of these submanifolds.

We take the following construction from section 35 of \cite{KMS}.
A \emph{Weil algebra} is a finite dimensional associative,
commutative algebra $A$ over $\bbR$ with a unit such that
$A=\bbR\oplus N$ where $N$ is the ideal of nilpotent elements.
Equivalently it could be described as a quotient algebra of
$J^r_0(\bbR^m,\bbR)$ for some $r$ and $m$. Therefore we can get
all ideals of $\mcR$ with a spectrum consisting of a single point
as kernels of surjective algebra homomorphisms $\mcR\ra A$ for
some Weil algebra $A$ (namely $A$ is the quotient of $\mcR$ by
that ideal). We give the set of all such homomorphisms (surjective
or not) a smooth structure in such a way that the map sending such
a homomorphism to the spectrum of its kernel is a bundle
projection. This bundle is called the \emph{Weil bundle}
associated to $A$.

We first give a construction of this bundle and then show that its
points can be indeed identified with homomorphism $\mcR\ra A$. We
start with the restriction
 \begin{equation} \label{int_jet_bundle}
  J^r_{0,\mathrm{diff}}(\bbR^m,M)\ra M
 \end{equation}
of the jet bundle $J^r(\bbR^m,M)\ra M$ to the subspace of all
invertible jets with source~$0$. Setting
$G^r_m:=J^r_{0,\mathrm{diff}}(\bbR^m,\bbR^m)_0$, the Lie group of
all invertible jets with source and target~$0$, we see that
(\ref{int_jet_bundle}) is a principal $G^r_m$-bundle and so we can
define
 \begin{equation} \label{int_definition_of_TAM}
  T_AM:=J^r_{0,\mathrm{diff}}(\bbR^m,M)\times_{G^r_m}\Hom_{\mathrm{alg}}(J^r_0(\bbR^m,\bbR),A)
 \end{equation}
This clearly expresses $T_AM$ as a smooth bundle over $M$ with
fibre
 \[\Hom_{\mathrm{alg}}(J^r_0(\bbR^m,\bbR),A)\cong N^m\]
Moreover we have a bijection defined in terms of
(\ref{int_definition_of_TAM}) by the formula
 {\setlength\arraycolsep{2pt}
 \begin{eqnarray*}
  T_AM & \xra{\cong} & \Hom_{\mathrm{alg}}(\mcR,A)
 \\
  {}[j^r_xg,\varphi] & \mapsto & \left(\mcR\xra{j^r_x}
  J^r_x(M,\bbR)\xra{g^*}J^r_0(\bbR^m,\bbR)\xra{\varphi} A\right)
 \end{eqnarray*}}
It is easily seen to be a bijective correspondence (that a kernel
of any $\mcR\ra A$ has a spectrum consisting of only a single
point follows from the fact that in $A$, the only idempotents are
$0$ and $1$). We have a subbundle
 \[
  \check{T}_AM:=J^r_{0,\mathrm{diff}}(\bbR^m,M)\times_{G^r_m}\SurHom_{\mathrm{alg}}(J^r_0(\bbR^m,\bbR),A)
 \]
which then corresponds to surjective algebra homomorphisms
$\mcR\ra A$.

One says that an ideal $I$ is of type $A$ if $\mcR/I\cong A$. An
ideal $I$ of type $A$ can then be identified with a class of
surjective homomorphisms $\mcR\ra A$, namely with the class of all
those homomorphisms that have kernel $I$. In this way we get a
space $J^AM$ of all ideals of type $A$ as a certain quotient of
$\check{T}_AM$. A crucial observation in \cite{Alonso} is that the
action of $G^r_m$ on
$\SurHom_{\mathrm{alg}}(J^r_0(\bbR^m,\bbR),A)$ is transitive and
so, after choosing some
$\alpha_0\in\SurHom_{\mathrm{alg}}(J^r_0(\bbR^m,\bbR),A)$, one can
identify it with the quotient of $G^r_m$ by the stabilizer of
$\alpha_0$. In the same terminology an ideal in
$J^r_0(\bbR^m,\bbR)$ of type $A$ is a class of $G^r_m$ modulo the
stabilizer of $\ker\alpha_0$. Therefore the space of ideals of
type $A$ can be identified with the smooth bundle
 \[
  J^AM\cong J^r_{0,\mathrm{diff}}(\bbR^m,M)/\mathrm{St}(\ker\alpha_0)\lra M
 \]
and clearly the smooth structure does not depend on the choice of
$\alpha_0$.

\begin{proposition} \label{int_Weil_bundle_inclusion}
Let $D\pitchfork\scM_{d+1}(M)$ be a finite dimensional linear
subspace of $\mcR$ and let $A$ be a $d$-dimensional Weil algebra.
Then the inclusion
 \[\iota:J^AM\subseteq \scM_d(M)\subseteq G_d(D)\]
is an injective immersion.
\end{proposition}

\begin{proof}
This is another application of Lemma
\ref{int_smooth_maps_into_Md}. Consider the map
 {\setlength\arraycolsep{2pt}
 \begin{eqnarray*}
  \varphi:\mcR\times J^r_{0,\mathrm{diff}}(\bbR^m,M) & \ra & A
 \\
  (f,j^r_x(g)) & = & \alpha_0(j^r_x(fg))
 \end{eqnarray*}}
Clearly this map is surjective linear in the first and smooth in
the second variable and hence in the sense of Lemma
\ref{int_smooth_maps_into_Md} it defines
 \[\psi:J^r_{0,\mathrm{diff}}(\bbR^m,M)\ra\scM_d(M)\subseteq G_d(D)\]
which is also smooth. As we have a commutative diagram
 \[\xymatrix{
  J^r_{0,\mathrm{diff}}(\bbR^m,M) \ar[r]^{\psi} \ar@{->>}[d] &
  **{!<-25pt,0pt>+}\scM_d(M)\subseteq G_d(D)
 \\
  J^r_{0,\mathrm{diff}}(\bbR^m,M)/\mathrm{St}(\ker\alpha_0)
  \ar@{-->}[ru]_{\iota}
 }\]
in order to show that the dashed arrow is an immersion we need to
identify $\ker\psi_*$. Lemma \ref{int_smooth_maps_into_Md} gives
an answer in terms of the kernel of the differential (say at
$j^r_xg$) of the map $\varphi(f,-)$ which can be decomposed as
 \[J^r_{0,\mathrm{diff}}(\bbR^m,M)\xlra{f_*}
 J^r_0(\bbR^m,\bbR)\xlra{\alpha_0} A\]
To give a tangent vector in
$T_{j^r_0g}J^r_{0,\mathrm{diff}}(\bbR^m,M)$ is the same as to give
an element of $T_{\id}J^r_{0,\mathrm{diff}}(\bbR^m,\bbR^m)$ and then
compose with $g$. The elements of
$T_{\id}J^r_{0,\mathrm{diff}}(\bbR^m,\bbR^m)$ arise from vector
fields. Therefore let $X:\bbR^m\ra T\bbR^m$ be a local vector field
with a local flow
 \[\gamma:\bbR^m\times\bbR\lra \bbR^m\]
Under our identifications it defines a tangent vector
 \[\hat{X}:=\left.\frac{d}{dt}\right|_{t=0}j^r_0(g\gamma(-,t))
 \in T_{j^r_0g}J^r_{0,\mathrm{diff}}(\bbR^m,M)\]
Then for each $f\in\mcR$ we get
 \[d(\alpha_0f_*)(\hat{X})=
 \left.\frac{d}{dt}\right|_{t=0}\alpha_0(j^r_0(fg\gamma(-,t)))=
 \alpha_0(j^r_0(X(fg)))\]
Suppose that $X(0)\neq 0$. Then we claim that there exists an
$f\in\psi(j^r_xg)$ for which this expression is nonzero as well. In
other words such $\hat{X}$ can never lie in $\ker\psi_*$. We
postpone the proof of this claim and thus assume that the only
$\hat{X}$ which could produce an element in this kernel are the
vectors tangent to the submanifold
$J^r_{0,\mathrm{diff}}(\bbR^m,M)_x$. The Lie group $G^r_m$ acts
simply transitively on this space. Let $Y$ be an element of the Lie
algebra of $G^r_m$. Then we obtain a vector field (with $p$ running
over $J^r_{0,\mathrm{diff}}(\bbR^m,M)_x$)
 \[Y^+(p):=\left.\frac{d}{dt}\right|_{t=0}p\cdot\exp(tY)\]
and we also have similar vector fields on $J^r_0(\bbR^m,\bbR)$. The
restriction of $\varphi(f,-)$ is simply the composition
 \[J^r_{0,\mathrm{diff}}(\bbR^m,M)_x\xlra{f_*}J^r_0(\bbR^m,\bbR)
 \xlra{projection}J^r_0(\bbR^m,\bbR)/\ker\alpha_0\cong A\]
Therefore $Y^+(j^r_0g)\in\ker\psi_*$ if and only if for each
$f\in\psi(j^r_0g)$ we have $df_*(Y^+(j^r_0g))\in\ker\alpha_0$. As
the map $f_*$ is $G^r_m$-equivariant we can rewrite
 \[df_*(Y^+(j^r_0g))=Y^+(f_*j^r_0g)=Y^+(j^r_0(fg))\]
for the corresponding canonical vector field on
$J^r_0(\bbR^m,\bbR)$. Now observe that we get all possible values
$j^r_0(fg)\in\ker\alpha_0$ by varying $f$ over $\psi(j^r_0g)$.
Therefore $Y^+(j^r_0g)\in\ker\psi_*$ if and only if
$Y^+(\ker\alpha_0)\subseteq\ker\alpha_0$. These $Y$ clearly
constitute the Lie algebra of the stabilizer
$\mathrm{St}(\ker\alpha_0)$ and therefore we conclude that
$\ker\psi_*$ is exactly the vertical tangent bundle of
 \[J^r_{0,\mathrm{diff}}(\bbR^m,M)\lra
 J^r_{0,\mathrm{diff}}(\bbR^m,M)/\mathrm{St}(\ker\alpha_0)\cong J^AM\]
Consequently $\psi$ induces on the quotient $J^AM\cong
J^r_{0,\mathrm{diff}}(\bbR^m,M)/\mathrm{St}(\ker\varphi_0)$ an
immersion $\iota:J^AM\lra G_d(D)$.

Now we prove the remaining claim. Because we assume that $X(0)\neq
0$ we can find a local diffeomorphism $h:(\bbR^m,0)\ra(\bbR^m,0)$
such that $h^*X=\partial_{x_1}$. Then
 \[X(fg)=\partial_{x_1}(fgh)\circ h^{-1}\]
We denote by $K$ the kernel of
 \[J^r_0(\bbR^m,\bbR)\xlra{(h^{-1})^*}J^r_0(\bbR^m,\bbR)\xlra{\varphi_0}A\]
and we are looking for $f\in\mcR$ such that $j^r_0(fgh)\in K$ but
$j^r_0(\partial_{x_1}(fgh))\not\in K$. Let
$K\subseteq(x_1^k,x_2,\ldots,x_m)$ but
$K\not\subseteq(x_1^{k+1},x_2,\ldots,x_m)$. Both $g$ and $h$ being
diffeomorphisms there exists $f\in\mcR$ such that $j^r_0(fgh)\in
K-(x_1^{k+1},x_2,\ldots,x_m)$. Then $fgh=x_1^k\cdot\lambda$ modulo
$(x_2,\ldots,x_m)$ with $\lambda(0)\neq 0$ and it is easy to see
that $\partial_{x_1}(fgh)=x_1^{k-1}\mu$ modulo $(x_2,\ldots,x_m)$
with $\mu(0)\neq 0$ so that $j^r_0(\partial_{x_1}(fgh))$ does not
lie in $(x_1^k,x_2,\ldots,x_m)$ and in particular it does not lie in
$K$.
\end{proof}

\begin{question}
Is the inclusion $\iota:J^AM\hookrightarrow\scM_d(M)$ an
embedding? Quite easily (reducing to a local question and using
polynomials) one can reduce this problem to the question of the
canonical map
 \[G^r_m/(G^r_m\cap\mathrm{St}(\ker\varphi_0))\hookrightarrow Gl(J^r_0(\bbR^m,\bbR))/\mathrm{St}(\ker\varphi_0)\]
being an embedding
.
\end{question}

An easy generalization to the case of finitely many Weil algebras
$A_i$, $i=1,\ldots,n$, produces a bundle
 \[
  \check{T}_{A_1,\ldots,A_n}M=
  \left.\left(\check{T}_{A_1}M\times\cdots\times\check{T}_{A_n}M\right)\right|_{M^{(n)}}\lra
  M^{(n)}
 \]
togeter with a bijection
$\check{T}_{A_1,\ldots,A_n}M\cong\SurHom_{\mathrm{alg}}(\mcR,A_1\times\cdots\times
A_n)$. Every ideal $I\in\scM_d(M)$ of type $A_1\times\cdots\times
A_n$ can be clearly recovered as a kernel of such surjective
homomorphism. In this way we get a space $J^{A_1,\ldots,A_n}M$ of
ideals of a fixed type $A_1\times\cdots\times A_n$ as a quotient of
$\check{T}_{A_1,\ldots,A_n}M$. Moreover we can again identify
$J^{A_1,\ldots,A_n}M$ with a quotient\footnote{Here
$(J^r_{0,\mathrm{diff}}(\bbR^m,M))^{(n)}_{\ushort{M}}$ denotes the
restriction of the power $(J^r_{0,\mathrm{diff}}(\bbR^m,M))^n\ra
M^n$ of the jet bundle to the subspace $M^{(n)}\subseteq M^n$. In
particular as the original bundle was a principal $G^r_m$-bundle the
resulting bundle over $M^{(n)}$ will be a principal
$(G^r_m)^n$-bundle on which there is an action of the symmetric
group $\Sigma_n$. Taking the quotient by this action one gets a
principal $(G^r_m\wr\Sigma_n)$-bundle
$(J^r_{0,\mathrm{diff}}(\bbR^m,M))^{(n)}_{\ushort{M}}\ra M^{[n]}$.}
 \[(J^r_{0,\mathrm{diff}}(\bbR^m,M))^{(n)}_{\ushort{M}}
 \times_{G^r_m\wr\Sigma_n}(G^r_m\wr\Sigma_n)/\mathrm{St}(\ker\alpha_0)\]
for some (any) surjective homomorphism
$\alpha_0:(J^r_0(\bbR^m,\bbR))^n\ra A_1\times\cdots\times A_n$. The
canonical inclusion map
$\iota:J^{A_1,\ldots,A_n}M\hookrightarrow\scM_d(M)$ with $d=\dim
A_1+\cdots+\dim A_n$ is an injective immersion. For a proof observe
that locally $(J^r_{0,\mathrm{diff}}(\bbR^m,M))^{(n)}_{\ushort{M}}$
is just a product of $J^r_{0,\mathrm{diff}}(\bbR^m,M)$ and so one
can almost copy the proof of Proposition
\ref{int_Weil_bundle_inclusion} (also see the proof of Proposition
\ref{int_configuration_space_inclusion}).


\section{A criterion for topological finite generation}

In this section we give a criterion for $\mcV$ to be a topologically
finitely generated $\mcR$-module in terms of its local structure
allowing us to prove that some interesting topological
$\mcR$-modules are topologically finitely generated.

First we observe that $\mcR$ is the space of global sections of the
sheaf $\mcR(U)=C^\infty(U,\bbR)$ of (locally convex) topological
$\bbR$-algebras. The sheaf property in the realm of topological
spaces is proved in Lemma \ref{fun_sheaf_property} whereas the
algebraic part is clear. Suppose that $\mcV(U)$ is a sheaf of
topological $\mcR(U)$-modules and denote $\mcV=\mcV(M)$. We have the
following result.

\begin{proposition}
An $\mcR$-module $\mcV$ is topologically finitely generated if there
is an open covering $\scU$ of $M$ such that for every $U\in\scU$,
$\mcV(U)$ is a topologically finitely generated $\mcR(U)$-module.
\end{proposition}

\begin{proof} First we construct an auxiliary map enabling us to
extend sections. Let $U$ be an open subset of $M$ and let
$\lambda:M\rightarrow\bbR$ be a smooth function with
$\supp(\lambda)\subseteq U$. We define a $\lambda$-extension map
$e_\lambda:\mcV(U)\rightarrow\mcV$ which is simply multiplication by
$\lambda$ and extending by $0$. To be more precise let
$V=M-\supp(\lambda)$. By the sheaf property we have a pullback
diagram
\[\xymatrix{
{}\mcV \ar[r]^{res} \ar[d]_{res} & {}\mcV(U) \ar[d]^{res} \\
{}\mcV(V) \ar[r]^-{res} & {}\mcV(U\cap V)}\] Now we can define
$e_\lambda:\mcV(U)\rightarrow\mcV$ by
\[(l_\lambda,0)^T:\mcV(U)\lra\mcV(U)\times_{\mcV(U\cap V)}\mcV(V)\cong\mcV\]
the first component being multiplication by $\lambda$. By the
universal property of pullbacks the composition
$\mcV\xlra{res}\mcV(U)\xlra{e_\lambda}\mcV$ is multiplication by
$\lambda$ as is the other composition $res\circ e_\lambda$.

As an application we prove that if $A$ is a closed subset of $M$
such that $A\subseteq U$ then $\mcR/\mfn_A\cong\mcR(U)/\mfn_A$. A
map in one direction is induced by the restriction map and in the
other direction by $e_\lambda$ where $\lambda:M\rightarrow\bbR$ is
any smooth map supported in $U$ such that $\lambda=1$ on a
neighbourhood of $A$.

Now we proceed with the actual proof. Let
$\scU=\{U_1,\ldots,U_l\}$ be any open covering of $M$ and let
$\lambda_1,\ldots,\lambda_l$ be a subordinate partition of unity.
We define a gluing map
\[(e_{\lambda_1},\ldots,e_{\lambda_l}):\mcV(U_1)\times\cdots\times\mcV(U_l)\twoheadrightarrow\mcV\]
which is a quotient map as it has a section
$(res,\ldots,res)^T:\mcV\rightarrow\mcV(U_1)\times\cdots\times\mcV(U_l)$.
Suppose now that for each $i$ there is a topological quotient map
\[\varphi_i:\mcR(U_i)E_i\twoheadrightarrow\mcV(U_i)\]
with $E_i$ a finite dimensional real vector space. Denoting by
$\psi_i$ the composition
\[\mcR E_i\xlra{res}\mcR(U_i)E_i\xlra{\varphi_i}\mcV(U_i)\xlra{e_{\lambda_i}}\mcV\]
we claim that the map $(\psi_1,\ldots,\psi_l):\mcR
E_1\times\cdots\times\mcR E_l\lra\mcV$ is also a topological
quotient map. This follows from the diagram
\[\xymatrix{
{}\prod\mcR E_i \ar[r]^-{\prod res} \ar@{->>}[d] &
{}\prod\mcR(U_i)E_i \ar@{->>}[r]^-{\prod\varphi_i} \ar@{->>}[d] &
{}\prod\mcV(U_i) \ar@{->>}[d]^{(e_{\lambda_1},\ldots,e_{\lambda_l})} \\
{}\prod(\mcR/\mfn_{A_i})E_i \ar[r]^-\cong &
{}\prod\bigl(\mcR(U_i)/\mfn_{A_i}\bigr)E_i \ar@{-->}[r] &
{}\mcV}\] where $A_i=\supp(\lambda_i)\subseteq U_i$. Here the
dashed arrow exists and is clearly a quotient map. As $\mcR
E_1\times\cdots\times\mcR E_l\cong\mcR(E_1\times\cdots\times E_l)$
this finishes the proof.
\end{proof}

%

\begin{example}
\label{int_example_of_topologically_finitely_generated_modules}
Let $M$ be a compact manifold, let $N\subseteq M$ be a closed
submanifold and suppose that either $N$ is neat or that
$N=\partial M$. Let $r$ be a positive integer. We define a
submodule $\mcV\subseteq\mcR$ as
\[\mcV=\{f\in C^\infty(M,\bbR)\ |\ j^{r-1}_z f=0\ \forall z\in N\}\]
Then $\mcV$ is a topologically finitely generated $\mcR$-module.
\end{example}

\begin{proof} We will give the proof in the case $\partial
M=\emptyset$ to simplify the notation. Clearly $\mcV$ is a subsheaf
of $\mcR$ and thus a sheaf itself. We cover $M$ by coordinate charts
$\bbR^m\cong U\subseteq M$ under which either $N=\emptyset$ or
\[N=\bbR^{m_0}\subseteq\bbR^{m_0}\times\bbR^{m_1}\cong\bbR^m\]
We are left to show that each $\mcV(U)$ is a topologically finitely
generated $\mcR(U)$-module. Using the charts we can assume that
$U=\bbR^m$ and in fact we will be assuming that $M=U=\bbR^m$ to
avoid writing $U$ everywhere.

When $N=\emptyset$ we have $\mcV=\mcR$ and there is nothing to
prove. Hence let us assume that $N=\bbR^{m_0}$. We define $E$ to
be the space of homogenous polynomials on $\bbR^m$ of degree $r$
that belong to $\mcV$, i.e. those polynomials that do not depend
on the $\bbR^{m_0}$-coordinate. This gives us an embedding
$E\rightarrow\mcV$ and we define
\[\varphi:\mcR E\rightarrow\mcV\]
to be the unique $\mcR$-module map extending this embedding. In
other words, if we think of $\mcR E$ as $C^\infty(\bbR^m,E)$ then
$g:\bbR^m\rightarrow E$ is send by $\varphi$ to the function
$f:z\mapsto(g(z))(z)$. To prove that it is a topological quotient
map we will construct a continuous section (which will happen to
be linear but will not be an $\mcR$-module map. In fact there is
no $\mcR$-section).

Let $(x,y)\in\bbR^{m_0}\times\bbR^{m_1}\cong\bbR^m$. We define an
$r$-simplex $\sigma_{(x,y)}$ to be
\[\sigma_{(x,y)}=[(x,0),\ldots,(x,0),(x,y)]\]
and then for $f\in\mcR=C^\infty(\bbR^m,\bbR)$ we set
\[\mcJ(f,x,y):=\frac{1}{r!}\cdot
\int_{\Delta^r}f^{(r)}_{\bbR^{m_1}}\sigma_{(x,y)}\in\Hom(S^r\bbR^{m_1},\bbR)\]
where $f^{(r)}_{\bbR^{m_1}}$ is the $r$-fold derivative of $f$ in
the direction $\bbR^{m_1}$.

If we assume that $f\in\mcV$ then by Lemma
\ref{int_integration_formula}
\[f(x,y)=\mcJ(f,x,y)(y,\ldots,y)\]
Hence if we identify $\Hom(S^r\bbR^{m_1},\bbR)$ with the
polynomial space $E$ then we get
\[\mcJ(f,x,y)\in E\]
with the property that
\[f(x,y)=\mcJ(f,x,y)(x,y)\]
Clearly each $\mcJ(f,-):\bbR^m\rightarrow E$ is a smooth function
and so we obtain a section $s:\mcV\ra\mcR E$ of $\varphi$ by
defining
\[s:f\mapsto\mcJ(f,-)\]
The only thing that remains to be checked is the continuity of
this map. Note that the same formula defines a map
$s:\mcR\rightarrow\mcR E$ and hence it is enough to prove that
this extension is continuous. This is an easy exercise in
differential topology: roughly speaking, if the $(r+k)$-jet of $f$
is small then the $k$-jet of $s(f)$ is small as well.
\end{proof}

This readily generalizes to sections of a smooth vector bundle
$F\rightarrow M$:
\[\mcV=\{f\in\Gamma F\ |\ j^{r-1}_z f=0, z\in N\}\]
is a topologically finitely generated $\mcR$-module. Here the
condition $j^{r-1}_z f=0$ means that $f$ and the zero section
$0:M\rightarrow F$ have the same $(r-1)$-jet at $z$. Another
generalization is possible. If we are given a finite collection
$N_i$ of submanifolds (such that either $N_i$ is neat or
$N_i=\partial M$) in general position (i.e. the tangent spaces at
any intersection point are in general position) and positive
integers $r_i$ then the $\mcR$-module
\[\mcV=\bigcap\limits_i\ \{f\in \Gamma F\ |\ j^{r_i-1}_z f=0\ \forall z\in N_i\}\]
is a topologically finitely generated $\mcR$-module. We will
indicate the necessary changes in the case of two submanifolds
(that intersect transversely). We can assume that $F=\bbR\times M$
and
\[M=\bbR^m\cong\bbR^{m_0}\times\bbR^{m_1}\times\bbR^{m_2}\]
and that $N_1=\bbR^{m_0}\times\{0\}\times\bbR^{m_2}$,
$N_2=\bbR^{m_0}\times\bbR^{m_1}\times\{0\}$. Take $E$ to be the
space of homogenous polynomials on $\bbR^m$ of degree $r_1+r_2$
that belong to $\mcV$. Again we get an $\mcR$-module map
\[\varphi:\mcR E\rightarrow\mcV\]
which under the identification $\mcR E\cong C^\infty(M,E)$ sends
$g:M\rightarrow E$ to $f:z\mapsto(g(z))(z)$. To construct the
section of $\varphi$ we define for a point
$(w,x,y)\in\bbR^{m_0}\times\bbR^{m_1}\times\bbR^{m_2}\cong\bbR^m$
an affine map
\[\sigma_{(w,x,y)}:\Delta^{r_1}\times\Delta^{r_2}\cong\{*\}\times\Delta^{r_1}\times\Delta^{r_2}
\xlra{w\times[0,\ldots,0,x]\times[0,\ldots,0,y]}
\bbR^{m_0}\times\bbR^{m_1}\times\bbR^{m_2}\] and we set
\[\mcJ(f,w,x,y)=\int_{\Delta^{r_1}\times\Delta^{r_2}}f^{(r_1),(r_2)}_{\bbR^{m_1},\bbR^{m_2}}\sigma_{(w,x,y)}
\in\Hom(S^{r_1}\bbR^{m_1}\otimes S^{r_2}\bbR^{m_2},\bbR)\] where
$f^{(r_1),(r_2)}_{\bbR^{m_1},\bbR^{m_2}}$ denotes the
$(r_1+r_2)$-fold derivative, $r_1$-times in the direction
$\bbR^{m_1}$ and $r_2$-times in the direction $\bbR^{m_2}$. Again
we can identify the target space with $E$ and get the desired
section.

\begin{example}
\label{int_finite_generation_of_sections_of_jet_bundle}
Let $F\ra M$ be a vector bundle and $J^rF\ra M$ its jet
prolongation. There are two actions of $\mcR$ on $\Gamma(J^rF)$. The
first one does not use the jet structure and is given in terms of
$f\in\mcR$ and $s\in\Gamma(J^rF)$ by the formula $(f\cdot
s)(x)=f(x)s(x)$. The second is defined in the following way: let
$x\in M$ and let $s(x)=j^r_xg$. Then we set $(f\cdot
s)(x)=j^r_x(fg)$. To distinguish the two actions we denote
$\Gamma(J^rF)$ with the latter one by $j^r\Gamma(J^rF)$.

We claim that $j^r\Gamma(J^rF)$ is a topologically finitely
generated $\mcR$-module.
\end{example}

\begin{proof}
Locally we can reduce to the case $F=\bbR\times M$ and $M=\bbR^m$.
Using trivialization $J^r(\bbR^m,\bbR)\cong
J^r_*(\bbR^m,\bbR)\times\bbR^m$ we have $\Gamma(J^rF)\cong
C^\infty(M,J^r_*(\bbR^m,\bbR))$ and think of
$J^r_*(\bbR^m,\bbR)\subseteq\Gamma(J^rF)$ as the subspace of
constant maps. Clearly with the first action $\Gamma(J^rF)$ is free
on $J^r_*(\bbR^m,\bbR)$. We will show that the same is true for
$j^r\Gamma(J^rF)$. First note that on the kernel of the canonical
map $\alpha:\Gamma(J^rF)\ra\Gamma(J^{r-1}F)$ the two actions
coincide. Denoting by $K_r$ the kernel of $J^r_*(\bbR^m,\bbR)\ra
J^{r-1}_*(\bbR^m,\bbR)$ we obtain a commutative diagram with
$\varphi_r$ the unique $\mcR$-homomorphism extending the above
mentioned inclusion $J^r_*(\bbR^m,\bbR)\subseteq j^r\Gamma(J^rF)$:
 \[\xymatrix{
   0 \ar[r] &
   \mcR K_r \ar[r] \ar[d]_{\psi_r}^\cong &
   \mcR J^r_*(\bbR^m,\bbR) \ar[r] \ar[d]_{\varphi_r} &
   \mcR J^{r-1}_*(\bbR^m,\bbR) \ar[r] \ar[d]_{\varphi_{r-1}} &
   0
  \\
   0 \ar[r] &
   \ker\alpha \ar[r] &
   j^r\Gamma(J^rF) \ar[r]^-\alpha &
   j^{r-1}\Gamma(J^{r-1}F) \ar[r] &
   0
 }\]
As the map on the left is the same as for the first action it must
be an isomorphism. Using a splitting of the top row and an induction
on $r$ it is easy to produce a (continuous!) inverse of $\varphi_r$.
\end{proof}

\section{Noncompact manifolds}

This section is more of an informative one with no complete
proofs. The idea is that for a noncompact $M$ everything should go
the same except the compactness of $\scM_d$ has to be replaced by
the properness of the canonical map $\scM_d\ra\scS_d$.

Let us try to summarize what possible complications could arise.
The first one is that we do not have a characterization of the
maximal ideals in $\mcR$. There are certainly more than just ones
of the form $\mfm_x$ for $x\in M$, the other correspond to some
behaviour ``at infinity'' (the second part of the next proof shows
that such ideals do not contain any proper function; in fact not
even a function with a compact spectrum). Nevertheless the
additional maximal ideals have infinite codimension:

\begin{proof}[Proof (based on the proof of lemma 35.8. of \cite{KMS})]
Suppose that $I$ is a maximal ideal of finite codimension and such
that $\spec(I)=\emptyset$. Take any proper function $h\in\mcR$. As
$\mcR/I$ is finite dimensional, the vectors $h+I,h^2+I,\ldots$
must be linearly dependent and so there is a nonzero finite linear
combination
 \[g=\sum_{i=1}^na_ih^i\in I\]
One can see easily that $g$ is itself proper as it is a
composition
 \[g:M\xlra{h}\bbR\xlra{p}\bbR\]
with $p=\sum a_ix^i$ a non-constant polynomial, hence proper. In
particular $g^{-1}(0)$ is compact. The collection of closed
subsets $\spec(f)=f^{-1}(0)$ with $f$ running over elements of $I$
has by assumption an empty intersection. Therefore (remember
$g^{-1}(0)$ is compact) there is a finite number of elements
$f_1,\ldots,f_n\in I$ such that
 \[\spec(f_1)\cap\cdots\cap\spec(f_n)=\emptyset\]
Then $f_1^2+\cdots+f_n^2\in I$ is positive on $M$ and is therefore
a unit in $\mcR$, a contradiction to maximality.
\end{proof}

Therefore we obtain the same results about the structure of
submodules of finite codimension. We can even find a finite
dimensional transversal (the proofs did not use compactness).

%
Problems arise in the proof of Proposition
\ref{int_injectivity_proposition} (as there exist maximal ideals of
infinite codimension now. Nevertheless the conclusions -- $U_D$
injective and reflecting inclusions -- can be proved easily for
$D\pitchfork\scM_{d+1}$ using intersections) and of continuity of
the natural map $\pi:U_D\scM_d\ra\scS_d$ in Theorem
\ref{int_compactness_for_RE}.
%



\chapter{Simplex on a space} \label{simplex}

The purpose of this chapter is to derive a spectral sequence for
computing the \Cech cohomology with compact supports of the target
of a surjective finite-to-one proper map $f:X\ra Y$. The idea (taken
from \cite{Vassiliev}, \cite{Vassiliev2}) is to replace $Y$ by its
``resolution'' $Rf$. It is a space with the same \Cech cohomology as
$Y$ and with a natural finite filtration such that $R^0f=X$ and
$R^pf-R^{p-1}f$ has an interpretation in terms of $(p+1)$-tuples of
points mapping to a single point. The spectral sequence is then just
a spectral sequence associated with this filtration. The resolution
$Rf$ is obtained from the following construction.

Let $\mcF$ denote the category of finite sets and all maps between
them. There is a functor $\Delta:\mcF\rightarrow\Top$ sending $K$
to
\[\Delta K=\textrm{a ``free'' convex hull of $K$}\]
To be more precise it is a convex hull of $K$ inside the free
$\bbR$-vector space $\bbR^K$ on $K$. We think of $\mcF$ as a
(full) subcategory of the category $\Top$ of compactly generated
Hausdorff spaces. The topological left Kan extension of $\Delta$
along the inclusion yields a functor $\Top\rightarrow\Top$ which
we will also denote by $\Delta$
\[\Delta X:=\int^{K\in\mcF}\map(K,X)\times\Delta K\]
(for coends see Section IX.6 of \cite{MacLane}). The space $\Delta
X$ is what one would call a simplex on $X$. To give an evidence we
describe how one thinks of elements of $\Delta X$ as (free) convex
combinations of points in $X$. For start, points in $\Delta K$ are
convex combinations of elements of $K$ whereas a point in
$\map(K,X)$ identifies these elements with some points in $X$.
Hence the points in the product $\map(K,X)\times\Delta K$ can be
thought of as convex combinations of points in $X$ together with a
labeling of these points by elements of $K$. The coend then
quotients out the relations we would expect: bijections identify
all the possible labelings (so that we are left just with convex
combinations of points in $X$), inclusions tell us that we can
leave out any summand of the form $0\cdot x$ and surjections give
the relation $s\cdot x+t\cdot x=(s+t)\cdot x$. Note that for a
finite discrete space $K$ this agrees with our previous definition
of $\Delta K$.

\begin{remark} There is a different description of this
construction. The functor $\Delta:\Top\rightarrow\ConvTop$ is the
left adjoint of the forgetful functor $\ConvTop\rightarrow\Top$
defined on the category of convex topological spaces (a
topological space $X$ together with a map $I\times X\times
X\rightarrow X$ thought of as the map $(t,x,y)\mapsto(1-t)x+ty$
satisfying certain relations).
\end{remark}

A fibrewise version of this construction for spaces
$f:X\rightarrow Y$ over $Y$ needs just a small modification:
instead of allowing all convex combinations one only considers
combinations of points in a single fibre. In effect one replaces
$\map(K,X)$ by $X_Y^K$, the subspace of those maps $g:K\rightarrow
X$ for which $fg:K\rightarrow Y$ is constant
\[\Delta_YX:=\int^{K\in\mcF}X_Y^K\times\Delta K\]
We have the following pullback diagram (as $\Delta_YX$ is a
subspace of $\Delta X$
)
\[\xymatrix{
{}\Delta_YX \ar[r] \pullback \ar[d] & \Delta X \ar[d]\\
Y \ar[r] & \Delta Y}\]

$\mcF$ has a filtration by full subcategories $\mcF_n$ of sets of
cardinality at most $n+1$. Hence we get a filtration of
$\Delta_YX$ by the $n$-skeleta:
\[\Delta_YX\cong\colim\biggl(X=\Delta^{0}_YX\hookrightarrow\cdots\hookrightarrow
\Delta^{n}_YX\hookrightarrow\cdots\biggr)\] where
\[\Delta^n_YX:=\int^{K\in\mcF_n}X_Y^K\times\Delta K\cong\int^{K\in\mcF}X_Y^K\times\Delta^{n}K\]
All the maps $\Delta^{n-1}_YX\hookrightarrow\Delta^n_YX$ are
easily seen to be closed inclusions\footnote{the idea is that one
can think of $\Delta^{n-1}_YX\hookrightarrow\Delta^n_YX$ as
$\int^{K\in\mcF_n}X_Y^K\times\Delta^{n-1}K\lra\int^{K\in\mcF_n}X_Y^K\times\Delta^nK$
and show explicitly that it is closed from the fact that the maps
$f^*:X_Y^L\ra X_Y^K$ are closed for $f:K\ra L$ surjective.} and
\begin{equation} \label{sim_filtration_quotient}
\Delta^{(n)}_YX=\Delta^{n}_YX-\Delta^{n-1}_YX\cong
X_Y^{(n+1)}\times_{\Sigma_{n+1}}\textrm{int}\Delta^n
\end{equation}
where $X_Y^{(n+1)}$ denotes the subspace of $X_Y^{n+1}$ of
($n+1$)-tuples of distinct points in $X$ (lying in one fibre over
$Y$) and $\textrm{int}\Delta^n$ denotes the interior of the
standard $n$-dimensional simplex. From now on we assume that all
the spaces admit a proper map $g:X\rightarrow\bbR_+$ (the
nonnegative real numbers). Such spaces are necessarily locally
compact and paracompact (and the converse is true for connected
spaces).

\begin{lemma} \label{sim_support_of_compact_subset}
Let $K\subset\Delta X$ be a compact subset. Then $\supp K$, i.e.
the closure of the set of points in $X$ which appear in some
element of $K$ with a non-zero coefficient, is a compact subset.
\end{lemma}

\begin{proof}
For this proof we introduce a useful tool - an extension map
\[E:\map(X,\bbR)\rightarrow\map(\Delta X,\bbR)\] defined as an
adjoint of
\[E^\sharp:\Delta X\times\map(X,\bbR)\rightarrow\bbR\]
sending $(\sum t_i x_i,g)$ to $\sum t_i g(x_i)$. It can be easily
seen to be continuous.

Let us start with the proof now. Assume that $\supp K$ is not
compact. Then there is a sequence of points $x_n\in X$ such that
$g(x_n)\rightarrow\infty$ and such that $x_n$ appears in some
element of $K$ with a nonzero coefficient $t_n$. Let
$k:\bbR\rightarrow\bbR$ be any map such that $k(g(x_n))>n/t_n$.
Then $E(kg)$ is unbounded on $K$ giving a contradiction%
.
\end{proof}

\begin{corollary}
Suppose that $f:X\rightarrow Y$ is a proper map such that the
number of preimages $|f^{-1}(y)|$ is bounded on every compact
subset of $Y$. Then the induced map $\Delta X\rightarrow\Delta Y$
is also proper.
\end{corollary}

\begin{proof}
Suppose that $K\subset\Delta Y$. Then $K\subset\Delta^{n}Y$ for
some $n$. On a compact subset $\supp K\subset Y$, $|f^{-1}(y)|\leq
k$ and finally $(\Delta f)^{-1}(K)$ is a closed subset of
$\Delta^{kn}(f^{-1}(\supp K))$ which is a compact space.
\end{proof}

Let us restrict ourselves to the case of a
surjective\footnote{When applying the following to a
non-surjective map we just need to replace $Y$ by $\im f$ - as an
image of a proper map it is a closed subspace and so it possesses
a proper map to $\bbR_+$: simply restrict any proper map
$Y\ra\bbR_+$ to $\im f$.} proper map $f:X\ra Y$. As we noted
earlier we have a pullback diagram
\[\xymatrix{
**{!<17pt,0pt>}Rf:=\Delta_YX \ar@{c->}[r] \pullback
\ar[d]_{\hat{f}} & \Delta X \ar[d]^{\Delta f}
\\
Y \ar@{c->}[r] & \Delta Y}\] and for simplicity we denote
$Rf:=\Delta_YX$. It is easy to check that $\hat{f}$, as a pullback
of a proper map, is again proper. Therefore $Rf$ admits a proper
map to $\bbR_+$. Moreover, the fibre of $\hat{f}$ at each point is
a (finite-dimensional) simplex. In particular, by Vietoris-Begle
Theorem \ref{she_vietoris_begle_compact_supports}, $\hat{f}$
induces an isomorphism in \Cech cohomology with compact supports.
As we have seen above, $Rf$ has a natural filtration by
\[R^nf=\Delta^n_YX=\int^{K\in\mcF}X_Y^K\times\Delta^{n}K\]
We then get a spectral sequence with
\[E_1^{pq}=\check{H}_c^{p+q}(R^{p}f,R^{p-1}f)\cong\check{H}_c^{p+q}(R^{p}f-R^{p-1}f)\]
(the isomorphism coming from Proposition II.12.2.~and II.10.2.~of
\cite{Bredon}). When $|f^{-1}(y)|$ is bounded on $Y$ then the
filtration is finite and the spectral sequence converges to
$\check{H}_c^{p+q}(Rf)\cong\check{H}_c^{p+q}(Y)$. We have an
identification (\ref{sim_filtration_quotient})
\[R^{(p)}f=R^{p}f-R^{p-1}f\cong X_Y^{(p+1)}\times_{\Sigma_{p+1}}\textrm{int}\Delta^p\]
and therefore, by the Thom isomorphism,
$E_1^{pq}\cong\check{H}_c^q(X_Y^{[p+1]};\sign)$ where
$X_Y^{[p+1]}=X_Y^{(p+1)}/{\Sigma_{p+1}}$ and $\sign$ denotes the
system of coefficients on $X_Y^{[p+1]}$ given by the composition
\[X_Y^{[p+1]}\lra B\Sigma_{p+1}\xlra{B\sign}B\bbZ/2\]
where the first map represents the principal $\Sigma_{p+1}$-bundle
$X_Y^{(p+1)}\rightarrow X_Y^{[p+1]}$. Therefore we have a result:

\begin{theorem} \label{sim_spectral_sequence}
For a surjective finite-to-one proper map $f:X\rightarrow Y$ we
have a spectral sequence
\[E_1^{pq}=\check{H}_c^q(X_Y^{[p+1]};\sign)\Rightarrow\check{H}_c^{p+q}(Y)\]
which converges if $|f^{-1}(y)|$ is bounded on $Y$. \qed
\end{theorem}

\begin{remark}
There is a different way to define the topology on $Rf$ (when
$Rf\ra Y$ is proper, the case that we considered). Remember our
extension map from the proof of Lemma
\ref{sim_support_of_compact_subset}
\[E^\sharp:\Delta X\times\map(X,\bbR)\rightarrow\bbR\]
and consider its other adjoint
\[\widetilde{E}:\Delta X\rightarrow\map(\map(X,\bbR),\bbR)\]
I do not know whether $\widetilde{E}$ expresses $\Delta X$ as a
subspace of $\map(\map(X,\bbR),\bbR)$ but I doubt it. Nevertheless
for a locally compact Hausdorff (or completely regular for this
matter) space $X$ the restriction of $\widetilde{E}$ to $X$ is a
subspace inclusion
\[X\hookrightarrow\map(\map(X,\bbR),\bbR)\]
Therefore we have a commutative diagram
\[\xymatrix{
Rf \ar[r]^{g} \ar[d]^{\hat{f}} & {}\widetilde{Rf} \ar@{c->}[r]
\ar@{-->}[d]^{\tilde{f}} & {}\map(\map(X,\bbR),\bbR) \ar[d]
\\
Y \ar@{=}[r] & Y \ar@{c->}[r] & {}\map(\map(Y,\bbR),\bbR)}\] where
$\widetilde{Rf}$ denotes the image of $Rf$ in
$\map(\map(X,\bbR),\bbR)$, i.e. the same set $Rf$ but with the
subspace topology. As $\hat{f}$ is proper, so are $\tilde{f}$ and
$g$. Then $g$ is a bijective continuous proper map between locally
compact Hausdorff spaces and is therefore a homeomorphism.
\end{remark}

%



\chapter{A transversality theorem} \label{transversality}

In this chapter we present a transversality theorem for affine
maps from topological affine spaces to sections of an affine
bundle claiming that a certain subset is residual. For this to
have some weight we prove that the affine spaces that we consider
- the topologically finitely generated affine $\mcR$-modules - are
Baire spaces, i.e. residual subsets are dense.

We will be considering maps $\varphi:P\rightarrow C^\infty(M,N)$.
We denote by $\varphi^\sharp$ the adjoint
\[\varphi^\sharp:P\times M\rightarrow N\]
For the following let us denote by $C_\partial^\infty(M,N)$ the
subspace of $C^\infty(M,N)$ of those maps $f$ for which
$f^{-1}(\partial N)=\partial M$.

\begin{lemma} \label{ttt_basic_lemma}
Let $M$, $N$ be smooth manifolds and $A\subseteq N$ a submanifold
without boundary such that either $A\subseteq N-\partial N$ or
$A\subseteq\partial N$. Let there be given two open coverings:
$\scU$ of $M$ and $\scV$ of $A$. Let $P$ be a topological space
and $\varphi:P\rightarrow C_\partial^\infty(M,N)$ a continuous map
where $C_\partial^\infty(M,N)$ is given the weak topology. Assume
that for every $p_0\in P$ and every $U\in\scU$, $V\in\scV$ there
is a finite dimensional manifold $Q$ and a continuous map
$g:Q\rightarrow P$ with $p_0$ in its image such that
\[Q\times U\xlra{g\times incl}P\times M\xlra{\varphi^\sharp}N\]
is smooth and transverse to $V$. Then the set
\[\mfX:=\bigl\{p\in P\ \bigl|\ \varphi(p)\pitchfork A\bigl\}\subseteq P\]
is residual in $P$.
\end{lemma}

\begin{proof}
Following the proof of the Theorem 4.9. of Chapter 4 of
\cite{GolubitskyGuillemin}, let us cover $M$ by a countable family
of compact disks $Y_i$ that have a neighbourhood $U_i\in\scU$ and
at the same time we choose a covering of $A$ by a countable family
of compact subsets $Z_j$ that have a neighbourhood $V_j\in\scV$.
Then the set $\mfX$ is a countable intersection of the sets
\[\mfX_{ij}:=\bigl\{p\in P\ \bigl|\ \varphi(p)\pitchfork Z_j\textrm{ on }Y_i\bigl\}\]
and it is enough to show that each $\mfX_{ij}$ is open and dense.
The set $\hat{\mfX}_{ij}$ of maps $M\rightarrow N$ transverse to
$Z_j$ on $Y_i$ is open in $C^\infty(M,N)$ and
$\mfX_{ij}=\varphi^{-1}(\hat{\mfX}_{ij})$ so it is also open.

To prove the denseness we fix $p_0\in P$ and choose a map
$g:Q\rightarrow P$ with $p_0=g(q_0)$ such that the map
\[h:Q\times U_i\rightarrow N\]
from the statement is smooth and transverse to $V_j$. By the
parametric transversality theorem\footnote{cf.~Theorem 2.7.~of
Chapter 3 of \cite{Hirsch}. For manifolds with boundary
$h^{-1}(V_j)$ is still a submanifold.} the points $q\in Q$ for
which $h(q,-)\pitchfork V_j$ is dense in $Q$. In particular $q_0$
lies in the closure of this set and hence $p_0$ lies in the
closure of its image in $P$. But this image certainly lies in
$\mfX_{ij}$.
\end{proof}

For the next corollary we denote by $\varphi^{(k)}$ the following
map
\[\varphi^{(k)}:P\lra C^\infty(M,N)\lra C^\infty(M^{(k)},N^k)\]
The second map sends $f$ to $f^k|_{M^{(k)}}$ and it is clear that
it is continuous in the weak topology.

\begin{corollary}
Let $M$, $N$ be smooth manifolds and $A\subseteq N^k$ a
submanifold without boundary lying in a single depth $\partial_i
N-\partial_{i+1}N$ of the boundary of $N$. Let there be given two
open coverings: $\scU$ of $M$ and $\scV$ of $A$. Let $P$ be a
topological space and $\varphi:P\rightarrow
C_\partial^\infty(M,N)$ a continuous map where
$C_\partial^\infty(M,N)$ is given the weak topology. Assume that
for every $p_0\in P$ and every $U\in\scU$, $V\in\scV$ there is a
finite dimensional manifold $Q$ and a continuous map
$g:Q\rightarrow P$ with $p_0$ in its image such that
\[Q\times U\xlra{g\times incl}P\times M^{(k)}\xlra{(\varphi^{(k)})^\sharp}N^k\]
is smooth and transverse to $V$. Then the set
\[\mfX:=\bigl\{p\in P\ \bigl|\ \varphi^{(k)}(p)\pitchfork A\bigl\}\subseteq P\]
is residual in $P$. \qed
\end{corollary}

\begin{proposition}
Let $\mcV$ be a topological affine space (over $\bbR$), $M$ a smooth
manifold and $E\rightarrow M$ a smooth finite dimensional affine
bundle. Let $A\subseteq E$ be a submanifold. Let
$\varphi:\mcV\rightarrow\Gamma E$ be a continuous affine map, where
the set $\Gamma E\subseteq C^\infty(M,E)$ of sections of $E$ is
given the weak topology. We denote by $\varphi_x$ the map
\[\varphi^\sharp(-,x):\mcV\rightarrow E_x\]
where $E_x$ is the fibre of $E$ over $x$. Let us assume that for
each $x\in M$ the image $\varphi_x(\mcV)$ is either the whole
fibre $E_x$ or is disjoint with $A$. Then the set
\[\mfX:=\bigl\{v\in\mcV\ \bigl|\ \varphi(v)\pitchfork A\bigl\}\subseteq\mcV\]
is residual in $\mcV$.
\end{proposition}

\begin{proof}
Because of the assumptions we can replace $M$ by its open subset
where $\varphi_x$ is surjective. In effect we can therefore assume
that $\varphi^\sharp$ is surjective.

As $\varphi_x:\mcV\ra E_x$ is a surjective affine map and $E_x$ is
finite dimensional we can find a splitting. Let us denote the
image of such a splitting by $Q_0$ and give it the Euclidean
topology. Then the restricted map
\[Q_0\times M\lra\mcV\times M\xlra{\varphi^\sharp}E\]
is affine and moreover isomorphic on the fibre over $x$. Therefore
there is a coordinate disk $K$ around $x$ such that this map is
isomorphic over $K$. If $v\in V$, then taking the affine hull
$Q=\langle Q_0,v\rangle$ of $Q_0\cup\{v\}$ (again with the
Euclidean topology) we see that the hypothesis of the last lemma
is satisfied.
\end{proof}

The proper setup for the multi-version is the affine bundle
\[E^{(k)}_{\ushort{M}}\lra M^{(k)}\]
which is just the restriction of $E^k\lra M^k$ to $M^{(k)}$. There
is a canonical continuous affine map $\Gamma
E\lra\Gamma\bigl(E^{(k)}_{\ushort{M}}\bigr)$ and we denote the
composition
\[\mcV\xlra{\varphi}\Gamma E\lra\Gamma\bigl(E^{(k)}_{\ushort{M}}\bigr)\]
by $\varphi_k$. As a corollary of the last proposition we get:

\begin{corollary}
Let $\mcV$ be a topological affine space (over $\bbR$), $M$ a smooth
manifold and $E\rightarrow M$ a smooth finite dimensional affine
bundle. Let $A\subseteq E^{(k)}_{\ushort{M}}$ be a submanifold. Let
$\varphi:\mcV\rightarrow\Gamma E$ be a continuous affine map, where
$\Gamma E\subseteq C^\infty(M,E)$ is given the weak topology. Let us
assume that for each $(x_1,\ldots,x_k)\in M^{(k)}$ the image
$(\varphi_{x_1}\times\cdots\times\varphi_{x_k})(\mcV)$ is either the
whole product $E_{x_1}\times\cdots\times E_{x_k}$ of fibres or is
disjoint with $A$. Then the set
\[\mfX:=\bigl\{v\in\mcV\ \bigl|\ \varphi_k(v)\pitchfork A\bigl\}\subseteq\mcV\]
is residual in $\mcV$. \qed
\end{corollary}

\begin{proposition}
Under the assumptions of the last proposition let
$v_i\in\dirsp\mcV$, $i=1,\ldots,n$ be any elements of the
underlying vector space of $\mcV$. Define a map
\[\tilde{\varphi}:\mcV\rightarrow C^\infty(M\times\bbR^n,E)\]
or rather its adjoint
\[\tilde{\varphi}^\sharp:\mcV\times M\times\bbR^n\rightarrow E\]
by the formula
\[\tilde{\varphi}^\sharp(v,x,(\lambda_i))=(\varphi(v+\lambda_1 v_1+\cdots+\lambda_n v_n)(x))\]
Then the set
\[\mfX:=\bigl\{v\in\mcV\ \bigl|\ \tilde{\varphi}(v)\pitchfork A\bigl\}\subseteq\mcV\]
is residual in $\mcV$.
\end{proposition}

\begin{proof}
Let $E\times\bbR^n\rightarrow M\times\bbR^n$ be the product bundle
and define $\hat{\varphi}:\mcV\rightarrow \Gamma(E\times\bbR^n)$
by
\[\hat{\varphi}^\sharp(v,x,(\lambda_i))=(\varphi(v+\lambda_1 v_1+\cdots+\lambda_n v_n)(x),(\lambda_i))\]
Obviously $\mfX=\bigl\{v\in\mcV\ \bigl|\
\hat{\varphi}(v)\pitchfork A\times\bbR^n\bigl\}$ and to apply the
last proposition it is enough to verify that $\hat{\varphi}$ is
continuous. We can express $\hat{\varphi}^\sharp$ as
\[\left(\left(\varphi(v)\times\id\bigr)+\bigl(\varphi(\lambda_1 v_1+\cdots+\lambda_n
v_n)\times\id\right)\right)(x,(\lambda_i))\] and therefore
$\hat{\varphi}$ is a sum of two terms the first of which is the
composition
\[\mcV\xlra{\varphi}\Gamma E\lra\Gamma(E\times\bbR^n)\]
while the second is constant.
\end{proof}

Again we have a multi-version:

\begin{theorem} \label{ttt_main_theorem}
Under the assumptions of the last corollary let
$v_i\in\dirsp\mcV$, $i=1,\ldots,n$ be any elements of the
underlying vector space of $\mcV$. Define a map
\[\tilde{\varphi}_k:\mcV\rightarrow C^\infty\left(M^{(k)}\times\bbR^n,E^{(k)}_{\ushort{M}}\right)\]
or rather its adjoint
\[\tilde{\varphi}_k^\sharp:\mcV\times M^{(k)}\times\bbR^n\rightarrow E^{(k)}_{\ushort{M}}\]
by the formula
\[\tilde{\varphi}_k^\sharp(v,(x_j),(\lambda_i))=(\varphi(v+\lambda_1 v_1+\cdots+\lambda_n v_n)(x_j))\]
Then the set
\[\mfX:=\bigl\{v\in\mcV\ \bigl|\ \tilde{\varphi}_k(v)\pitchfork A\bigl\}\subseteq\mcV\]
is residual in $\mcV$.
\end{theorem}


\begin{proposition}
Let $M$ be a smooth manifold. Then every topologically finitely
generated affine $C^\infty(M,\bbR)$-module is a Baire space, i.e.
every residual subset is dense.
\end{proposition}

\begin{proof}
This is almost a triviality. Let $\mcR=C^\infty(M,\bbR)$ and
$\mcV$ be a topologically finitely generated affine $\mcR$-module.
It is well-known that $\mcR E=C^\infty(M,E)$ is a Baire space. We
have a quotient map
\[p:\mcR E\rightarrow\mcV\]
As it is essentially a quotient by the action of a group - the
kernel of this map - it is an open map. Hence if
$U_i\subseteq\mcV$, $i=1,\ldots$, are open dense subsets then so
are $p^{-1}(U_i)$ and so
\[\mcV=p\left(\overline{\bigcap p^{-1}(U_i)}\right)\subseteq
\overline{p\left(\bigcap p^{-1}(U_i)\right)}=\overline{\bigcap
U_i}\]
\end{proof}

\begin{remark}
This works well for both the strong and the weak topology on
$\mcR=C^\infty(M,\bbR)$ but note that in the strong topology while
$\mcR$ is a topological ring it is not a topological
$\bbR$-algebra unless $M$ is compact.
\end{remark}


\chapter{A spectral sequence} \label{h_principle}

This chapter is the heart of the whole thesis. It is based on the
article \cite{Vassiliev} of Vassiliev. First we have to explain what
our setup is. We have a (topologically finitely generated) affine
$\mcR$-module $\mcV$ and its ``representation'' on an affine bundle
$E$ over $M$ which is a certain map $\varphi:\mcV\ra\Gamma E$ to the
space of smooth sections of $E$.

Typically $E$ is a jet prolongation of some affine bundle $F$ and
$\mcV$ consists either of all sections of $E$ or of all integrable
sections - those that come from sections of the bundle $F$ (so
that in this case one can think of $\mcV$ as $\Gamma F$). The next
piece of data is a ``prohibited'' closed subset $A\subseteq E$.
The main interest of the thesis is in the subspace $\mcV_A$ of
$\mcV$ consisting of those $v\in\mcV$ for which
$\varphi(v)\in\Gamma E$ is a section that does \emph{not}
intersect $A$.

We closely follow \cite{Vassiliev}, \cite{Vassiliev2} and construct
(in this abstract setting) a spectral sequence computing homology of
$\mcV_A$. Here we have to assume that $A$ is a stratified subset of
$E$ of codimension at least $\dim M+1$. The main idea is to
approximate $\mcV$ by finite dimensional affine subspaces $D$ which
``have good transversality and interpolation properties''.

These properties ensure for example that for all $v\in D$ the
section $\varphi(v)$ has only finitely many intersections with
$A$. Let us denote by $X$ the space of pairs $(v,x)$ where $v\in
D$ and $x\in M$ is a point such that the section $\varphi(v)$
meets $A$ at $x$. Then we can reformulate this finiteness
condition as the projection map $X\ra D$ being finite-to-one. This
is where the spectral sequence of Chapter \ref{simplex} comes up.
It converges to the cohomology of the image of this projection
which is clearly $D-\mcV_A$. Alternatively, by Alexander duality,
it converges to the (reduced) \emph{homology} of $D\cap\mcV_A$.
Depending on the ``interpolation quality'' of $D$ one can identify
a range of entries on the $E_1$-page.

With a bit of work one can glue these individual spectral
sequences (for various finite dimensional affine subspaces $D$) to
a single one where one is able to identify the whole $E_1$-page
and it turns out that it does not depend on $\mcV$ too much.
Therefore an affine $\mcR$-homomorphism $\alpha:\mcU\ra\mcV$ is
very likely to be a homology isomorphism provided that both
spectral sequences (for $\mcU$ and $\mcV$) converge. Two criteria
are given at the end of the chapter together with the above
mentioned fundamental example of sections of a jet bundle and the
submodule of integrable sections.

\section{The spectral sequence}

As was said in the introduction the main object of our study is a
representation of a topological $\mcR$-module on a vector bundle
$E\rightarrow M$ and the main example that we have in mind (and
which should therefore satisfy our definition) is the map
$j^r:\Gamma F\rightarrow\Gamma(J^rF)$ sending a section of a
vector bundle $F$ to the corresponding integrable section of the
jet prolongation $E=J^rF$.

Thinking of $E$ as a vector bundle, $\Gamma E$ becomes an
$\mcR$-module \textit{but} the map $j^r$ is not $\mcR$-linear. An
$\mcR$-homomorphism $\varphi:\mcV\rightarrow\Gamma E$ can be
equivalently described as an $\bbR$-linear map that has the
property $\varphi(I_x\mcV)\subseteq I_x\Gamma E$. The map $j^r$
satisfies a weaker condition $j^r((I_x)^r\Gamma F)\subseteq
I_x\Gamma E$. This is the most economical description of what our
notion of a representation should be which is even sufficient for
the construction of the spectral sequence. In next chapters
however we will need more structure associated with
representations. In order to describe it we make the following
observations.

We consider the bundle $J^r(M,\bbR)\rightarrow M$ of algebras
(over the trivial bundle $\bbR\times M\rightarrow M$) and denote
$J^r\mcR=\Gamma(J^r(M,\bbR))$. The inclusion $\bbR\times
M\subseteq J^r(M,\bbR)$ then induces an algebra homomorphism
$i^r:\mcR\rightarrow J^r\mcR$. Every vector bundle $E\rightarrow
M$ has a canonical fibrewise action of $\bbR\times M\rightarrow M$
and this is how the action of $\mcR$ on $\Gamma E$ arises. Suppose
now that there is a fibrewise action of $J^r(M,\bbR)$ on $E$
extending the canonical action of $\bbR\times M$. Then $\Gamma E$
becomes a $J^r\mcR$-module and restricting the action of $J^r\mcR$
along $i^r:\mcR\ra J^r\mcR$ we recover the original $\mcR$-module
structure of $\Gamma E$. On the other hand we also have an algebra
homomorphism $j^r:\mcR\rightarrow J^r\mcR$ sending a section $f$
to its $r$-jet prolongation $j^rf$ and we denote by $j^r\Gamma E$
the $\mcR$-module obtained by restriction along $j^r$. It is
generally different from $\Gamma E$.

%
\begin{definition}
By a \emph{$j^r$-representation} of an $\mcR$-module $\mcV$ on a
vector bundle $E$ over $M$ we mean an action of $J^r(M,\bbR)$ on
$E$ as above together with an $\mcR$-linear map
$\varphi:\mcV\rightarrow j^r\Gamma E$. Note that always
$\varphi((I_x)^r\mcV)\subseteq I_x\Gamma E$ as the action is
fibrewise.

By an \emph{affine $j^r$-representation} of an affine
$\mcR$-module $\mcV$ on an affine bundle $E$ over $M$ we mean an
action of $J^r(M,\bbR)$ on the underlying vector bundle of $E$
together with an affine $\mcR$-homomorphism
$\varphi:\mcV\rightarrow j^r\Gamma E$.
\end{definition}

\begin{definition}
For the purpose of this chapter we use the following definition of
a stratified subset. Let $M$ be a smooth manifold. We say that
$X\subseteq M$ is a stratified subset if there is given a finite
decomposition
\[X=\coprod_{i=0}^kX_i\]
of $X$ into disjoint subsets $X_i$ that are submanifolds of $M$
with either $X_i\subseteq M-\partial M$ or $X_i\subseteq\partial
M$ and such that each partial union $\bigcup_{i=0}^jX_i$ is closed
in $M$, $j=0,\ldots,k$. We call this decomposition a
stratification of $X$.
\end{definition}

The advantage of this definition is that it is closed under
transverse pullbacks, i.e. if $f:N\rightarrow M$ is a smooth map
such that $f\pitchfork X_i$ for all $i=0,\ldots,k$ then
$f^{-1}(X_i)$ provide a stratification of $f^{-1}(X)$.

On the other hand this definition is strong enough to say
something about cohomological properties of $X$. We define the
dimension $\dim X$ of $X$ to be the maximal dimension of its
stratum.

\begin{lemma} $\check{H}_c^n(X;\mcA)=0$ for all $n>\dim X$ and
any coefficient system $\mcA$ on $X$.
\end{lemma}

\begin{proof} We have a spectral sequence
$E_2^{pq}=\check{H}_c^{p+q}(X_p;\mcA|_{X_p})\Rightarrow\check{H}_c^{p+q}(X;\mcA)$
associated with the filtration of $X$ by $\bigcup_{i=0}^jX_i$ and
by our assumptions this spectral sequence vanishes outside of
$0\leq p+q\leq\dim X$. As it converges the same holds for
$\check{H}_c^{p+q}(X;\mcA)$.
\end{proof}

\begin{notation}
In this chapter let $\mcV$ a topologically finitely generated
affine $\mcR$-module, $\xi:E\rightarrow M$ an affine bundle with
$e$ the dimension of the fibre, $\varphi:\mcV\rightarrow j^r\Gamma
E$ an affine $j^r$-representation of $\mcV$ on $E$, $\hat{M}$ the
subset of all points $x\in M$ for which $\varphi_x$ is surjective,
$A\subseteq E$ a closed stratified subset such that outside of
$\hat{M}$ the image of $\varphi_x:\mcV\rightarrow E_x$ is disjoint
with $A$. We denote
\[c=\bigl(\textrm{codimension of $A$ in $E$}\bigr)-m\]
From now on we will assume that $c>0$ or in other words that the
codimension of $A$ in $E$ is at least $m+1$. We are interested in
the space
\[\mcV_A:=\{v\in\mcV\ |\ \im(\varphi(v))\cap A=\emptyset\}\]
here particularly in its homology groups.
\end{notation}

\begin{construction}
We recall the multi-version of $\varphi$
\[\varphi_k:\mcV\xlra{\varphi}\Gamma E\lra\Gamma\bigl(E^{(k)}_{\ushort{M}}\bigr)\]
We will also use its adjoint $\varphi_k^\sharp:\mcV\times
M^{(k)}\rightarrow E^{(k)}_{\ushort{M}}$. We claim that for
$(x_1,\ldots,x_k)\in\hat{M}^{(k)}$ the map
\[(\varphi_k)_{(x_1,\ldots,x_k)}:\mcV\lra\bigl(E^{(k)}_{\ushort{M}}\bigr)_{(x_1,\ldots,x_k)}\cong E_{x_1}\times\cdots\times E_{x_k}\]
is surjective. This is because for any $i\in\{1,\ldots,k\}$ and
for any $w\in E_{x_i}$ there is $v\in\mcV$ with
$\varphi(v)(x_i)=w$; multiplying $v$ by any function $\lambda$
which is $1$ near $x_i$ and $0$ near the remaining points we get
\[\varphi_k(\lambda v)(x_1,\ldots,x_k)=(0,\ldots,0,\overset{\overset{\scriptstyle i}{\downarrow}}{w},0,\ldots,0)\]
We also have a multi-version $A^{(k)}_{\ushort{M}}$ of the
stratified subset $A$ in $E^{(k)}_{\ushort{M}}$, its codimension
is $k(c+m)$ and on the fibres over points outside of
$\hat{M}^{(k)}$ the image of $\varphi_k^\sharp$ is disjoint with
$A^{(k)}_{\ushort{M}}$.
\end{construction}

\begin{proposition} There is a sequence of finite dimensional
affine subspaces $D_d\subseteq\mcV$ satisfying
\begin{itemize}
\item[$(1_d)$]{For each $k$ the map $\varphi_k^\sharp:D_d\times M^{(k)}\rightarrow E^{(k)}_{\ushort{M}}$
is transverse to (each stratum of) $A_k$.}
\item[$(2_d)$]{$D_d\pitchfork\scS_{dr}(\mcV)$.}
\item[$(3)$]{$D_d\subseteq D_{d+1}$ and the union $\bigcup_dD_d$ is dense in
$\mcV$.}
\end{itemize}
\end{proposition}

\begin{proof} First note that conditions $(2_d)$ and $(3)$
depend only on the direction spaces of $D_d$. Theorem
\ref{int_existence_of_transversal} ensures, for each $d$, an
existence of a finite dimensional $D'_d$ satisfying $(2_d)$. It is
also well-known\footnote{It is certainly well-known for free
$\mcR$-modules $\mcR E=C^\infty(M,E)$ of finite rank and this
property clearly passes to quotients.
} that a sequence of finite dimensional $D''_d$ satisfying $(3)$
exists. Therefore
\[D'''_d=D'_1+\cdots+D'_d+D''_d\]
is clearly a sequence of direction spaces satisfying both $(2_d)$
and $(3)$. According to Theorem \ref{ttt_main_theorem} we can find
$v\in\mcV$ such that $D_d=v+D'''_d$ will fulfill even $(1_d)$: the
number of conditions imposed on such a $v$ is countable (one for
each $d$ and $k$) and each such condition is satisfied by a
residual subset.
\end{proof}

Let us try to explain now what these conditions are good for. We
start with the condition $(3)$ which roughly says that $D_d$
approximate $\mcV$ well. This is made more precise in the
following lemma.

\begin{lemma} Whenever $D_d$ is an increasing sequence of affine
subspaces of a locally convex topological affine space $\mcV$ such
that $D_\infty:=\bigcup_dD_d$ is dense in $\mcV$ then for any open
subset $U\subseteq\mcV$ we have the following isomorphisms
\[\colim\ H_*(D_d\cap U)\cong H_*(D_\infty\cap U)\cong H_*(U)\]
\end{lemma}

\begin{proof} The idea of the proof is as follows. Define
$H_*^{\textrm{aff}}(X)$ for any open subspace $X$ of a locally
convex topological vector space in the same way as singular homology
but only allowing affine maps $\Delta^n\rightarrow X$. We prove that
the natural inclusion induces an isomorphism
$H_*^{\textrm{aff}}(X)\cong H_*(X)$. To do so, for a covering $\scU$
of $X$ by open convex subsets, we introduce
$H_*^{\textrm{aff}}(X,\scU)$ where image of each simplex has to lie
in one of the elements of $\scU$. It is a standard fact that the
inclusion induces an isomorphism $H_*(X,\scU)\cong H_*(X)$ and the
same holds for the affine version hence we are left to show that in
the diagram
\[\xymatrix{
H_*^{\textrm{aff}}(X,\scU) \ar[r]^-\cong \ar[d] & H_*^{\textrm{aff}}(X) \ar[d] \\
H_*(X,\scU) \ar[r]^-\cong & H_*(X)}\] the map
$H_*^{\textrm{aff}}(X,\scU)\rightarrow H_*(X,\scU)$ is an
isomorphism. One can easily check that the inverse is induced by
"straightening" the simplices, i.e. by replacing every singular
simplex by the unique affine simplex having the same vertices.
Hence we can replace the singular homology $H_*$ in the statement
by its affine version $H_*^{\textrm{aff}}$. This proves the first
isomorphism. The second one could be proved by an easy observation
that one can perturb slightly (without changing the homology
class) any affine chain in $U$ to one in $D_\infty\cap U$.
\end{proof}

This leads us to investigating $H_*(D_d\cap\mcV_A)$. The
conditions $(1_d)$ and $(2_d)$ that we impose on these spaces will
allow us to construct a spectral sequence for it. Hence let us fix
a finite dimensional affine subspace $D\subseteq\mcV$ satisfying
both $(1_d)$ and $(2_d)$. The proof of the last lemma also shows
that $H_*(D\cap\mcV_A)$, being isomorphic to a purely algebraic
group $H_*^{\textrm{aff}}(D\cap\mcV_A)$, does not depend on the
topology of $D$ as long as this topology makes $D$ into a locally
convex topological vector space in which $\mcV_A$ is open.
Therefore we can and will assume that the topology on $D$ is the
Euclidean topology.

Let us denote by $X\subset D\times M$ the stratified subset of
points $(v,x)$ for which $\varphi_x(v)\in A$. We denote by $\pi$
the composition $X\hookrightarrow D\times M\ra D$. Now we can draw
some consequences of $(1_d)$. First we need a multi-version of the
stratified subset $X$. For every $k$ we have the following diagram
where the square is a pullback square
\begin{equation} \label{hpr_definition_of_DXk}
\entrymodifiers={!!<0pt,\the\fontdimen22\textfont2>+}
\xymatrixnocompile{
& X_D^{(k)} \ar[ld]_{\pi_k} \ar[r] \pullback \ar@{c->}[d] & {}A^{(k)}_{\ushort{M}} \ar@{c->}[d]\\
D & D\times M^{(k)} \ar[r]^-{\varphi_k^\sharp} \ar[l] &
E^{(k)}_{\ushort{M}}}
\end{equation}
(for $k=1$ we are getting $X_D^{(k)}=X$ and $\pi_k=\pi$). The
condition $(1_d)$ then implies that $X_D^{(k)}$ is a
$\Sigma_k$-equivariantly stratified subset of dimension $\dim
X_D^{(k)}=\dim D-ck$. This will be important later, for now we only
need to know that for $k>\dim D$, $X_D^{(k)}=\emptyset$ or in other
words $|\pi^{-1}(v)|\leq\dim D$ for every $v\in D$.
Theorem \ref{sim_spectral_sequence} provides a spectral sequence
\[E_1^{pq}(D)=\check{H}_c^{p+q}(R^{(p)}\pi)\cong\check{H}_c^q(X_D^{[p+1]};\sign)\]
that converges to
\[\check{H}_c^{p+q}(\im(\pi))\cong\widetilde{H}_{\dim D-(p+q)-1}(D\cap\mcV_A)\]
(by Alexander duality). We will now identify
$\check{H}_c^q(X_D^{[p+1]};\sign)$ using the condition $(2_d)$.
Assuming that $p<d$ the map $\varphi_{p+1}^\sharp:D\times
\hat{M}^{(p+1)}\lra E^{(p+1)}_{\ushort{\hat{M}}}$ is fibrewise
surjective. This is because the projection on the fibre over
$(x_0,\ldots,x_p)$ has kernel
\[\varphi^{-1}\bigl(I_{x_0}\cdots I_{x_p}(\Gamma E)\bigr)\supseteq(I_{x_0})^r\cdots(I_{x_p})^r\mcV\]
and
$D\pitchfork(I_{x_0})^r\cdots(I_{x_p})^r\mcV\in\scS_{(p+1)r}(\mcV)$.
Note that then
\[\varphi_{p+1}^\sharp:D\times\hat{M}^{(p+1)}\lra E^{(p+1)}_{\ushort{\hat{M}}}\]
is a fibrewise surjective affine map between affine bundles over
$\hat{M}^{(p+1)}$ and so it is itself an affine bundle. By
restricting this bundle to $A_{p+1}$ we obtain a
$\Sigma_{p+1}$-equivariant affine bundle $X_D^{(p+1)}\rightarrow
A_{\ushort{\hat{M}}}^{(p+1)}$ as in the diagram
(\ref{hpr_definition_of_DXk}) and by further taking
$\Sigma_{p+1}$-orbits an affine bundle
\[h_{p+1}:X_D^{[p+1]}\rightarrow A_{\ushort{\hat{M}}}^{[p+1]}\]
By the Thom isomorphism
\[\check{H}_c^{q}(X_D^{[p+1]};\sign)\cong\check{H}_c^{q-t}\left(A_{\ushort{\hat{M}}}^{[p+1]};w_1(h_{p+1})+\sign\right)\]
where $t=\dim D-e(p+1)$. As we mentioned earlier, this is under
the assumptions that $p<d$.

\begin{note}
One can express $w_1(h_{p+1})$ in terms of $w_1(\xi)$ as follows.
We have the following composition
 \[H_*(\hat{M}^{[p+1]})\xlra{\tr}H_*(\hat{M}^{(p+1)})^{\Sigma_{p+1}}\xlra{pr}H_*(\hat{M})\]
where the first map is a transfer map for the covering
$\hat{M}^{(p+1)}\ra\hat{M}^{[p+1]}$ and the second map is a
projection on a factor (does not matter which one). Applying
$\Hom(-,\bbZ/2)$ we get a map in cohomology
 \[H^1(\hat{M};\bbZ/2)\lra H^1(\hat{M}^{[p+1]};\bbZ/2)\lra
 H^1(A_{\ushort{\hat{M}}}^{[p+1]};\bbZ/2)\]
the second map being induced by the obvious map
$A_{\ushort{\hat{M}}}^{[p+1]}\ra\hat{M}^{[p+1]}$. The image of
$w_1(\xi)$ is $w_1(h_{p+1})+e\cdot\sign$.
\end{note}

Our aim now is to glue the spectral sequences for all $D_d$ into
one and discuss its convergence. First we deal with some
naturality issues for a single spectral sequence as above. Suppose
that $D'\subseteq D$ and that $D'$ also satisfies $(1_d)$ and
$(2_d)$. We have the following diagram of spaces
 \[\xymatrix{
  X'=R^{0}\pi'
   \ar[r] \ar@/^20pt/[rrr]^-{\pi'} \ar@{c->}[d] &
  R^{p}\pi'
   \ar[r] \ar@{c->}[d] &
  R\pi'
   \ar[r] \ar@{c->}[d] &
  D'
   \ar@{c->}[d]
 \\
  X=R^{0}\pi
   \ar[r] \ar@/_20pt/[rrr]_(.4){\pi} &
  R^{p}\pi
   \ar[r] \ar[rrd]_(.7){\rho_p} &
  R\pi
   \ar[r] \ar[rd]^(.6){\rho_\infty} &
  D
   \ar[d]^-\rho
 \\
   & & &
  {}\bbR_+
 }\]
Here $\rho:D\rightarrow\bbR_+$ is any seminorm on $D$ with kernel
$D'$ (equivalently a norm on $D/D'$). Thinking of the middle row
as a diagram of spaces over $\bbR_+$ via $\rho$ the top row
consists of fibres (over $0$) of these spaces.

\begin{proposition} \label{hpr_naturality_of_spectral_sequence}
There is a natural way of defining $i^!_p$ and $i^!_{(p)}$ in such
a way that $i^!_{\infty}=i^!$ and the following diagram commutes
(a morphism of exact couples extending $i^!$ on the colimit)
 \[\xymatrix{
   \check{H}_c^{*-1}(R^{p-1}\pi')
    \ar[d]^{i^!_{p-1}} \ar[r] \ar[d]^{i^!_{p-1}} &
   \check{H}_c^*(R^{(p)}\pi')
    \ar[r] \ar[d]^{i^!_{(p)}} &
   \check{H}_c^*(R^{p}\pi')
    \ar[r] \ar[d]^{i^!_{p}} &
   \check{H}_c^*(R^{p-1}\pi')
    \ar[d]^{i^!_{p-1}}
  \\
   \check{H}_c^{*+\delta-1}(R^{p-1}\pi)
    \ar[r] &
   \check{H}_c^{*+\delta}(R^{(p)}\pi)
    \ar[r] &
   \check{H}_c^{*+\delta}(R^{p}\pi)
    \ar[r] &
   \check{H}_c^{*+\delta}(R^{p-1}\pi)}
 \]
where $\delta=\dim D-\dim D'$. Moreover, if $p<d$ then $i^!_{(p)}$
is the Thom isomorphism.
\end{proposition}

\begin{proof} The idea is as follows. The map
$i^!:\check{H}_c^*(R\pi')\rightarrow\check{H}_c^{*+\delta}(R\pi)$
can be obtained by taking a colimit $\varepsilon\rightarrow 0$ of
\begin{multline*}
\check{H}_c^*(R\pi')\leftarrow\check{H}_c^*(\rho_\infty^{-1}[0,\varepsilon])
\rightarrow\check{H}_c^{*+\delta}(\rho_\infty^{-1}([0,\varepsilon],\{\varepsilon\}))\cong
\\
\cong\check{H}_c^{*+\delta}(\rho_\infty^{-1}(\bbR_+,[\varepsilon,\infty)))
\rightarrow\check{H}_c^{*+\delta}(R\pi)
\end{multline*}
which makes the first arrow an isomorphism. This is proved in the
Appendix \ref{dual_map} in Theorem \ref{dua_main_theorem}. The
second map is a cup product with the pullback along
$R\pi\rightarrow D$ of the generator of
$\check{H}^\delta(\rho^{-1}([0,\varepsilon],\{\varepsilon\})$. We
can now define $i^!_p$ and $i^!_{(p)}$ by the same formula
replacing $\infty$ by $p$ or $(p)$ and the commutativity of the
diagram then follows from various naturality properties of the
maps involved.

Taking a different seminorm one can construct a comparison map
from the above diagram for $\rho$ and $\varepsilon$ to the same
diagram for $\rho'$ and $\varepsilon'$. It is easy to see that in
the colimit this map is an isomorphism - its inverse is
constructed in exactly the same way. Therefore our construction
does not depend on the seminorm $\rho$.

The last part follows from the fact that in this case
$R^{(p)}\pi'\hookrightarrow R^{(p)}\pi$ is an inclusion of (the
image of) a section of an affine bundle and in this case the above
constructed map is precisely the Thom isomorphism.
\end{proof}

This ensures that we get a morphism of the corresponding spectral
sequences. First we have to rewrite the spectral sequence in a
homological way by transformation $p\leftrightarrow -p-1$ and
$q\leftrightarrow\dim D-q$
\[E^1_{pq}(D)=\check{H}_c^{\dim D-(p+q)-1}(R^{(-p-1)}\pi)\cong\check{H}_c^{\dim D-q}(X_D^{[-p]};\sign)\]
converging to $\widetilde{H}_{p+q}(D\cap\mcV_A)$. Stably (when
$p>-d-1$) we have
\[E^1_{pq}(D)\cong\check{H}_c^{-ep-q}(A_{\ushort{\hat{M}}}^{[-p]};w_1(h_{-p})+\sign)\]
Now we can define a colimit spectral sequence
\begin{equation} \label{hpr_colimit_spectral_sequence}
E^1_{pq}(\mcV)=\colim_d\check{H}_c^{\dim
D_d-(p+q)-1}(R^{(-p-1)}\pi_d)\Rightarrow\widetilde{H}_{p+q}(\mcV_A)
\end{equation}
using Proposition \ref{hpr_naturality_of_spectral_sequence}. Also
from this proposition it is clear that
\[E^1_{pq}(\mcV)\cong\check{H}_c^{-ep-q}(A_{\ushort{\hat{M}}}^{[-p]};w_1(h_{-p})+\sign)\]
for all $p$ and $q$. As we discussed before $E^1_{pq}(D)$ always
converges. Now we want to give some sufficient conditions for the
convergence of
$E^1_{pq}(\mcV)\Rightarrow\widetilde{H}_{p+q}(\mcV_A)$.

\begin{proposition}
\label{hpr_zero_spectral_sequence_implies_convergence}
If $E^1_{pq}(\mcV)=0$ then the spectral sequence converges. In
other words in this case $\widetilde{H}_*(\mcV_A)=0$.
\end{proposition}

\begin{proof} Up to a change in the indices
$\widetilde{H}_*(D_d\cap\mcV_A)$ is identified with
$\check{H}_c^*(R\pi_d)$ and for $d'\gg d$ the latter is mapped to
the stable part of $\check{H}_c^*(R\pi_{d'})$ which by our
assumptions is $0$. In other words the map
$\widetilde{H}_*(D_d\cap\mcV_A)\rightarrow\widetilde{H}_*(D_{d'}\cap\mcV_A)$
is $0$ for $d'\gg d$. As
$\widetilde{H}_*(\mcV_A)\cong\colim\widetilde{H}_*(D_d\cap\mcV_A)$
this clearly proves the claim.
\end{proof}

\begin{proposition} If $c>1$, or in other words if the
codimension of $A$ in $E$ is at least $m+2$, then the spectral
sequence $E^1_{pq}(\mcV)\Rightarrow\widetilde{H}_{p+q}(\mcV)$
converges.
\end{proposition}

\begin{proof} The proof is similar. Clearly
$E^1_{pq}(D_d)\cong\check{H}_c^{\dim D_d-q}(X_{D_d}^{[-p]};\sign)$
is $0$ for $p>0$. Hence we assume $p\leq 0$ and we recall that
$X_{D_d}^{(-p)}$ is a $\Sigma_{-p}$-equivariantly stratified
subset of $D_d\times M^{(-p)}$ of dimension at most $\dim D_d+cp$.
Therefore we also get $E^1_{pq}(D_d)=0$ for $q<-cp$. Consequently
for a fixed $n$ and for $d\gg 0$, $r\gg 0$ both
$\widetilde{H}_n(D_d)$ and
$\Tot\bigl(E^r_{**}(D_d)\bigr)_n=\bigoplus_{p+q=n}E^r_{pq}(D_d)$
do not depend on $d$ and $r$. One concludes that the colimit
spectral sequence converges from the convergence of each
$E^1_{pq}(D_d)$.
\end{proof}

\begin{theorem} Suppose that $c>1$ and that
$\alpha:\mcU\subseteq\mcV$ is a topological submodule such that
$(\varphi\alpha)_x$ is still surjective for all $x\in\hat{M}$.
Then the map $\alpha_A:\mcU_A\rightarrow\mcV_A$ is a homology
isomorphism.
\end{theorem}

\begin{proof} We have a map of spectral sequences\footnote{We are
cheating here a bit. But it is not difficult to arrange the
subspaces $D_d$ in $\mcU$ and $\mcV$ in such a way that $\alpha$
preserves them. Also later we prove independently a stronger
result.} $E^r_{pq}(\mcU)\rightarrow E^r_{pq}(\mcV)$ and on the
$E^1$-page this map is an isomorphism. Hence so is the map on the
limits.
\end{proof}


\section{Examples}

The following is a fundamental example.

\begin{example} Let $F\rightarrow M$ be a finite
dimensional affine bundle and set $E=J^rF$. Then there is an obvious
action of $J^r(M,\bbR)$ on $E$ and we have a topological
$\mcR$-module $\mcV=j^r\Gamma E$ (which is topologically finitely
generated by Example
\ref{int_finite_generation_of_sections_of_jet_bundle}). Obviously
$\mcU=\Gamma F$ can be identified via $j^r$ with an $\mcR$-submodule
of $\mcV$ of integrable sections. Clearly
$\varphi=\id:\mcV\rightarrow\Gamma E$ is surjective on fibres as is
$\varphi j^r$. Therefore by the last theorem we know that if $c>1$
then the inclusion map $(\Gamma F)_A\hookrightarrow(\Gamma(J^rF))_A$
is a homology isomorphism (this will be improved to having acyclic
fibres in Theorem~\ref{thm_main_theorem}).
\end{example}

\begin{example}
To illustrate what the spaces appearing in the spectral sequence
$E^1_{pq}=E^1_{pq}(\mcV)$ are let us take $\mcV=\Gamma E$, the space
of sections of a vector bundle $\xi:E\ra M$ of dimension at least
$\dim M+2$ and let $A=M$ be the zero section. Then with respect to
the representation $\id:\mcV\ra\Gamma E$ we have $\hat{M}=M$ and
 \[E^1_{pq}\cong\check{H}_c^{-ep-q}(A_{\ushort{\hat{M}}}^{[-p]};w_1(h_{-p})+\sign)\cong
 H_c^{-ep-q}(M^{[-p]};(pr\cdot\tr)^*w_1(\xi)+(e+1)\sign)\]
converging to the reduced homology of the space of sections of the
associated unit sphere bundle.
\end{example}



\chapter{Constructing the homotopy fibre} \label{homotopy_fibre}

In the previous chapter we saw that an affine $\mcR$-homomorphism
$\alpha:\mcU\ra\mcV$ induces a homology isomorphism
$\alpha_A:\mcU_A\ra\mcV_A$ if it is ``enough surjective''. In this
chapter we make the first step towards the proof of a stronger
result that all the homotopy fibres of $\alpha_A$ are acyclic
(have zero homology). This step consists of constructing the
homotopy fibres in our realm of $j^r$-representations of
topologically finitely generated affine modules. A simple
computation in the next chapter then shows that the first page of
the basic spectral sequence (\ref{hpr_colimit_spectral_sequence})
is zero.

In the first section we construct the underlying module and show
that it is topologically finitely generated. In the next section
the $j^r$-representation is described while the biggest part is
devoted to the proof that what we get is indeed the homotopy
fibre.

The actual construction is best understood on the example of the
space $\mcV=\Gamma(F)$ of sections of an affine bundle $F$ and its
submodule $\mcU$. Classically the homotopy fibre would be a
certain subspace of the space $\map(I,\mcV)$ of continuous paths
in $\mcV$. Remembering that $\mcV=\Gamma(F)$ we can identify the
path space $\map(I,\mcV)$ with some subspace of the space
$\Gamma^0(F\times I)$ of continuous sections of $F\times I\ra
M\times I$. We follow the same idea but take only smooth sections.
This corresponds in abstract to an extension of scalars. We have a
topological algebra $\mcR^I=C^\infty(M\times I,\bbR)$ and an
inclusion map $\mcR\hookrightarrow\mcR^I$. Our space of sections
of $F\times I$ is then nothing but $\mcR^I\tensor{\mcR}\mcV$.

The topological tensor product is not as well-behaved as is the
algebraic one in general. Its basic properties on the category of
locally convex topological vector spaces are gathered in Appendix
\ref{LCTVS}. In Appendix \ref{function_spaces} a proof that $\mcR$
is locally convex can be found. As locally convex topological
vector spaces are closed under products and quotients every
topologically finitely generated $\mcR$-module is also locally
convex. Therefore we do not lose anything in assuming that all
$\mcR$-modules in this chapter are locally convex.

\section{The algebraic part}

We consider the smooth manifold $M\times I$ and set
$\mcR^I=C^\infty(M\times I,\bbR)$. There is a natural inclusion
$\varepsilon:\mcR\subseteq\mcR^I$ induced by the projection
$M\times I\rightarrow M$. For any $t\in I$ one gets a splitting
$res_t:\mcR^I\rightarrow\mcR$ of the above inclusion
$\mcR\subseteq\mcR^I$ induced by
\[M\cong M\times\{t\}\subseteq M\times I\]

Therefore every topological $\mcR^I$-module is canonically a
topological $\mcR$-module and for every topological $\mcR$-module
$\mcV$ we have an induced topological $\mcR^I$-module
$\mcV^I:=\mcR^I\tensor{\mcR}\mcV$ together with a split inclusion
$\varepsilon:\mcV\hookrightarrow\mcV^I$. Also if $\mcV$ is
topologically finitely generated then so is $\mcV^I$. Indeed,
tensoring the quotient map $\mcR E\rightarrow\mcV$ with $\mcR^I$
yields a quotient map
\[\mcR^I E\cong \mcR^I\tensor{\mcR}\mcR
E\rightarrow\mcR^I\tensor{\mcR}\mcV=\mcV^I\]

\begin{lemma} \label{mod_splitting_of_RI}
If we denote the $I$-coordinate function by $t$ then there is a
split short exact sequence
\[0\longrightarrow \mcR^I\xlra{t(t-1)\cdot}\mcR^I\xlra{(res_0,res_1)^T}\mcR\times\mcR\longrightarrow0\]
of $\mcR$-modules.
\end{lemma}

\begin{note} As the category of topological $\mcR$-modules (or
topological vector spaces) is not an abelian category we have to
be a little careful here:
\[0\longrightarrow A\xlra{i}B\xlra{q}C\longrightarrow0\]
is a short exact sequence if $A$ is the kernel (in the sense of
topological $\mcR$-modules) of $q$ and $C$ the cokernel of $i$. In
other words if it is a short exact sequence of $\mcR$-modules such
that $A$ has the subspace topology and $C$ the quotient topology
induced from $B$. The formulation of a split short exact sequence
we are going to use here is that of a commutative diagram
\[\xymatrix{
A \ar[rr]^{0} \ar[dr]^{i} \ar[dd]_{\id} & & C \\
& B \ar[ur]^{q} \ar[dl]^{p} \\
A & & C \ar[ll]^{0} \ar[ul]^{j} \ar[uu]_{\id} }\] such that
$ip+jq=\id$. In such a situation we say that $B$ is a biproduct of
$A$ and $C$.
\end{note}

\begin{proof} The second map has a section sending $(f_0,f_1)$
to $f(x,t)=(1-t)f_0(x)+tf_1(x)$. Now we construct a section of the
first map. For $(x,t)\in M\times I$ we define a $2$-simplex
$\sigma_{(x,t)}:\Delta^2\rightarrow M\times I$ to be
\[\sigma_{(x,t)}=[(x,0),(x,t),(x,1)]\]
and then for $f\in\mcR^I$ we set
\[\mcJ(f,x,t):=\frac{1}{2}\cdot
\int_{\Delta^2}f^{''}_{tt}\sigma_{(x,t)}\] The map
$p:\mcR^I\rightarrow\mcR^I$ defined by $p(f)(x,t)=\mcJ(f,x,t)$ is
continuous and by Lemma \ref{int_integration_formula} has the
property that
\[t(t-1)p(f)(x,t)=f(x,t)-(1-t)f(x,0)-tf(x,1)\]
Therefore it is a splitting of the
``multiplication by $t(t-1)$'' map (as this multiplication map is
clearly injective). The same formula proves that these splittings
describe $\mcR^I$ as a biproduct of $\mcR^I$ and $\mcR\times\mcR$.
\end{proof}

\begin{corollary} \label{mod_finite_generation_of_homotopy_fibre}
Let $\mcV,\mcV_0,\mcV_1$ be topologically finitely generated
$\mcR$-modules, let $\alpha_i:\mcV_i\rightarrow\mcV$ be
$\mcR$-homomorphisms. Then the $\mcR^I$-module $\mcW$ defined as
the pullback
\[\xymatrix{
{\mcW} \ar[rr] \ar[d] & & {\mcV_0\times\mcV_1}
\ar[d]^{\alpha_0\times\alpha_1}
\\
{\mcV^I} \ar[rr]^{(res_0,res_1)^T} & & {\mcV\times\mcV}}\] is also
topologically finitely generated. Here the $\mcR^I$-module
structures on $\mcV\times\mcV$ and $\mcV_0\times\mcV_1$ are
induced by $\mcR^I\xlra{(res_0,res_1)}\mcR\times\mcR$.
\end{corollary}

\begin{proof} Starting with the split short exact sequence from
the previous lemma we tensor it with $\mcV$ to get
\[0\longrightarrow\mcV^I\xlra{t(t-1)\cdot}\mcV^I
\xlra{(res_0,res_1)}\mcV\times\mcV\longrightarrow0\] Taking the
pullback along $\alpha_0\times\alpha_1$ we get a split short exact
sequence of $\mcR$-modules
\[0\longrightarrow\mcV^I\xlra{i}\mcW
\longrightarrow\mcV_0\times\mcV_1\longrightarrow0\] Hence if
$j:\mcV_0\times\mcV_1\rightarrow\mcW$ is the splitting and
$\varphi:\mcR E\rightarrow\mcV$ and $\varphi_i:\mcR
E_i\rightarrow\mcV_i$ are topological presentations then the
following is a topological presentation
\[(i\varphi^I,\overline{j\varphi_0},\overline{j\varphi_1}):
\mcR^I E\times\mcR^I E_0\times\mcR^I E_1\longrightarrow\mcW\] of
$\mcW$
. Here $\overline{j\varphi_i}$ is the unique $\mcR^I$-linear
extension of the $\mcR$-homomorphsisms $j\varphi_i:\mcR
E_i\rightarrow\mcW$ to $(\mcR E_i)^I\cong\mcR^I E_i$.
\end{proof}

\begin{note}
The same is true for affine $\mcR$-modules and affine maps between
them. The pullback is clearly nonempty and so we can reduce to the
linear case by choosing a common origin (take the images of any
element in the pullback). As the forgetful functor
$U:\Rmod\ra\AffRmod$ preserves limits\footnote{There is an
adjunction $F:\AffRmod\rightleftarrows\Rmod:U$ where $F$ sends an
affine $\mcR$-module $\mcV$ to an affine $\mcR$-module spanned
freely by $\mcV$ and an extra point $0$ serving as an origin (this
is a general fact about comma categories - the forgetful functor
$c\downarrow\mcC\rightarrow\mcC$ has a left adjoint
$d\mapsto(c\hookrightarrow d\coprod c)$).}  the affine case then
follows from the linear version.
\end{note}

\section{The topological part}

Suppose that $\varphi:\mcV\rightarrow\Gamma E$ is a
$j^r$-representation. We construct a $j^r$-representation of $\mcV^I$ on
$E\times I$ as follows. There is a canonical homomorphism
\[J^r(M\times I,\bbR)\rightarrow J^r(M,\bbR)\times I\]
of bundles (of algebras) over $M\times I$ which simply forgets
some of the jet information. On the other hand we have a product
action
\[(J^r(M,\bbR)\times I)\underset{M\times I}{\times}(E\times I)\rightarrow E\times I\]
Composing these two maps gives the desired action of $J^r(M\times
I,\bbR)$ on $E\times I$. There is an obvious map
$\varepsilon:\Gamma E\rightarrow\Gamma(E\times I)$ (sending a
section $f:M\rightarrow E$ to the section $f\times\id:M\times
I\rightarrow E\times I$ constant in the $I$ direction). This map
is $\mcR$-linear with respect to the $j^r$-structures. Therefore
there is a unique $\mcR^I$-linear map $(j^r\Gamma E)^I\rightarrow
j^r\Gamma(E\times I)$ giving a factorization $j^r\Gamma
E\rightarrow(j^r\Gamma E)^I\rightarrow j^r\Gamma(E\times I)$ of
$\varepsilon$ and we define $\varphi^I$ to be the composition
$\mcV^I\rightarrow(j^r\Gamma E)^I\rightarrow\Gamma(E\times I)$.
\[\xymatrix{
{\mcV} \ar[r]^{\varepsilon} \ar[d]_{\varphi} & {\mcV^I} \ar[r]^{=} \ar[d] & {\mcV^I} \ar[d]^{\varphi^I} \\
{\Gamma E} \ar[r] & {(j^r\Gamma E)^I} \ar[r] & {\Gamma(E\times
I)}}\]

Now we need few auxiliary constructions. Proving that they are all
given by continuous maps would destroy the flow of this section
and so these proofs are only given in Appendix
\ref{function_spaces}.

Let $\lambda:(I,\partial I)\rightarrow (I,\partial I)$ be a smooth
function. We define a reparametrization map
\[\lambda^*:\mcR^I\rightarrow\mcR^I\]
by $f(x,t)\mapsto f(x,\lambda(t))$. It is easily seen to be a
continuous $\mcR$-homomorphism. Moreover thinking of $\mcR^I$ as a
bundle over $\mcR\times\mcR$ via $(res_0,res_1)$ it is a fibrewise
map (not necessarily over $\id$ but there is only a discrete
choice of base maps). If $\lambda_0,\lambda_1$ are two such maps
with $\lambda_0=\lambda_1$ on $\partial I$ then the maps
$\lambda_0^*$ and $\lambda_1^*$ are homotopic through
reparametrization maps and in particular by a fibrewise homotopy.
For an $\mcR$-module $\mcV$ we define a reparametrization map
\[\lambda^*:\mcV^I=\mcR^I\tensor{\mcR}\mcV\xlra{\lambda^*\otimes\id}\mcR^I\tensor{\mcR}\mcV=\mcV^I\]
and observe that for any $j^r$-representation of $\mcV$ the
diagram
\[\xymatrix{
{\mcV^I} \ar[r] \ar[d]_{\lambda^*} & {\Gamma(E\times I)} \ar[d]^{\lambda^*} \\
{\mcV^I} \ar[r] & {\Gamma(E\times I)}}\] commutes. Let us denote
by $(\mcR^I)^{\hat{n}}$ the following subset of $(\mcR^I)^n$
\[(\mcR^I)^{\hat{n}}=\bigl\{(f_1,\ldots,f_n)\in(\mcR^I)^n\ \bigl|\ res_1f_{i}=res_0f_{i+1}, i=1,\ldots,n-1\bigr\}\]
i.e. the space of $n$-tuples of functions that ``can be
concatenated''. The problem is of course that one cannot
concatenate them straight away, the result would not be smooth. We
have a projection
\[(\sigma,\tau):(\mcR^I)^{\hat{n}}\rightarrow\mcR\times\mcR\]
sending $(f_1,\ldots,f_n)$ to $(res_0f_1,res_1f_n)$ (here $\sigma$
stands for the source map and $\tau$ for the target map). We think
of $(\mcR^I)^{\hat{n}}$ as a space over $\mcR\times\mcR$ via this
map. Let $0=s_0<\cdots<s_n=1$ be a sequence of numbers and
$\lambda_1,\ldots,\lambda_n:I\rightarrow I$ smooth nondecreasing
maps with $\lambda_i=0$ near $0$ and $\lambda_i=1$ near $1$ (but
where we only require $\lambda_1(0)=0$ and $\lambda_n(1)=1$ to
allow $\id=\mu_1((0,1),\id)$). We define a concatenation map
\[\mu_n=\mu_n((s_i),(\lambda_i)):(\mcR^I)^{\hat{n}}\longrightarrow\mcR^I\]
by sending $(f_i)$ to the concatenation of the maps
$\lambda_i^*f_i$ shrunk to $M\times[s_{i-1},s_i]$. It is clearly
an $\mcR$-linear map over $\mcR\times\mcR$. We also have the
canonical inclusion $\varepsilon:\mcR\rightarrow\mcR^I$ and a
reverse map
\[\iota:\mcR^I\xlra{(1-t)^*}\mcR^I\]
(which is not a map over $\mcR\times\mcR$). All the maps
$\mu_n((s_i),(\lambda_i))$ (for all possible $s_i$, $\lambda_i$
but with $n$ fixed) are homotopic through maps of the form $\mu_n$
(and in particular over $\mcR\times\mcR$). One easily verifies
that $\mu_n(\mu_{k_1}\times\cdots\times\mu_{k_n})$ is of the form
$\mu_{k_1+\cdots+k_n}$; $\mu_{n+1}(\varepsilon\sigma,\id)$ and
$\mu_{n+1}(\id,\varepsilon\tau)$ are of the form $\mu_n$. In
particular $\mu_2(\id\times\mu_2)\simeq\mu_2(\mu_2\times\id)$
through maps of the form $\mu_3$ and similarly
$\mu_2(\varepsilon\sigma,\id)\simeq\id$,
$\mu_2(\id,\varepsilon\tau)\simeq\id$. Also
\[\mu_2(f,\iota f)=(\mu_2(t,1-t))^*f\]
and thus $\mu_2(\id,\iota)\simeq\sigma$ over $\mcR\times\mcR$ and
through reparametrization maps, similarly $\mu_2(\iota,\id)\simeq
\tau$.

Now we extend these maps to $\mcV^I$ by taking the tensor product
with $\mcV$. We identify $(\mcR^I)^{\hat{n}}\tensor{\mcR}\mcV$
with the obvious generalization $(\mcV^I)^{\hat{n}}$ of
$(\mcR^I)^{\hat{n}}$ using Lemma \ref{mod_splitting_of_RI}. We get
maps $\mu_n:(\mcV^I)^{\hat{n}}\longrightarrow\mcV^I$,
$\varepsilon:\mcV\rightarrow\mcV^I$ and
$\iota:\mcV^I\rightarrow\mcV^I$. They commute with the
representation $\varphi^I:\mcV^I\rightarrow\Gamma(E\times I)$. The
crucial property of these maps on $\Gamma(E\times I)$ is that for
any subset $A\subseteq E$ these maps preserve
\[\Gamma(E\times I)_{A\times I}=\{f\in\Gamma(E\times I)\ |\ \im f\cap (A\times I)=\emptyset\}\]
where for $\varepsilon$ this is understood as $\varepsilon(\Gamma
E_A)\subseteq\Gamma(E\times I)_{A\times I}$. Therefore the maps on
$\mcV^I$ preserve
\[(\mcV^I)_{A\times I}=(\varphi^I)^{-1}\Gamma(E\times I)_{A\times I}\]

\begin{theorem} \label{mod_fibration_theorem}
The following map is a Serre fibration if $A$ is closed
\[p:(\mcV^I)_{A\times I}\xlra{(res_0,res_1)}\mcV_A\times\mcV_A\]
\end{theorem}

\begin{proof} The map
\[\mcV^I\xlra{(res_0,res_1)}\mcV\times\mcV\]
is a projection on a direct summand by Lemma
\ref{mod_splitting_of_RI}. As both $(\mcV^I)_{A\times I}$ and
$\mcV_A\times\mcV_A$ are open subsets (when $A$ is closed) we see
that for every $u_0\in\mcV_A$ there is a neighbourhood $U$ of
$u_0$ and a section
\[j_{u_0}:U\cong \{u_0\}\times U\rightarrow(\mcV^I)_{A\times I}\]
with $j_{u_0}(u_0)=\varepsilon(u)$. Therefore if
$(u_0,v_0)\in\mcV_A\times\mcV_A$ and $U$, $V$ are the
neighbourhoods as above we can define two maps over $U\times V$
\begin{eqnarray*}
p^{-1}(U\times V) & \rightarrow & U\times V\times p^{-1}(u_0,v_0) \\
U\times V\times p^{-1}(u_0,v_0) & \rightarrow & p^{-1}(U\times V)
\end{eqnarray*}
The first map sends $w$ to
\[(res_0w,\ res_1w,\ j_{u_0}(res_0w)*w*\iota (j_{v_0}(res_1w)))\]
and the second sends $(u,v,w)$ to $\iota(j_{u_0}(u))*w*j_{v_0}(v)$
where we use the notation with $*$ to denote the multiplication
$\mu_n$.

They are easily seen to be homotopy inverse to each other over
$U\times V$. Therefore $p$ has the weak homotopy lifting property.
As an open subset of a fibration it is also a microfibration.
These two properties together easily imply that it is a Serre
fibration\footnote{For a map $f:D^k\times I\ra\mcV_A\times\mcV_A$
and a partial lift $L_0$ over $D^k\times 0$ there is a fibrewise
homotopy $h$ from $L_0$ to a lift that can be extended to $L_1$
defined on $D^k\times I$. We think of $h$ and $L_1$ as giving a
partial lift of $D^k\times I\times I\ra D^k\times
I\xra{f}\mcV_a\times\mcV_A$ over $D^k\times 0\times I\cup
D^k\times I\times 1$. By a microfibration property there is a
neighbourhood $D^k\times [0,\varepsilon]\times I\cup D^k\times
I\times[1-\varepsilon,1]$ and a lift $L_2$ over it extending both
$h$ and $L_1$. We define $L:D^k\times I\ra(\mcV^I)_{A\times I}$ to
be $L_2(x,t,t/\varepsilon)$ for $t\leq\varepsilon$ and
$L_2(x,t,1)$ for $t\geq\varepsilon$.}.
\end{proof}

\begin{note}
We have an affine version of all the constructions. To define
\[\mcV^I=\mcR^I\tensor{\mcR}\mcV\]
we just choose an origin in $\mcV$, take the ordinary tensor
product and then forget the origin of the result. The canonical
inclusion $\varepsilon:\mcV\ra\mcV^I$ as well as the maps $\mu_n$,
$\iota$ do not depend on this choice. The same is true for the
induced $j^r$-representation
\[\varphi^I:\mcV^I\ra\Gamma(E\times I)\]
of an affine $j^r$-representation $\varphi:\mcV\ra\Gamma E$: the
choice of an origin in $\mcV$ gives a section in $E$ which then
makes $E$ into a vector bundle. After extending to $M\times I$ we
forget the origin/section and the result is independent of our
choice. Clearly Theorem \ref{mod_fibration_theorem} remains true
in this case.
\end{note}

Let $\alpha:\mcU\rightarrow\mcV$ be an affine $\mcR$-homomorphism,
$\varphi:\mcV\rightarrow\Gamma E$ an affine $j^r$-representation,
$A\subseteq E$ a closed subset. By
$\alpha_A:\mcU_A\rightarrow\mcV_A$ we denote the restriction of
$\alpha$ to the open subset $\mcU_A$. As this is defined via the
representation $\varphi\alpha:\mcU\rightarrow\Gamma E$ its image
necessarily lies in $\mcV_A$.

\begin{theorem} \label{mod_homotopy_fibre}
In the notation from above let $v\in\mcV_A$ be a point and let
$\mcW$ be the pullback
\[\xymatrix{
{}\mcW \ar[rr] \ar[d] & & {}\{v\}\times\mcU \ar[d]^{incl\times\alpha} \\
{}\mcV^I \ar[rr]^{(res_0,res_1)^T} & & {}\mcV\times\mcV}\] Then
the homotopy fibre $\hofib_v\alpha_A$ of
$\alpha_A:\mcU_A\rightarrow\mcV_A$ over $v$ has the weak homotopy
type of $\mcW_{A\times I}$ defined in terms of the affine
$j^r$-representation
$\psi:\mcW\ra\mcV^I\xra{\varphi^I}\Gamma(E\times I)$.
\end{theorem}

\begin{proof} We have a diagram where all the squares are
pullback squares
\[\xymatrix{
& & {\mcW} \ar[rr] \pullback \ar[d] & & \{v\}\times\mcU_A \ar[rr]
\pullback \ar@{c->}[d] & & \{v\} \ar@{c->}[d]
\\
{\mcU_A} \ar[rr]^-{\simeq} & & {\mcP} \ar@{->>}[rr]^F \pullback
\ar[d]_G & & {\mcV_A\times\mcU_A}
\ar@{->>}[rr]^-{\textrm{projection}} \ar[d]^{\id\times\alpha_A} &
& {\mcV_A}
\\
& & (\mcV^I)_{A\times I} \ar@{->>}[rr]^{(res_0,res_1)^T} & &
{\mcV_A\times\mcV_A}}\] The maps denoted by
$\xymatrix@1{{}\ar@{->>}[r]&{}}$ are fibrations and the map
denoted by $\simeq$ is a homotopy equivalence constructed from the
universal property of a pullback from maps
$\mcU_A\xlra{(\alpha_A,\id)}\mcV_A\times\mcU_A$ and
$\mcU_A\xlra{\varepsilon\alpha_A}(\mcV^I)_{A\times I}$ (there is a
deformation retraction
 \[(\lambda_s^*G,(res_sG,pr_{\mcU_A}F)):\mcP\lra\mcP=((\mcV^I)_{A\times I})\times_{(\mcV_A\times\mcV_A)}(\mcV_A\times\mcU_A)\]
whose main part is $\lambda_s^*:(\mcV^I)_{A\times
I}\rightarrow(\mcV^I)_{A\times I}$, $s\in I$, with
$\lambda_s(t)=t+(1-t)s$). As the composition across the middle row
is $\alpha_A$ we see that $\mcW$ is the homotopy fibre
$\hofib_v\alpha_A$.
\end{proof}


\begin{lemma} \label{mod_image_of_the_pullback_representation}
For the affine $j^r$-representation $\psi$ from the last theorem
 \[\psi_{(x,t)}(\mcW)=\left\{\begin{array}{rcl}
  \varphi_x(v) & \textrm{for} & t=0 \\
  \varphi_x(\mcV) & \textrm{for} & 0<t<1 \\
  (\varphi\alpha)_x(\mcU) & \textrm{for} & t=1
 \end{array}\right.\]
\end{lemma}

\begin{proof}
We have a diagram
 \[\xymatrix{
  & &
  {}\mcW \ar[rr] \ar[d] & &
  {}\{v\}\times\mcU \ar[d]^{incl\times\alpha}
 \\
  {}\mcV \ar[d]^{\varphi} & &
  {}\mcV^I \ar[rr]^{(res_0,res_1)^T} \ar[ll]_{res_t} \ar[d]^{\varphi^I} & &
  {}\mcV\times\mcV \ar[d]^{\varphi\times\varphi}
 \\
  \Gamma E & &
  {}\Gamma(E\times I) \ar[rr]^{(res_0,res_1)^T} \ar[ll]_{res_t} & &
  {}\Gamma E\times\Gamma E
 }\]
This proves easily the claims for $t=0$ and $t=1$. Also it is
clear that for $0<t<1$ the image can be at most
$\varphi_x(\mcV)=(\varphi^I)_{(x,t)}(\mcV^I)$. But according to
the proof of Corollary
\ref{mod_finite_generation_of_homotopy_fibre} we have a
decomposition $\overline{\mcW}\cong\overline{\mcV^I}\times
\overline{\{v\}}\times\overline{\mcU}$ of the underlying
$\mcR$-modules and the linear part of the composition
$\mcV\subseteq\mcV^I\xra{\psi}\Gamma(E\times I)$ maps $v$ to a
section
 \[(x,t)\mapsto t(t-1)\varphi_x(v)\]
and as $t(t-1)\neq 0$ the image is the same as that of
$\varphi_x$.
\end{proof}



\chapter{The main theorem} \label{main_theorem}

This short chapter is devoted to outlining the relation between
cohomology of a configuration space of a product of two spaces and
cohomology of the two individual configuration spaces. The answer is
not at all satisfactory but it suffices for proving (under some
assumptions) the acyclicity of the homotopy fibre of
$\alpha_A:\mcU_A\lra\mcV_A$ for an affine $\mcR$-homomorphism
$\alpha:\mcU\lra\mcV$. An example is given at the end.

Let $X$ be a locally compact Hausdorff space. By $X^{[d]}$ we denote
the configuration space of $d$ distinct unordered points in $X$,
i.e. the quotient $X^{(d)}/\Sigma_d$ where $X^{(d)}\subseteq X^d$ is
the subspace of injective maps $d\rightarrow X$. In this chapter we
will be interested in the \Cech cohomology with compact supports of
the space $(X\times Z)^{[d]}$ where $Z$ is another locally compact
Hausdorff space. Note that there is a canonical map
\[f:(X\times Z)^{[d]}\rightarrow\scS_d(X)=X^d/\Sigma_d\]
We want to apply Leray spectral sequence to this map to get some
cohomological information about the configuration space $(X\times
Z)^{[d]}$ and so we need to identify the fibres of $f$. Let
$Y=\{(x_1,k_1),\ldots,(x_n,k_n)\}\in\scS_d(X)$. The fibre
$f^{-1}(Y)$ is then obviously homeomorphic to
$Z^{[k_1]}\times\cdots\times Z^{[k_n]}$. Let us suppose now that
for all $k>0$ and all locally constant sheaves $\mcA$ on $X$
\[\check{H}^*_c(Z^{[k]};\mcA)=0\]
as is for example the case for $Z=\bbR_+$, the closed half-line.
Here $(\bbR_+)^{[k]}$ is topologically a $k$-simplex with all its
$(k-1)$-faces but one removed. Hence its one-point compactification
is homeomorphic to the $k$-dimensional disk $D^k$ and thus
contractible. The K\"unneth formula then says that for all
$k_1,\ldots,k_n>0$ and all locally constant sheaves $\mcA$ we have
\[\check{H}^*_c(Z^{[k_1]}\times\cdots\times Z^{[k_n]};\mcA)=0\]
Therefore in the Leary spectral sequence of $f$ with compact
supports we have $E_2^{pq}=0$ and therefore also
$\check{H}^*_c((X\times Z)^{[d]};\mcA)=0$ for any locally compact
Hausdorff space $X$ and any locally constant sheaf $\mcA$ on
$(X\times Z)^{[d]}$.

Now we apply these ideas to the homotopy fibre of $\alpha_A$ for
some affine $\mcR$-homomorphism $\alpha:\mcU\ra\mcV$, affine
$j^r$-representation $\varphi:\mcV\ra\Gamma E$ and a stratified
subset $A\subseteq E$ of codimension at least $\dim M+2$. Recall
that by Theorem \ref{mod_homotopy_fibre} this homotopy fibre can
be identified with $\hofib_{v}\alpha_A\simeq\mcW_{A\times I}$. See
this theorem for the explanation of $\mcW$ and the affine
$j^r$-representation $\psi$ via which $\mcW_{A\times I}$ is
defined.

Provided that $(\varphi\alpha)_x$ is surjective on $\hat{M}$,
Lemma \ref{mod_image_of_the_pullback_representation} guarantees
that $\psi_{(x,t)}$ is either surjective or disjoint with $A\times
I$ and that $\widehat{M\times I}=\hat{M}\times(0,1]$. As the
codimension of $A\times I$ in $E\times I$ is at least
$\dim(M\times I)+1$ we have a spectral sequence
(\ref{hpr_colimit_spectral_sequence}) for the affine
$j^r$-representation $\psi$ from the theorem
\[E^1_{pq}(\mcW)\cong\check{H}_c^{-ep-q}\left((A\times I)_{\ushortw{\widehat{M\times I}}}^{[-p]};w_1(h_{-p})+\sign\right)\Rightarrow\widetilde{H}_{p+q}(\mcW_{A\times I})\]
For the case of the configuration space of $\hat{M}\times(0,1]$
situation is the same as above: the fibres of $(A\times
I)_{\ushortw{\hat{M}\times(0,1]}}^{[-p]}\rightarrow
A_{\ushort{\hat{M}}}^{-p}/\Sigma_{-p}$ are products of
$(0,1]^{[k]}$ and we are getting $E^1_{pq}(\mcW)=0$. According to
Proposition \ref{hpr_zero_spectral_sequence_implies_convergence}
this spectral sequence converges and thus
$\widetilde{H}_*(\hofib_v\alpha_A)=0$. Therefore we have a theorem

\begin{theorem} \label{thm_main_theorem}
Let $M$ be a compact smooth manifold, let $\alpha:\mcU\ra\mcV$ be
an affine $\mcR$-homomorphism between two topologically finitely
generated affine $\mcR$-modules, let
$\varphi:\mcV\rightarrow\Gamma E$ be an affine
$j^r$-representation and $A\subseteq E$ a stratified subset of
codimension at least $\dim M+2$ such that outside the set
\[\hat{M}=\{x\in M\ |\ (\varphi\alpha)_x\textrm{ is surjective}\}\]
we have $\im\varphi_x\cap A=\emptyset$. Then each homotopy fibre
$\hofib_v\alpha_A$ of $\alpha_A$ is acyclic, i.e.
$\tilde{H}_*(\hofib_v\alpha_A)=0$. \qed
\end{theorem}

We finish with an example (and its refinement) of use of our main
Theorem \ref{thm_main_theorem}.

\begin{example}
Let us consider smooth functions $M\ra\bbR$ and let $A\subseteq
J^3(N,\bbR)$ be the complement of the set of $3$-jets which have at
most $A_2$-singularity (i.e. are either regular, nondegenerate
critical or have a singularity of type $A_2$). In other words the
complement of $A$ consists precisely of those jets which take in
some coordinate chart the form
 \[f(x_1,\ldots,x_m)=c\pm x_1^k\pm x_2^2\pm\cdots\pm x_m^2\]
with $1\leq k\leq 3$. Then according to Igusa's theorem (Theorem 9.1
of \cite{Igusa}) the canonical map
 \[(j^3)_A:C^\infty(M,\bbR)_A\ra\Gamma(J^3(M,\bbR))_A\]
is $(\dim M)$-connected and therefore its homotopy fibres are $(\dim
M-1)$-connected. In particular when $\dim M>1$ they are simply
connected and by Theorem \ref{thm_main_theorem} and Whitehead's
theorem they are weakly contractible making $(j^3)_A$ a weak
homotopy equivalence.
\end{example}

\begin{example}
To demonstrate an example in which the setting of topological
$\mcR$-modules can be easily applied we modify the previous example
slightly by considering a closed submanifold $N\subseteq M$ which is
either neat or $N=\partial M$ and fix a Morse function $g:M\ra\bbR$
(or just its germ at $N$). Next we consider the following
topologically finitely generated affine $\mcR$-modules
 \begin{eqnarray*}
  \mcU & = & \{f\in C^\infty(M,\bbR)\ |\ j^3_zf=j^3_zg\ \forall z\in N\} \\
  \mcV & = & \{s\in \Gamma(J^3(M,\bbR))\ |\ s(z)=j^3_zg\ \forall z\in N\}
 \end{eqnarray*}
According to Theorem \ref{thm_main_theorem} the homotopy fibres of
the restricted jet prolongation $(j^3)_A:\mcU_A\ra\mcV_A$ (with $A$
from the previous example) are again acyclic. Here $\mcU_A$ is the
space of functions $M\ra\bbR$ with at most $A_2$-singularities which
agree with $g$ up to order $3$ at $N$ and similarly for $\mcV_A$ and
we obtain a relative version of the previous example.
\end{example}




\newpage
\appendix


\chapter{Locally convex topological vector spaces} \label{LCTVS}

This appendix serves as a source of results about locally convex
topological vector spaces for the purpose of the thesis. The main
source used was \cite{Treves} where one can find everything (in a
great more detail) with the exception of Theorem
\ref{LCTVS_monoidal_structure}.

We work here with vector spaces over the field $\bbR$ of real
numbers. A topological vector space $\mcV$ is called
\textit{locally convex} if every neighbourhood of $0$ contains a
convex neighbourhood of $0$.

A subset $V\subseteq\mcV$ is called \textit{balanced} if
$rV\subseteq V$ for all $|r|\leq 1$ and \textit{absorbing} if for
every $v\in\mcV$ there is $r>0$, such that the line segment from
$-v$ to $v$ is contained in $rV$. It is called a \textit{barrel}
if it closed, convex, balanced and absorbing. Equivalently, $\mcV$
is locally convex if every neighbourhood $V$ of $0$ contains a
barrel neighbourhood of $0$.

A seminorm on a vector space $\mcV$ is a function
$\rho:\mcV\ra\bbR$ such that
\begin{itemize}
\item[(a)]{$\rho(v)\geq 0$}
\item[(b)]{$\rho(v+v')\leq\rho(v)+\rho(v')$}
\item[(c)]{$\rho(\lambda\cdot v)=|\lambda|\cdot\rho(v)$}
\end{itemize}
Every continuous seminorm on a topological vector space $\mcV$
gives a closed ball $B_\rho=\rho^{-1}[0,1]$. They are barrel
neighbourhoods of $0$. If on the other hand $V$ is a barrel
neighbourhood of $0$ in $\mcV$ then there exists a continuous
seminorm $\rho:\mcV\ra\bbR$, the so called Minkowski functional,
such that $V=\rho^{-1}[0,1]$. It is defined by the formula
\[\rho(v)=\inf\{\lambda\in\bbR_+\ |\ 1/\lambda\cdot v\in V\}\]
Therefore, if $\mcV$ is locally convex its topology can be
completely described by continuous seminorms.

A family $\mcP$ of continuous seminorms on a LCTVS $\mcV$ is
called a \emph{basis of continuous seminorms} if for any
continuous seminorm $\sigma$ on $\mcV$ there is $r>0$ and
$\rho\in\mcP$ such that $\sigma\leq r\rho$, or equivalently if the
closed $\varepsilon$-balls
$B_\rho(\varepsilon)=\rho^{-1}[0,\varepsilon]$ form a
neighbourhood basis of $0$. In such a case we say that the family
$\mcP$ \emph{define} the topology of $\mcV$.

For example, if $\mcV$ is a (not necessarily locally convex) TVS,
the family of all continuous seminorms define some locally convex
topology on $\mcV$. It is always weaker than the one that we
started with and this procedure provides a left adjoint of the
forgetful functor $\LCTVS\hookrightarrow\TVS$ (hence $\LCTVS$ is a
reflective subcategory of $\TVS$).

If $\mcU$, $\mcV$ are two vector spaces,
$\alpha:\mcU\times\mcV\rightarrow\mcU\otimes\mcV$ the canonical
map and $U\subseteq\mcU$, $V\subseteq\mcV$ any convex subsets we
define their convex tensor product $U\otimes
V\subseteq\mcU\otimes\mcV$ to be the convex hull of the set
$\alpha(U\times V)=\{u\otimes v\ |\ u\in U, v\in V\}$. If $\rho$,
$\sigma$ are seminorms on $\mcU$, $\mcV$ we define a seminorm
$\rho\otimes\sigma:\mcU\otimes\mcV\ra\bbR$ by the requirement
$B_{\rho\otimes\sigma}=B_\rho\otimes B_\sigma$. Therefore if
$\mcU$ and $\mcV$ are both LCTVS we can define a locally convex
topology on their tensor product $\mcU\otimes\mcV$ via the
seminorms $\rho\otimes\sigma$ where $\rho$, $\sigma$ range over
all continuous seminorms on $\mcU$, $\mcV$ (or for that purpose
ranging over any bases of seminorms on the respective spaces).
Together with this topology we call $\mcU\otimes\mcV$ the
projective tensor product of $\mcU$ and $\mcV$.

\begin{proposition}
If $\mcU$ and $\mcV$ are both LCTVS then the projective topology
on $\mcU\otimes\mcV$ is the strongest LCTVS topology for which
$\alpha:\mcU\times\mcV\rightarrow\mcU\otimes\mcV$ is continuous.
It enjoys the following universal property: every
\emph{continuous} bilinear map
$\beta:\mcU\times\mcV\rightarrow\mcW$ to a LCTVS $\mcW$ factorizes
uniquely through $\alpha$ as
\[\beta:\mcU\times\mcV\xlra{\alpha}\mcU\otimes\mcV\xlra{\gamma}\mcW\]
with $\gamma$ a \emph{continuous} linear map. In other words there
is a natural isomorphism
\[\LCTVS(\mcU\otimes\mcV,\mcW)\cong\Bilin(\mcU,\mcV;\mcW)\]
\end{proposition}

\begin{proof}
The first part follows straight from the definition. According to
the algebraic properties of a tensor product we only need to show
that the continuity of $\beta$ implies the continuity of $\gamma$.
Therefore let $W\subseteq\mcW$ be any convex neighbourhood of $0$.
By continuity of $\beta$ there are $\rho$, $\sigma$ such that
$\beta(B_\rho\times B_\sigma)\subseteq W$. As $W$ is convex it
also contains
\[\ch(\beta(B_\rho\times B_\sigma))=
\ch(\gamma\alpha(B_\rho\times B_\sigma))=
\gamma(\ch\alpha(B_\rho\times B_\sigma))=
\gamma(B_{\rho\otimes\sigma})\] where $\ch$ denotes the convex
hull.
\end{proof}

\begin{theorem} \label{LCTVS_monoidal_structure}
The category $\LCTVS$ of locally convex topological vector spaces
has the following structure
\begin{itemize}
\item[(i)]{$\LCTVS$ together with $\otimes$ is a symmetric monoidal category.}
\item[(ii)]{The tensor product functor $\mcU\otimes-$ commutes with
finite 
colimits.
}
\item[(iii)]{If $\mcU,\mcV,\mcW\in\LCTVS$, $P$ is a locally compact
space and $P\times(\mcU\times\mcV)\rightarrow\mcW$ is a continuous
pointwise $\bbR$-bilinear map then the induced map
$P\times(\mcU\otimes\mcV)\rightarrow\mcW$ is also continuous. In
particular the tensor product preserves homotopies.}
\end{itemize}
\end{theorem}

\begin{proof}
The associativity and commutativity of $\otimes$ follow from the
description of the topology on the tensor product, e.g. the
topology on
$(\mcU\otimes\mcV)\otimes\mcW\cong\mcU\otimes(\mcV\otimes\mcW)$ is
generated by neighbourhoods of $0$ of the form $(B_\rho\otimes
B_\sigma)\otimes B_\tau=B_ \rho\otimes(B_ \sigma\otimes B_\tau)$.

Tensoring a biproduct diagram with $\mcU$ clearly produces again a
biproduct diagram (mere additivity of the tensor product suffices)
so $\mcU\otimes-$ preserves finite coproducts. The fact that
$\mcU\otimes-$ preserves cokernels follows from the corresponding
property of the product $\mcU\times-$
.

To prove (iii) we use adjunction to reduce everything to showing
that $\map(P,\mcW)$ is locally convex. But a subbasis for the
topology of $\map(P,\mcU)$ is given by convex sets
$\{f:P\rightarrow\mcU\ |\ f(K)\subseteq U\}$ with $K$ ranging over
all compact and $U$ over all convex
open subsets.
\end{proof}

\begin{note} 

Any locally convex topological $\bbR$-algebra $\mcA$ is a monoid in
the monoidal structure on $\LCTVS$ and topological $\mcA$-modules
are precisely
$\mcA$-modules 
in the monoidal category sense. Consequently we can use monoidal
category techniques to deal with $\mcA$-modules. For simplicity we
assume that $\mcA$ is commutative. For $\mcA$-modules $\mcU$ and
$\mcV$ we get a tensor product $\mcU\tensor{\mcA}\mcV$ as a
coequalizer
\[\xymatrix{{\mcU\otimes\mcA\otimes\mcV} \ar[r]<3pt> \ar[r]<-3pt> & {\mcU\otimes\mcV}}\]
the two maps being the two structure maps. In this way
$\mcA$-modules become a symmetric monoidal category (here the
right exactness of the tensor product is necessary in order to
define the action of $\mcA$ on $\mcU\tensor{\mcA}\mcV$).

Suppose that $\varphi:\mcA\rightarrow\mcB$ is a continuous algebra
homomorphism between commutative topological $\bbR$-algebras
$\mcA$ and $\mcB$ (i.e. a homomorphism of commutative monoids in
the monoidal category $\LCTVS$). Then the forgetful functor
$\varphi^*:\Bmod\rightarrow\Amod$ has a left adjoint
$\mcB\tensor{\mcA}-$.
\end{note}



\chapter{Differential topology, function spaces} \label{function_spaces}

The purpose of this appendix is to give a brief overview of the
differential topology used and to prove some auxiliary results
that would, if included in the main text, disturb the flow of
exposition. The two main sources used were
\cite{GolubitskyGuillemin} and \cite{Hirsch}.

Let $E\ra M$ be a smooth fibre bundle. There are two kinds of
topologies on the space $\Gamma(E)$ of smooth sections, the weak
(or sometimes called compact-open) topology and strong (or
Whitney) topology defined for each degree $0\leq r\leq\infty$ of
differentiability. We start with the definition for $r=0$. The
weak topology on the space $\Gamma^0(E)$ of continuous sections of
$E$ is just the usual compact-open topology. We denote the
resulting space by $\Gamma^0_W(E)$. The basis for the strong
topology on $\Gamma^0(E)$ is indexed by open subsets $U$ of $E$
for which we have a generating open subset
\[\{s\in\Gamma^0(E)\ |\ \im(s)\in U\}\subseteq\Gamma^0(E)\]
The resulting topological space is denoted by $\Gamma^0_S(E)$. For
finite $r$ and $*=W,S$ we obtain the topology on $\Gamma^r_*(E)$
via the jet prolongation map
\begin{equation} \label{fun_jet_prolongation_map}
j^r:\Gamma^r(E)\lra\Gamma^0_*(J^r E)
\end{equation}
and the above definition for $r=0$. More precisely this map is
injective and we give $\Gamma^r_*(E)$ the subspace topology. In
the diagram
\[\xymatrix{
\Gamma^{r+1}_*(E) \ar@{c->}[r] \ar@{-->}[d] & \Gamma^0_*(J^{r+1}E)
\ar[d]
\\
\Gamma^r_*(E) \ar@{c->}[r] & \Gamma^0_*(J^rE)}\] the dashed arrow is
continuous by the properties of subspaces (or one can say that this
arrow exists in $\Top$) and is clearly the canonical inclusion. We
define the topology on $\Gamma^\infty_*(E)$ as the limit topology
\[\Gamma^\infty_*(E)\cong\lim\limits_{r\ra\infty}\Gamma^r_*(E)\]
In other words one just takes all the open subsets from all
$\Gamma^r_*(E)$, $r=0,1,\ldots$, together. There is a description of
this topology similar to (\ref{fun_jet_prolongation_map}) using the
infinite jet bundle $J^\infty E$. As a space over $M$ it is a limit
of
\[\xymatrix{
{}\cdots \ar[r] & J^rE \ar[r] \ar[d] & {}\cdots \ar[r] & J^1E \ar[r]
\ar[d] & E \ar[d]
\\
{}\cdots \ar@{=}[r] & M \ar@{=}[r] & {}\cdots \ar@{=}[r] & M
\ar@{=}[r] & M}\] We have a similar jet prolongation map
\[j^\infty:\Gamma^\infty_*(E)\lra\Gamma^0_*(J^\infty E)\]
and it is again an inclusion of a subspace. For a \emph{compact}
manifold $M$ the two topologies agree
$\Gamma^r_W(E)=\Gamma^r_S(E)$, as can be easily deduced from the
case $r=0$.

In the case when $E=N\times M$ is the trivial bundle (so that
$\Gamma^r(E)$ is the space of all $C^r$-maps $M\rightarrow N$) we
write $C^r_*(M,N)$ instead of $\Gamma^r_*(E)$. For the
differential topology one of the most important features of these
topologies is that they are all Baire spaces, i.e. countable
intersections of open dense subsets are dense.

As we are only interested in the case $r=\infty$ this is what we
are going to assume from now on. The following results explain how
could smooth maps be defined by specifying then on subsets such
that they agree on the intersections.

\begin{lemma} \label{fun_sheaf_property}
Let $E\ra M$ be a smooth bundle and let $\scU$ be an open covering
of $M$. Then the restriction maps $\Gamma^\infty_W(E)\ra
\Gamma^\infty_W(U,E)$ for $U\in\scU$ induce a homeomorphism
 \[\Gamma^\infty_W(E)\xra{\cong}\lim_{U\in\scU}\Gamma^\infty_W(U,E)\]
\end{lemma}

\begin{proof}
Firstly the above map is a homeomorphism when $\infty$ is replaced
by $0$, i.e. on the level of spaces of continuous sections with
the compact-open topology. As for each $r$ we have
 \[\Gamma^0_W(J^rE)\xra{\cong}\lim_{U\in\scU}\Gamma^0_W(U,J^rE)\]
and subspaces commute with limits we get a homeomorphism
 \[\Gamma^r_W(E)\xra{\cong}\lim_{U\in\scU}\Gamma^r_W(U,E)\]
by restriction. Finally, taking a limit for $r\ra\infty$ gives the
result (or we could take $r=\infty$ straight away).
\end{proof}

\begin{lemma} \label{fun_pullback}
If $M=M_1\cup M_2$ is a union of two submanifolds of codimension 0
meeting at their common boundary $M_0=M_1\cap M_2=\partial
M_1\cap\partial M_2$ then the following diagram is a pullback
square
\[\xymatrix{
\Gamma^\infty_*(E) \ar[r] \ar[d] & \Gamma^\infty_*(M_1,E) \ar[d]
\\
\Gamma^\infty_*(M_2,E) \ar[r] & \Gamma^0_*(M_0,J^\infty E)}\] Here
$\Gamma^\infty_*(M_i,E)$ denotes the space of smooth sections of $E$
over $M_i$.\qed
\end{lemma}

Before we give some applications of this lemma we switch our
attention to the case of sections of a \emph{vector} bundle $E$.
In this case it can be shown that the addition or subtraction of
sections is continuous so that $\Gamma_*(E)$ is actually a
topological abelian group. More care has to be taken with the
multiplication of sections. The map
\[C^\infty_*(M,\bbR)\times\Gamma_*(E)\lra\Gamma_*(E)\]
is continuous so that $\Gamma_*(E)$ is a topological module over
the topological ring $C^\infty_*(M,\bbR)$. The bad news is that
the inclusion $\bbR\ra C^\infty_*(M,\bbR)$ is continuous only for
$*=W$ or when $M$ is compact (in which case the two topologies
agree anyway).

From now on we assume that $M$ is compact and write $\Gamma(E)$ for
$\Gamma^\infty_*(E)$ (we do not have to distinguish between the weak
and strong topology now). As $\Gamma(E)$ is a topological vector
space we only need to understand the neighbourhood basis at $0$. If
one chooses a metric on each vector bundle $J^rE$ we can construct
one by taking
\[\{s\in\Gamma(E)\ |\ \forall x\in M:||j^r_xs||<\varepsilon\}\]
Therefore one easily sees that $\Gamma(E)$ is a LCTVS which is
defined by seminorms
\[s\mapsto\sup_{x\in M}||j^r_xs||\]
Using Lemma \ref{fun_pullback} we prove continuity of various maps
that are used in the main chapters.

\begin{corollary}
Let us denote $\mcR=C^\infty(M,\bbR)$ and $\mcR^I=C^\infty(M\times
I,\bbR)$, let $\mcS$ be the subspace of $\mcR^I\times\mcR^I$
consisting of those $(f,g)$ for which there exists $\varepsilon>0$
such that
\[f(x,1-t)=f(x,1)=g(x,0)=g(x,t)\]
for all $0\leq t\leq\varepsilon$. The ``concatenation'' map
\[\mu:(0,1)\times\mcS\lra\mcR^I\]
sending $(s,f,g)$ to the function\footnote{Geometrically just $f$
reparametrized to $M\times[0,s]$ patched with $g$ reparametrized to
$M\times[s,1]$.}
\[h(x,t)=\left\{
\begin{array}{ll}
f(x,t/s) & \textrm{when }0\leq t\leq s
\\
g(x,(t-s)/(1-s)) & \textrm{when }s\leq t\leq 1
\end{array}\right.\]
is continuous\footnote{This is true even on a bigger subspace of
those $(f,g)$ whose infinite jets are compatible.}.
\end{corollary}

\begin{proof}
We write $I\times(0,1)$ as a union of two submanifolds
\[P_1=\{(t,s)\ |\ s\leq t\} \qquad
P_2=\{(t,s)\ |\ s\geq t\}\] of codimension $0$ and define maps
$\alpha_i:P_i\rightarrow I$
\[\alpha_1(t,s)=t/s \qquad
\alpha_2(t,s)=(t-s)/(1-s)\] and therefore maps
\[\mcR^I=C^\infty(M\times I,\bbR)\xlra{(\id\times\alpha_i)^*}
C^\infty_W(M\times P_i,\bbR)\] (for the continuity of this map see
the first paragraph after the proof). Taking a product of these
maps we obtain the bottom map in the diagram
\[\xymatrix{
{}\mcS \ar@{-->}[r]^-{\beta} \ar[d] & C^\infty_W(M\times
I\times(0,1),\bbR) \ar[d]
\\
{}\mcR^I\times\mcR^I \ar[r] & C^\infty_W(M\times P_1,\bbR)\times
C^\infty_W(M\times P_2,\bbR)}\] According to Lemma
\ref{fun_pullback} the dashed arrow is continuous. We can write
the concatenation map $\mu$ as a composition
\[(0,1)\times\mcS\xlra{\id\times\beta}(0,1)\times C^\infty_W(M\times
I\times(0,1),\bbR)\xlra{pev}C^\infty(M\times I,\bbR)=\mcR^I\] with
$pev$ being the partial evaluation $(s,f)\mapsto f(-,s)$.
\end{proof}

The composition $C^\infty(N,P)\times C^\infty(M,N)\lra
C^\infty(M,P)$ is continuous in the weak topologies and for strong
topologies it is continuous on the subspace formed by $(f,g)$ with
$g$ proper. This proves most of

\begin{lemma}
Let $\lambda:I\ra I$ be a smooth function. The reparametrization
map
\[\lambda^*:\mcR^I\ra\mcR^I\]
defined by $f(x,t)\mapsto f(x,\lambda(t))$ is continuous. Moreover
for a smooth family $\lambda_\tau:I\ra I$ the induced homotopy
\[\lambda_\tau^*:\mcR^I\ra\mcR^I\]
is also continuous.
\end{lemma}

\begin{proof}
Again the homotopy can be written as
\[I\times C^\infty(M\times I,\bbR)\xlra{(\id\times\lambda)^*}
I\times C^\infty(M\times I\times I,\bbR)\xlra{pev}C^\infty(M\times
I,\bbR)\] with $\lambda:I\times I\ra I$ sending $(t,\tau)$ to
$\lambda_\tau(t)$.
\end{proof}

\begin{corollary}
Let $0<s<1$ and $\lambda:I\ra I$, $\rho:I\ra I$ two functions such
that $\lambda(0)=0$, $\lambda=1$ near $1$ and similarly
$\rho(1)=1$, $\rho=0$ near $0$. The concatenation map
\[\mu(s,\lambda,\rho):\mcR^I\times_\mcR\mcR^I\lra\mcR^I\]
given by $(f,g)\mapsto\mu(s,\lambda^*f,\rho^*g)$ is continuous.
Moreover for any other choice of $s$, $\lambda$ and $\rho$ the
resulting map is homotopic via maps of the same form.
\end{corollary}

\begin{proof}
Any two values of $s$ can be joined by a linear homotopy as can be
any two values of $\lambda$ and $\rho$.
\end{proof}


\chapter{Sheaves}

The main references for this appendix are the books \cite{Bredon}
and \cite{GelfandManin}. The idea of this appendix is to gather
facts about sheaves and in particular about sheaf cohomology that
are used in the thesis (mainly in Appendix \ref{dual_map} and in
Chapter~\ref{h_principle}).

The basic intuition behind sheaves is that they are systems of
abelian groups para\-met\-riz\-ed by the points of a topological
space $X$. In one approach they indeed are certain spaces over
$X$. In this sense they generalize the so-called systems of local
coefficients (also known as bundles of abelian groups). Moreover
they constitute an abelian category on which a variety of functors
is defined. For example, we have a functor of global sections
which associates to a sheaf the abelian group of its sections. Its
derived functors are the sheaf cohomology groups with coefficients
in the sheaf under consideration. This is roughly the content of
the first section.

There are other cohomology theories with coefficients in sheaves
(or in a special class of sheaves) and among them the sheaf
cohomology is ``universal''. For paracompact spaces most of these
theories agree. One exception is the singular cohomology which
does not behave well for spaces with local complexities. So far
for the section 2.

In the last section we pursue more advanced properties of sheaves
and in particular derive the Leray spectral sequence. Given a
continuous map $f:X\rightarrow Y$ it relates the sheaf cohomology
of $X$ to the sheaf cohomology of $Y$ and of the fibres
$f^{-1}(y)$ over the points $y\in Y$.

\section{Sheaves, sheaf cohomology}

Let $X$ be a topological space\footnote{we do not assume here that
it is compactly generated Hausdorff}. We denote by $\Op(X)$ the
poset of all open subsets of $X$ ordered by inclusion. A
\emph{presheaf} on $X$ is just a functor
$\Op(X)^{op}\rightarrow\Set$ to the category of sets. Therefore we
have a category of presheaves on $X$, namely the functor category
$\Pre_X\Set=\PreSheaf{X}{\Set}$. In a similar way we define a
presheaf in any category $\mcC$ as an object of
$\Pre_X\mcC=\PreSheaf{X}{\mcC}$. A presheaf $\mcA\in\Pre_X\mcC$ is
said to be a \emph{sheaf} if for every subset $\scU\subseteq\Op(X)$
which is closed under subsets the map
\[\mcA(\cup\scU)\lra\lim\limits_{U\in\scU}\mcA(U)\]
(with $\cup\scU=\bigcup\limits_{U\in\scU}U$) is an isomorphism,
i.e. if $\mcA$ preserves limits over all hereditary subsets. In
this way we get a full subcategory $\Sh_X\mcC$ of sheafs of
objects in $\mcC$ on $X$.

Suppose that $\mcC$ is
an ``algebraic category'' - a category
$\langle\Omega,E\rangle\textrm{-}\Alg$ of algebras\footnote{Here
an $\langle\Omega,E\rangle$-algebra is a set $S$ together with a
realization of $\Omega$ by operations $S^n\rightarrow S$ that
satisfy the relations $E$. The category
$\Sh_X\langle\Omega,E\rangle\textrm{-}\Alg$ can be alternatively
described as the category of $\langle\Omega,E\rangle$-algebras in
the category $\Sh_X\Set$ of sheaves of sets.} defined by
operations $\Omega$ and relations $E$.
Then the inclusion $\iota:\Sh_X\mcC\hookrightarrow\Pre_X\mcC$ has
a left adjoint $s:\Pre_X\mcC\rightarrow\Sh_X\mcC$ sometimes called
sheafification. It is given by the formula
\[s\mcA(U)=\colim\limits_{\scV\in\Cov_h(U)}\lim\limits_{V\in\scV}\mcA(V)\]
where the colimit is taken over all hereditary open coverings of
$U$.

A more illuminating description of the sheafification functor can
be obtained by introducing total spaces of presheaves. Let
$\mcA\in\Pre_X\Set$ be a presheaf of sets on $X$. We construct a
space $\pi_\mcA:\Tot(\mcA)\rightarrow X$ over $X$ called the
\emph{total space} of $\mcA$
\[\pi_\mcA:\Tot(\mcA)=\int\limits^{U\in\Op(X)}U\times\mcA(U)\lra\int\limits^{U\in\Op(X)}U=X\]
The projection $\pi_\mcA$ can be shown to be a local homeomorphism
and there is a bijection between $(s\mcA)(U)$ and the set
$\Gamma(U,\Tot(\mcA))$ of sections of $\Tot(\mcA)$ over $U$. In
this way the category $\Sh_X\Set$ of sheaves of sets on $X$
becomes equivalent to the category of local homeomorphisms $A\ra
X$. Moreover $\Tot$ is identified with the sheafification functor
$s$ and the functor $\Gamma(-,A)$ of local sections with the
inclusion functor $\iota$.

For the category $\mcC=\langle\Omega,E\rangle\textrm{-}\Alg$ we
get an equivalence of $\Sh_X\mcC$ with the category of
$\langle\Omega,E\rangle$-algebras in $\Sh_X\Set$, i.e. the
category of local homeomorphisms $A\ra X$ with a fibrewise
structure of an $\langle\Omega,E\rangle$-algebra. Of course
covering maps are local homeomorphisms and therefore they give
examples of sheaves, the so-called \emph{locally constant}
sheaves. Also systems of local coefficients (or bundles of groups)
on $X$ just correspond to locally constant sheaves of abelian
groups. The fibres $\mcA_x:=\pi_\mcA^{-1}(x)$ of the projection
$\pi_\mcA$ are called the \emph{stalks} of the sheaf $\mcA$.
Alternatively one can define them as the directed colimits
\[\mcA_x=\colim\limits_{\substack{U\in\Op(X)\\x\in U}}\mcA(U)\]

Now assume that in addition
$\mcC=\langle\Omega,E\rangle\textrm{-}\Alg$ is an abelian
category, main example being the category of abelian groups. Then
both $\Pre_X\mcC$ and $\Sh_X\mcC$ become also abelian with $\iota$
left exact and $s$ exact. Moreover if $\mcC$ has enough injectives
then so do both $\Pre_X\mcC$ and $\Sh_X\mcC$ and if $\mcC$ has
enough projectives then the same is true for $\Pre_X\mcC$
\emph{but not} for $\Sh_X\mcC$ in general. From now on we assume
that $\mcC$ and hence also $\Sh_X\mcC$ have enough injectives.

The abelian structure of the category $\Sh_X\mcC$ is reflected in
the stalks, namely a sequence of sheaves is exact if and only if
at each point the corresponding sequence of stalks is exact. This
is where the geometric description of a sheaf (in terms of the
total space) is very useful.

It is also this description that gives a name to the elements of
$\mcA(U)$ - they are called sections of $\mcA$ over $U$. Hence we
often write $\Gamma(U,\mcA)$ instead of $\mcA(U)$ and for global
sections simply $\Gamma(\mcA)$. In this way we get a ``global
sections'' functor $\Gamma:\Sh_X\mcC\rightarrow\mcC$. The essence
of the sheaf theory is that this functor is not exact but merely
left exact. As the category $\Sh_X\mcC$ has enough injectives one
can define the right derived functors of $\Gamma$. In this
appendix we denote $R^n\Gamma(\mcA)\in\mcC$ by $H^n(X;\mcA)$ and
call it the $n$-th \emph{sheaf cohomology} of $X$ with
coefficients in the sheaf $\mcA$. If $\varphi$ is a family of
supports on $X$ we can define $\Gamma_\varphi(\mcA)$, the sections
of $\mcA$ with supports in $\varphi$, to be the kernel of
\[\Gamma(\mcA)\lra\colim\limits_{F\in\varphi}\Gamma(X-F;\mcA)\]
Again $\Gamma_\varphi$ is left exact and there are right derived
functors $R^n\Gamma_\varphi(\mcA)=H^n_\varphi(X;\mcA)$, the $n$-th
sheaf cohomology of $X$ with supports $\varphi$ and with
coefficients in $\mcA$.

\section{The relation to other cohomology theories}
\label{sheaves_relation_to_other_cohomology_theories}

In this section we relate the sheaf cohomology to the singular
cohomology $H_\Delta^*$ and to the \Cech cohomology $\check{H}^*$.

First we treat the singular cohomology with constant coefficients
and any supports. For an open subset $U\in\Op(X)$ let us condiser
the singular cochain complex $C^*(U;\bbZ)$ on $U$. It provides a
cochain complex of  presheaves $C^*(-;\bbZ)$ of abelian groups on
$X$. By definition
\[H_\Delta^*(X)=H^*(C^*(X;\bbZ))\]
We denote by $\mcS^*\in\Sh_X\Ab$ the cochain complex of generated
sheaves, i.e.
\[\mcS^*=sC^*(-;\bbZ)\]
It turns out (see Section I.7 of \cite{Bredon}) that if $X$ is
paracompact then
\[H_\Delta^*(X)\cong H^*(\mcS^*(X))=H^*(\Gamma(\mcS^*))\]
We can define more generally the singular cohomology groups of $X$
with supports in $\varphi$ and coefficients in a locally constant
sheaf $\mcA$ as the cohomology of the cochain complex
\[C^p_\varphi(X;\mcA)=\{\alpha:\map(\Delta^p,X)\rightarrow\map(\Delta^p,\Tot(\mcA))\ |\
\pi_*\alpha=\id,\ \supp\alpha\in\varphi\}\] of presheaves of
abelian groups. If the family $\varphi$ is paracompactifying we
get again an isomorphism
\[H_{\Delta,\varphi}^*(X;\mcA)\cong H^*(\Gamma_\varphi(\mcS^*_\mcA))\]
where $\mcS^*_\mcA$ is the sheaf generated by the presheaf
$C^*(-;\mcA)$. Under the same assumptions the sheaves
$\mcS^*_\mcA$ (assuming that $\mcA$ is finitely generated, see
Section III.1 of \cite{Bredon}) turn out to be
$\Gamma_\varphi$-acyclic so that $H^*_{\Delta,\varphi}$ computes
the sheaf cohomology of $X$ with supports in $\varphi$ and with
coefficients in $\mcA$ provided that the cochain complex
\[\cdots\ra0\ra\mcA\ra\mcS^0_\mcA\ra\mcS^1_\mcA\ra\cdots\]
is acyclic. The stalk at $x$ of its $n$-th cohomology is clearly
\[\colim\limits_{\substack{U\in\Op(X)\\x\in U}}\tilde{H}_\Delta^*(U;\mcA|U)\]
Therefore for a paracompactifying family $\varphi$ the natural map
$H_{\Delta,\varphi}^*(X;\mcA)\ra H_\varphi^*(X;\mcA)$ is an
isomorphism if this colimit vanishes or in other words if all the
points $x\in X$ are taut with respect to the singular cohomology.
This is for example the case for all locally contractible spaces
$X$. Also it is implied by the condition HLC - homological local
connectedness
.

In general there is a convergent spectral sequence
\[E_2^{pq}=H^p_\varphi(X;sH^q_\Delta(-;\mcA))\Rightarrow H^{p+q}_{\Delta,\varphi}(X;\mcA)\]

Next comes the Alexander-Spanier cohomology.
There is a sequence of abelian groups
\[AS^p(U;\bbZ)=\map(U^{p+1},\bbZ)\]
It is easy to define the coface and codegeneracy maps making it
into a cosimplicial presheaf of abelian groups. The associated
cochain complex will be denoted again by $AS^*(-;\bbZ)$. Let
$\mcA\mcS^*$ be the sheaf generated by this presheaf. Then by
definition
\[H^*_{AS,\varphi}=H^*(\Gamma_\varphi(\mcA\mcS^*))\]
Again if $\varphi$ is paracompactifying then all the sheaves
$\mcA\mcS^*$ are $\Gamma_\varphi$-acyclic and
$\check{H}^*_\varphi$ is isomorphic to the corresponding sheaf
cohomology provided that $\mcA\mcS^*$ constitutes a resolution of
$\bbZ$. Unlike in the case of singular cohomology this is
\emph{always} the case.

Let $\mcA\in\Pre_X\Ab$ be a presheaf and $\scU\subseteq\Op(X)$ any
open covering of $X$. We define
\[\check{C}^p(\scU;\mcA)=\prod_{(U_0,\ldots,U_p)\in\scU^{p+1}}\mcA(U_0\cap\cdots\cap U_p)\]
Again one can make $\check{C}^*(\scU;\mcA)$ into a cosimplicial
abelian group and further into a cochain complex. We define the
\Cech cohomology of the covering $\scU$ with coefficients in
$\mcA$ as
\[\check{H}^*(\scU;\mcA)=H^*(\check{C}^*(\scU;\mcA))\]
If $\scV$ is a refinement of $\scU$ and
$\lambda:\scV\rightarrow\scU$ is a function such that
$V\subseteq\lambda V$ then we define a map
\[\lambda^*:\check{C}^p(\scU;\mcA)\lra\check{C}^p(\scV;\mcA)\]
using the restriction maps as in the diagram
\[\xymatrix{
{}\check{C}^p(\scU;\mcA) \ar@{-->}[r] \ar[d] &
{}\check{C}^p(\scV;\mcA) \ar[d]
\\
{}\mcA(\lambda V_0\cap\cdots\cap\lambda V_p) \ar[r]^{res} &
{}\mcA(V_0\cap\cdots\cap V_p)}\] It is a map of cosimplicial
abelian groups and therefore induces a map
\[\lambda^*:\check{H}^*(\scU;\mcA)\lra\check{H}^*(\scV;\mcA)\]
on cohomology. For a different choice of $\lambda$ one can show
that the resulting map of cosimplicial abelian groups is homotopic
to $\lambda^*$ and so on the level of cohomology one gets an
equality. Hence we can define
\[\check{H}^*(X;\mcA)=\colim\limits_{\scU\in\Cov(X)}\check{H}^*(\scU;\mcA)\]
where the colimit is taken over all open coverings of $X$ ordered
by refinement. This group is called the \Cech cohomology of $X$
with coefficients in the presheaf $\mcA$.

Consider the following composition
\[\Gamma:\Sh_X\Ab\xlra{\check{C}^*(\scU;-)}\Ch^{\geq 0}\Ab\xlra{H^0}\Ab\]
As the derived functors are easily seen to be $R^pH^0\cong H^p$
and
$(R^q\check{C}^*(\scU;-))(\mcA)\cong\check{C}^*(\scU;H^q(-;\mcA))$
the Grothendieck spectral sequence takes the form
\[E_2^{pq}(\scU)\cong\check{H}^p(\scU;H^q(-;\mcA))\Rightarrow H^{p+q}(X;\mcA)\]
Taking a colimit over the open coverings $\scU$ of $X$ one gets
again a convergent spectral sequence
\[E_2^{pq}\cong\check{H}^p(X;H^q(-;\mcA))\Rightarrow H^{p+q}(X;\mcA)\]
According to Spanier (with an additional assumption of $X$ being
paracompact) for any presheaf $\mcB\in\Pre_X\Ab$ we have an
isomorpism
\[\check{H}^*(X;\mcB)\xlra{\cong}\check{H}^*(X;s\mcB)\]
As clearly $sH^q(-;\mcA)=0$ for $q>0$ and $sH^0(-;\mcA)=\mcA$ the
colimit spectral sequence collapses to give an
isomorphism\footnote{even when $X$ is not paracompact this map is
an isomorphism for $*=0,1$ and a monomorphism for $*=2$.}
\[\check{H}^*(X;\mcA)\xlra{\cong}H^*(X;\mcA)\]

\section{The sheaf cohomology and continuous maps}

Now we describe how sheaves behave with respect to continuous maps.
If $f:X\rightarrow Y$ is a continuous map then we get an induced
functor $f^*:\Op(Y)\rightarrow\Op(X)$, given simply by
$f^*U=f^{-1}(U)$, and therefore a functor between the presheaf
categories called the \emph{direct image}
\[f_\bullet:\Pre_X\mcC\rightarrow\Pre_Y\mcC\]
It is easy to check that this functor preserves the subcategory of
sheaves and we have
\[f_\bullet:\Sh_X\mcC\rightarrow\Sh_Y\mcC\]
Explicitly $f_\bullet\mcA(U)=\mcA(f^{-1}(U))$. The functor
$f_\bullet:\Pre_X\mcC\rightarrow\Pre_Y\mcC$ has a left adjoint
$\Lan_{f^*}$, the left Kan extension along $f^*$. Therefore we
have a composition
\[f^\bullet:\Sh_Y\mcC\xlra{\iota}\Pre_Y\mcC\xlra{\Lan_{f^*}}\Pre_X\mcC\xlra{s}\Sh_X\mcC\]
It is called the \emph{inverse image} and we claim that it is a
left adjoint of $f_\bullet$. This is clear from
\[\Sh_X\mcC(f^\bullet\mcA,\mcB)\cong\Pre_X\mcC(\Lan_{f^*}\iota\mcA,\iota\mcB)
\cong\Pre_Y\mcC(\iota\mcA,f_\bullet\iota\mcB)\cong\Sh_Y\mcC(\mcA,f_\bullet\mcB)\]
A very nice description of the functor $f^\bullet$ is in terms of
the total spaces. It is plainly a pullback functor along the map
$f$, sending $A\rightarrow Y$ to $A\times_YX\rightarrow X$.
\[\xymatrixnocompile{
\smash[b]{A\times_YX} \ar[r] \pullback \ar[d] & A \ar[d] \\
X \ar[r] & Y}\] In particular $(f^\bullet\mcA)_x\cong\mcA_{f(x)}$
and therefore the functor $f^\bullet$ is exact. On the other hand
$f_\bullet$ is merely left exact and thus one can consider its
right derived functors called the \emph{higher direct images}.
Moreover as $f^\bullet$ is exact, $f_\bullet$ preserves injectives
and one gets a Grothendieck spectral sequence for the composition
$\Gamma(Y,f_\bullet\mcA)=\Gamma(X,\mcA)$
\[E_2^{pq}=H^p(Y;R^qf_\bullet(\mcA))\Rightarrow H^{p+q}(X;\mcA)\]
In this context it is known as the Leray spectral sequence. In
order to use this spectral sequence one has to have at least some
understanding of the higher direct images $R^nf_\bullet$. Clearly
$R^0f_\bullet=f_\bullet$ and under some conditions one can say
something at least about the stalks $(R^nf_\bullet(\mcA))_y$. Let
$\mcA\rightarrow\mcI^*$ be an injective resolution of $\mcA$ in
$\Sh_X\mcC$. Then $(R^nf_\bullet(\mcA))_y$ is by definition
\[(H^n(f_\bullet\mcI^*))_y\cong
H^n((f_\bullet\mcI^*)_y)\cong
H^n\biggl(\colim_{\substack{U\in\Op(Y)\\y\in
Y}}\mcI^*(f^*U)\biggr)\cong \colim_{\substack{U\in\Op(Y)\\y\in
Y}}H^n(\mcI^*(f^*U))\] and moreover $H^n(\mcI^*(f^*U))\cong
H^n(f^*U;\mcA|f^*U)$. If $f$ is closed (i.e. the image under $f$
of a closed subset is also closed) and the fibres $f^{-1}(y)$ over
points $y$ of $Y$ are ``taut'' then one has an isomorphism (proved
as Proposition IV.4.2.~in \cite{Bredon})
\[\colim_{\substack{U\in\Op(Y)\\y\in Y}}H^n(f^*U;\mcA|f^*U)\cong H^n(f^{-1}(y);\mcA|{f^{-1}(y)})\]

\begin{theorem}
Suppose that $f:X\rightarrow Y$ is a continuous map. Then there is
a convergent spectral sequence
\[E_2^{pq}=H^p(Y;R^qf_\bullet(\mcA))\Rightarrow H^{p+q}(X;\mcA)\]
If moreover $f$ is closed and each $f^{-1}(y)$ taut then there is
an isomorphism
\[(R^qf_\bullet(\mcA))_y\cong H^n(f^{-1}(y);\mcA|{f^{-1}(y)})\]
\end{theorem}

Of course, the stalks of a sheaf do not describe the sheaf
completely unless, for example, they are all zero. In this special
case we get the Vietoris-Begle mapping theorem

\begin{theorem} Suppose that $f:X\rightarrow Y$ is a closed continuous map
such that each $f^{-1}(y)$ is taut and suppose that
$H^n(f^{-1}(y);\mcA|f^{-1}(y))=0$ for all $y\in Y$ and for all
$n>0$. Then there is an isomorphism
\[H^*(Y;f_\bullet\mcA)\cong H^*(X;\mcA)\]
\end{theorem}

Now we will sketch a generalization for the sheaf cohomology with
supports. Let $\psi$ be a family of supports on $X$. There is a
``direct image with supports in $\psi$'' functor
\[f_\psi:\Sh_X\mcC\rightarrow\Sh_Y\mcC\]
Here $f_\psi\mcA$ is defined as the sheafification of the presheaf
\[U\mapsto\Gamma_{\psi\cap f^{-1}U}(f^{-1}U,\mcA)\]
and $f_\psi$ is easily seen to be left exact. If $\varphi$ is a
family of supports on $Y$ then there is another family called
$\varphi(\psi)$ on $Y$ such that
\[\Gamma_\varphi(f_\psi\mcA)=\Gamma_{\varphi(\psi)}(\mcA)\]
Under the additional hypothesis that $\varphi$ is
paracompactifying we get (in the same way but with a lot more work
caused by $f_\psi$ not preserving injectives) a convergent
spectral sequence
\[E_2^{pq}=H^p_\varphi(Y;R^qf_\psi(\mcA))\Rightarrow H^{p+q}_{\varphi(\psi)}(X;\mcA)\]
If for every $F\in\psi$ the image $f(F)$ is closed (i.e. if $f$ is
$\psi$-closed) and if each $f^{-1}(y)$ is $\psi$-taut then we can
identify the stalks as
\[(R^nf_\psi(\mcA))_y\cong H^n_\psi(f^{-1}(y);\mcA|f^{-1}(y))\]
We are particularly interested in the case where $X$ and $Y$ are
locally compact Hausdorff spaces and both $\varphi$ and $\psi$ are
the families of compact supports. In this case $\varphi(\psi)$ is
also the family $c$ of compact supports and clearly every
continuous $f$ is $c$-closed. Also every closed subset of $X$ is
$c$-taut and we have the following results

\begin{theorem} Let $f:X\rightarrow Y$ be a continuous map between
locally compact Hausdorff spaces. Then there is a convergent
spectral sequence
\[E_2^{pq}=H^p_c(Y;R^qf_c(\mcA))\Rightarrow H^{p+q}_c(X;\mcA)\]
Moreover $(R^qf_c(\mcA))_y\cong H^q_c(f^{-1}(y);\mcA|f^{-1}(y))$.
\end{theorem}

\begin{theorem} \label{she_vietoris_begle_compact_supports}
Let $f:X\rightarrow Y$ be a continuous map between locally compact
Hausdorff spaces such that $H^n_c(f^{-1}(y);\mcA|f^{-1}(y))=0$ for
all $y\in Y$ and for all $n>0$. Then there is an isomorphism
\[H^*_c(Y;f_c\mcA)\cong H^*_c(X;\mcA)\]
If moreover $f$ is proper and each $f^{-1}(y)$ connected then
$f=f_c$ and the canonical map $\mcB\xlra{\cong}f_\bullet
f^\bullet\mcB$ is an isomorphism. Consequently for an any
$\mcB\in\Sh_Y\mcC$
\[H^*_c(Y;\mcB)\cong H^*_c(X;f^\bullet\mcB)\]
\end{theorem}



\chapter{Alexander duality, transfer maps} \label{dual_map}

This appendix is concerned with Alexander duality for oriented
topological manifolds. If $f:M\ra N$ is a continuous map between
such manifolds and if $X\subseteq M$ and $Y\subseteq N$ are closed
subsets such that $f^{-1}(Y)\subseteq X$ then we have the
following diagram
 \[\xymatrix{
  H_*(M,M-X) \ar[r]^{f_*} \ar[d]_{\cong} &
  H_*(N,N-Y) \ar[d]^{\cong}
 \\
  {}\check{H}_c^{m-*}(X) \ar[r] &
  {}\check{H}_c^{n-*}(Y)
 }\]
and a natural question arises: is there a reasonable description
of the bottom map? We give a positive answer for the case of
inclusion of a zero section of a vector bundle (which easily leads
to the answer in the case of an embedding of a smooth
submanifold). This map is ``locally'' a cup product with a Thom
class and in this sense generalizes the Thom isomorphism. We start
with a general discussion of singular (co)homology with general
supports and at the end prove the promised Theorem
\ref{dua_main_theorem}.

Let $X$ be a topological space. We say that a family $\varphi$
of closed subsets of $X$ is a family of supports on $X$ if $\varphi$ is
closed under finite unions and taking closed subsets. We say that $\varphi$
is paracompactifying if all elements of $\varphi$ are paracompact and if
any $F\in\varphi$ has a neighbourhood $F\subseteq F'$ such that $F'\in\varphi$.

For any space $X$ we denote by $\phi$ the family of all closed subsets
and by $c$ the family of all compact subsets. Note that $\phi$ is
paracompactifying if $X$ is paracompact and $c$ is paracompactifying if $X$
is locally compact Hausdorff.

If $A\subseteq X$ is a subspace and $\varphi$ a family of supports on $X$
we denote by $\varphi|_A$ the family
\[\varphi|_A=\{F\in\varphi\ |\ F\subseteq A\}\]
It is a family of supports on $A$. Another useful family of supports on $A$ is
\[\varphi\cap A=\{F\cap A\ |\ F\in\varphi\}\]
Note that if $A$ is closed $\varphi|_A=\varphi\cap A$.

We define the (singular) cohomology of $X$ with supports in $\varphi$ as
\[H^*_\varphi(X)=\colim\limits_{F\in\varphi}H^*(X,X-F)\]
The dual construction for homology yields
\[H_*^\varphi(X)=\lim\limits_{F\in\varphi}H_*(X,X-F)\]
For an open subset $U\subseteq X$ we also define the relative versions
\[H^*_\varphi(X,U)=H^*_{\varphi|_{X-U}}(X) \qquad\textrm{and}\qquad
H_*^\varphi(X,U)=H_*^{\varphi|_{X-U}}(X)\]

Note that although the above cohomology groups deserve to be called ``with
supports in $\varphi$'', it is not the case for homology. Obviously
$H_*=H_*^\phi$ (resp. $H^*=H^*_\phi$) are isomorphic to the
singular homology (resp. cohomology) groups while $H^*_c$ is the usual
cohomology with compact supports and $H_*^c$ is the locally finite homology
(whereas the singular homology $H_*$ \textit{does} have compact supports).

\begin{remark}
OK, it is actually not true. The right thing to do (at least in
the case of locally compact spaces) is to take the above limit on
the level of chain complexes and then take homology. If the space
is moreover paracompact then one has a split short exact sequence
\[0\rightarrow\underset{F\in\varphi\;\;}{\lim\nolimits^1}(H_{*+1}(X,X-F))
\rightarrow H_*(\lim_{F\in\varphi}C_*(X,X-F))\rightarrow
\lim_{F\in\varphi}H_*(X,X-F)\rightarrow 0\]
Thus the group that we have defined is a quotient (unnaturally a direct
summand) of the usual one. In what proceeds we only use this $H_*^\varphi$
group in the case of a top dimensional locally finite homology of an orientable
$m$-dimensional manifold $M$ where we have
\[0\rightarrow\underset{F\in c\;\;}{\lim\nolimits^1}(H_{m+1}(M,M-F))
\rightarrow H_m^{lf}(M)\rightarrow H_m^c(M)\rightarrow 0\]
and the $\lim^1$ term vanishes trivially so that the two homology groups
are naturally isomorphic.
\end{remark}

If $f:X\rightarrow Y$ is a continuous map and $\psi$ is a family of
supports on $Y$ we get the induced family of supports $f^{-1}\psi$ on $X$ as
\[f^{-1}\psi=\{F\subseteq X\ |\ \overline{fF}\in\psi\}\]
the closure of $\{f^{-1}F\ |\ F\in\psi\}$ under taking closed subsets. Note
that if $f:X\subseteq Y$ is an inclusion of a subspace then
$f^{-1}\psi=\psi\cap A$.

Obviously if $f^{-1}\psi\subseteq\varphi$ for some family $\varphi$ of
supports on $X$ (or in other words if $F\in\psi\Rightarrow f^{-1}F\in\varphi$,
an obvious generalization of a proper map) we get induced maps
\begin{eqnarray*}
f^*:H^*_{\psi}(Y)\rightarrow H^*_{\varphi}(X) \\
f_*:H_*^{\varphi}(X)\rightarrow H_*^{\psi}(Y)
\end{eqnarray*}

Suppose that $\varphi$ is a paracompactifying family of supports on
$X$ and $A\subseteq X$ a closed subset. Then the inclusion induces a
map $H^*_\varphi(X)\rightarrow H^*_{\varphi|A}(A)$. If $U$ is an
open subset of $X$ then we obtain a map
\[H^*_{\varphi|_U}(U)=\colim_{F\in \varphi|_U}H^*(U,U-F)\cong
\colim_{F\in\varphi|_U}H^*(X,X-F)\rightarrow
\colim_{F\in\varphi}H^*(X,X-F)=H^*_\varphi(X)\]
going in the opposite direction. The isomorphism comes from excision. This is
where we use the assumption on $\varphi$.

Now let $\varphi$ and $\psi$ be two families of supports on $X$.
Then we have cup products
\[\cup:H^*_{\varphi}(X)\otimes H^*_{\psi}(X)\rightarrow
H^*_{\varphi\cap\psi}(X)\]
and cap products
\[\cap:H^*_{\varphi}(X)\otimes H_*^{\varphi\cap\psi}(X)\rightarrow
H_*^{\psi}(X)\]
They are defined as (co)limits of the corresponding products in singular
(co)homology.

A useful formulation of the Thom isomorphism is the following. If
$p:E\rightarrow B$ is an $n$-dimensional oriented vector bundle
there is an orientation class $\tau\in H^*_{c_p}(E)$ such that the
map
\[-\cup\tau:H^*(E)\rightarrow H^{*+n}_{c_p}(E)\]
is an isomorphism, where $c_p$ denotes the family of closed
subsets $C\subseteq E$ for which the restricted projection
$p|_C:C\rightarrow B$ is proper. Obviously if $B$ is compact then
$c_p=c$ and we get an isomorphism (with a shift) between ordinary
cohomology and cohomology with compact supports.

An immediate generalization is the following. Let $\varphi$ be a family of
supports on $B$. Then
\[-\cup\tau:H^*_{p^{-1}\varphi}(E)\rightarrow
H^{*+n}_{(p^{-1}\varphi)\cap c_p}(E)\] is an isomorphism. Also the
projection induces an isomorphism
\[p^*:H^*_{\varphi}(B)\rightarrow H^*_{p^{-1}\varphi}(E)\]
whose inverse is induced by the inclusion $i:B\rightarrow E$ of $B$
as a zero section of $E$ (for reasons that will become clear later we prefer
to use this inclusion instead of the projection map). If one takes $\varphi$
to be $c$, then easily $(p^{-1}c)\cap c_p=c$ and we get a Thom isomorphism
\[H^*_{c}(B)\underset{i^*}{\xlla{\cong}}
H^*_{p^{-1}c}(E)\underset{-\cup\tau}{\xlra{\cong}}H^{*+n}_{c}(E)\]
which is a model for our next proposition.

Dually we have a Thom isomorphism in homology
\[\tau\cap-:H_{*+n}^{(p^{-1}\varphi)\cap c_p}(E)\rightarrow
H_*^{p^{-1}\varphi}(E)\]
and again as a special case when $\varphi=c$
\[H_{*+n}^{c}(E)\underset{\tau\cap-}{\xlra{\cong}}
H_*^{p^{-1}c}(E)\underset{i_*}{\xlla{\cong}}H_*^c(B)\] If moreover
both $E$ and $B$ are oriented manifolds then their fundamental
classes have to correspond under this isomorphism, at least up to
a sign. One can verify that they correspond precisely (not just up
to a sign) if the orientation of the fibres of $p$ is chosen in
such a way that the (local) trvialization $E\cong\bbR^n\times B$
preserves orientations (with the fibre on the right we would have
to introduce a sign).

All of the above applies equally well to the \Cech cohomology
groups: let $X\subseteq M$ be a subspace of a (paracompact)
topological manifold $M$. It is well-known that
\[\check{H}^*(X)\cong\colim_{U}H^*(U)\]
where $U$ ranges over any cofinal system of neighbourhoods of $X$ in $M$.
Let $\varphi$ be a paracompactifying family of supports on $M$.
Then for the \Cech cohomology with supports in $\varphi$
\[\check{H}^*_{\varphi\cap X}(X)\cong\colim_{F\in\varphi}\check{H}^*(X,X-F)\]
we have similarly\footnote{This isomorphism follows from the
following facts: if $\varphi$ is paracompactifying, then it is
also paracompactifying for the pair $(M,X)$ (see Definition
II.9.14.~of \cite{Bredon}; the additional property that any $F\cap
X\in\varphi\cap X$ has a fundamental system of paracompact
neighbourhoods holds as in a metrizable space $M$ every subset is
paracompact). This implies (Proposition II.9.15.~of \cite{Bredon})
that $X$ is $\varphi$-taut (for the \Cech or sheaf cohomology;
again on paracompact spaces they agree), i.e.
\[\check{H}^*_{\varphi\cap X}(X)\cong\colim_U\check{H}^*_{\varphi\cap U}(U)\]
with $U$ running over a fundamental system of neighbourhoods.
Restricting to the cofinal system of open neighbourhoods one can
replace $\check{H}^*_{\varphi\cap U}(U)$ by $H^*_{\varphi\cap
U}(U)$. Here one uses the fact that $\varphi\cap U$ is again
paracompactifying and that $U$ is homologically locally connected
in all dimensions (see Section
\ref{sheaves_relation_to_other_cohomology_theories}). Hence the
canonical map $\check{H}^*_{\varphi\cap U}(U)\rightarrow
H^*_{\varphi\cap U}(U)$ is an isomorphism.}
 \[\check{H}^*_{\varphi\cap X}(X)\cong\colim_UH^*_{\varphi\cap U}(U)\]
with $U$ ranging over any cofinal system of neigbourhoods of $X$.
If $X$ is closed then the closed neighbourhoods of $X$ form such a
cofinal system and in this case
\[\check{H}^*_{\varphi|_X}(X)\cong\colim_UH^*_{\varphi|_U}(U)\]
In particular
\[\check{H}^*_c(X)\cong\colim_UH^*_c(U)\]
Here it is important to restrict to closed neighbourhoods (otherwise the maps
in this system would not have been even defined).

\begin{proposition} Let $p:E^n\rightarrow B^m$ be a vector bundle where both
$E$ and $B$ are oriented (paracompact) topological manifolds and let
$j:X\subseteq E$ be any closed subset. If we think of $B$ as a submanifold of
$E$ via the zero section $i:B\rightarrow E$ then the following diagram commutes
\[\xymatrix{
H_*(B,B-X) \ar[rr]^{i_*} & & H_*(E,E-X) \\
{\check H}^{m-*}_c(B\cap X) \ar[u]_{\cong}^{-\cap o_B} &
{\check H}^{m-*}_{p^{-1}c}(X) \ar[l]^-{k^*} \ar[r]_{-\cup j^*\tau}
\ar[ru]^(.4){-\cap i_*o_B} &
{\check H}^{n-*}_c(X) \ar[u]^{\cong}_{-\cap o_E}
}\]
where $\tau\in H^{n-m}_{c_p}(E)$ is the Thom class of $p$ and
$o_E\in H_n^c(E)$, $o_B\in H_m^c(B)$ are the fundamental classes.
The maps denoted by cap products are the Alexander duality maps and
will be described more closely in the proof.
\end{proposition}

\begin{proof} The map denoted by $-\cap o_E$ in the diagram should be
interpreted as follows. If $U$ is any closed neighborhood of $X$
then $o_E$ gives rise to a class in
\[H_n^c(E,E-X)\cong H_n^c(U,U-X)\]
Let us call this class $\xi_E$. Hence, using the cap product
\[H^{n-*}_c(U)\otimes H_n^c(U,U-X)\xlra{\cap}H_*(U,U-X)\]
we get a homomorphism
\[H^{n-*}_c(U)\xlra{-\cap\xi_E}H_*(U,U-X)\rightarrow H_*(E,E-X)\]
Taking a colimit over all closed neighborhoods $U$ of $X$ we get the map from
the statement. Hence the proof is reduced to the commutativity of the diagram
\[\xymatrix{
H_*(B\cap U,(B\cap U)-X) \ar[rr]^{g_*} & & H_*(U,U-X) \\
H^{m-*}_c(B\cap U) \ar[u]_{\cong}^{-\cap\xi_B} &
H^{m-*}_{p^{-1}c}(U) \ar[l]^{g^*} \ar[r]_{-\cup h^*\tau}
\ar[ru]^(.4){-\cap g_*\xi_B} & H^{n-*}_c(U)
\ar[u]^{\cong}_{-\cap\xi_E} }\] where $g:B\cap U\rightarrow U$ and
$h:U\rightarrow E$ are the inclusions. The part of the diagram on
the left commutes by the naturality of cap products. The
commutativity of the triangle on the right is proved using
\[(x\cup h^*\tau)\cap\xi_E=x\cap(h^*\tau\cap\xi_E)=x\cap g_*\xi_B\]
The second equality is an easy consequence of $\tau\cap
o_E=i_*o_B$, which was observed just before the proposition.
\end{proof}

If $k^*$ is an isomorphism this gives an intrinsic description of the
transfer map
\[i^!=P_E^{-1}i_*P_B:{\check H}^{m-*}_c(B\cap X)\rightarrow
{\check H}^{n-*}_c(X)\]
where $P_B$ and $P_E$ are the two duality maps in the statement of the
proposition.

In the general case let $D_\varepsilon E$ be the $\varepsilon$-disk subbundle
of $E$, $iD_\varepsilon E$ its interior and $S_\varepsilon E$ its boundary,
the $\varepsilon$-sphere bundle. For simplicity we use
\[X_A=A\cap X\]
whenever $A$ is a subspace of $E$. The proposition still holds with $E$
replaced by $iD_\varepsilon E$ and we have a diagram
\[\xymatrix{
H_*(B,B-X) \ar[rr] & & H_*(iD_\varepsilon E,iD_\varepsilon E-X) \\
{\check H}^{m-*}_c(X_B) \ar[u]_{\cong}^{-\cap o_B} &
{\check H}^{m-*}_{p^{-1}c}(X_{iD_\varepsilon E}) \ar[l]_{k_\varepsilon^*}
\ar[d] \ar[r]^{-\cup\tau_\varepsilon'} &
{\check H}^{n-*}_c(X_{iD_\varepsilon E})
\ar[u]^{\cong}_{-\cap o_{iD_\varepsilon E}} \ar[d]_{\cong} \\
& {\check H}^{m-*}_c(X_{D_\varepsilon E})
\ar[r]^-{-\cup\tau_\varepsilon} \ar[lu]^{l_\varepsilon^*} &
{\check H}^{n-*}_c(X_{D_\varepsilon E},X_{S_\varepsilon E}) }\]
with $\tau_\varepsilon'$ and $\tau_\varepsilon$ being the
restrictions of the Thom class
\[\tau\in H^{n-m}_{c_p}(E)\cong
H^{n-m}(D_\varepsilon E,S_\varepsilon E)\]

\begin{lemma} The following diagram commutes
\[\xymatrix{
H_*(iD_\varepsilon E,iD_\varepsilon E-X) \ar[rr] & & H_*(E,E-X) \\
{\check H}^{n-*}_c(X_{iD_\varepsilon E}) \ar[rr]
\ar[u]_{\cong}^{-\cap o_{iD_\varepsilon E}} \ar[d]^{\cong} & &
{\check H}^{n-*}_c(X) \ar[u]^{\cong}_{-\cap o_E} \\
{\check H}^{n-*}_c(X_{D_\varepsilon E},X_{S_\varepsilon E}) &
{\check H}^{n-*}_c(X,X-(iD_\varepsilon E)) \ar[ru] \ar[l]^{\cong}
}\]
where all the unlabeled maps are induced by respective inclusions.
\end{lemma}

\begin{proof} The commutativity of the top square is obtained from the
commutativity of the outer square of
\[\xymatrix{
H_*(U_{iD_\varepsilon E},U_{iD_\varepsilon E}-X) \ar[r]^-{\iota_*} &
H_*(U,U-X) \\
H^{n-*}_c(U_{iD_\varepsilon E}) \ar@{-->}[r]
\ar[u]_{\cong}^{-\cap\xi_{iD_\varepsilon E}} &
H^{n-*}_c(U) \ar[u]^{\cong}_{-\cap\xi_E} \\
H^{n-*}(U_{iD_\varepsilon E},U_{iD_\varepsilon E}-F) \ar[u] &
H^{n-*}(U,U-F) \ar[l]^-{\iota^*}_-{\cong} \ar[u]
}\]
(note that $\xi_E=\iota_*(\xi_{iD_\varepsilon E})$) first by passing to the
colimit over all compact $F\subseteq U_{iD_\varepsilon E}$ to get
commutativity of the upper square and then by passing to the colimit over
all closed neighbourhoods $U$ of $X$. The bottom part of the diagram from the
statement commutes obviously.
\end{proof}

By taking a (co)limit $\varepsilon\rightarrow 0$ both the maps
$k_\varepsilon^*$ and $l_\varepsilon^*$ induce isomorphisms
\[\xymatrix{
{\check H}^{m-*}_c(X_B) & &
\colim{\check H}^{m-*}_{p^{-1}c}(X_{iD_\varepsilon E})
\ar[ll]_-{\{k_\varepsilon^*\}}^-{\cong} \ar[d]^{\cong} \\
& & \colim{\check H}^{m-*}_c(X_{D_\varepsilon E})
\ar[llu]^{\{l_\varepsilon^*\}}_{\cong}
}\]
and we get at least some (intrinsic) description of the transfer map by
pasting the two diagrams together.

\begin{theorem} \label{dua_main_theorem}
The transfer map $i^!:{\check H}^{m-*}_c(B\cap X)\rightarrow
{\check H}^{n-*}_c(X)$ is the colimit of
\begin{multline*}
{\check H}^{m-*}_c(X_B)\xlla{l_\varepsilon^*} {\check
H}^{m-*}_c(X_{D_\varepsilon E})\xlra{-\cup\tau_\varepsilon}
{\check H}^{n-*}_c(X_{D_\varepsilon E},X_{S_\varepsilon
E})\xlla{\cong}
\\
\xlla{\cong}{\check H}^{n-*}_c(X,X-(iD_\varepsilon E))\lra {\check
H}^{n-*}_c(X)
\end{multline*}
where the first map becomes an isomorphism in the colimit. \qed
\end{theorem}



\chapter{Transversality of maps from preimages of other maps}

This appendix is rather unrelated to the rest of the thesis. It
seemed to be important for the proof of the main theorem but
turned out not to be. It gives an equivalent condition to a
transversality of a restriction of a fixed map $g:P\ra M$ to a
preimage $f^{-1}(A)$ of a submanifold along a map $f:P\ra N$ in
terms of transversality conditions of $f$ itself. An the end we
prove that this property is generic in the sense that such maps
form a residual subset.

\begin{lemma} \label{tra_equivalence1} Let
\[\xymatrix{
& f^{-1}(A) \ar[r] \ar@{c->}[d]^-j &
A \ar@{c->}[d] \\
g^{-1}(B) \ar@{c->}[r]^-i \ar[d] & P \ar[r]^{f} \ar[d]^{g} & M \\
B \ar@{c->}[r] & N}\] be a diagram of smooth manifolds and smooth
maps between them where the maps denoted by
$\xymatrix@1{{}\ar@{c->}[r]&{}}$ are inclusions of submanifolds
and where we assume that $f\pitchfork A$ and $g\pitchfork B$. Then
the following conditions are equivalent:
\begin{enumerate}
\item[(i)]{$gj\pitchfork B$}
\item[(ii)]{$fi\pitchfork A$}
\item[(iii)]{$f^{-1}(A)\pitchfork g^{-1}(B)$}
\item[(iv)]{$(f,g)\pitchfork A\times B$ where $(f,g):P\rightarrow M\times N$}
\end{enumerate}
\end{lemma}

\begin{proof} Let us start first by showing that (i), (ii) and
(iii) are equivalent. Because (iii) is symmetric it is enough to
show the equivalence of (i) and (iii). But (i) is equivalent to
the map
\[T_x f^{-1}(A)\longrightarrow T_{g(x)}N/T_{g(x)}B
\] induced by $g$ being surjective for every $x\in f^{-1}(A)\cap
g^{-1}(B)$. Because of the assumption $g\pitchfork B$, we have a
commutative diagram
\[\xymatrix{
T_x f^{-1}(A) \ar[r] \ar[d] & T_{g(x)}N/T_{g(x)}B \\
T_x P/T_x g^{-1}(B) \ar[ru]_\cong}\] and so (i) is equivalent to
(iii).

Let us deal with (iv) now (which we actually do not need in the
following). It can be phrased as a surjectivity of the map
(induced by $(f,g)$)
\[\varphi:T_xP\longrightarrow (T(M\times N)/T(A\times B))_{(f(x),g(x))}
\cong (TM/TA)_{f(x)}\oplus(TN/TB)_{g(x)}\] for every $x\in
f^{-1}(A)\cap g^{-1}(B)$. Assuming (i) and (ii) the image under
$\varphi$ of $T_xf^{-1}(A)$ is $0\oplus(TN/TB)_{g(x)}$ and the
image of $T_xg^{-1}(B)$ is $(TM/TA)_{f(x)}\oplus 0$ and so
$\varphi$ is surjective. Assuming the surjectivity of $\varphi$ on
the other hand we can find $u\in T_xP$ mapping to $(0,v)$ by
$\varphi$ for any choice of $v$. Necessarily $u\in T_xf^{-1}(A)$
and so the map
\[T_x f^{-1}(A)\longrightarrow T_{g(x)}N/T_{g(x)}B
\] is surjective. This is condition (i).
\end{proof}

\begin{construction} Let $D$ be a $d$-dimensional manifold and
\[\Diff(\bbR^d,0)=G_{0,\mathrm{diff}}(\bbR^d,\bbR^d)_0\]
the group\footnote{If we wanted to give $\Diff(\bbR^d,0)$ a
topology we could do so by inducing the topology via the map
$\Diff(\bbR^d,0)\ra J^\infty_0(\bbR^d,\bbR^d)_0$. Or one could
replace $\Diff(\bbR^d,0)$ by its image - the subspace of
invertible $\infty$-jets.} of germs at 0 of local diffeomorphisms
$\bbR^d\rightarrow\bbR^d$ fixing 0. Define a principal
$\Diff(\bbR^d,0)$-bundle
\[\Charts_D=G_{0,\mathrm{diff}}(\bbR^d,D)\xlra{\mathrm{ev}_0} D\]
of germs at 0 of local diffeomorphisms $\bbR^d\rightarrow D$. If
$F$ is any manifold with a (smooth in some sense) action of
$\Diff(\bbR^d,0)$ then we can construct an associated bundle
\[D[F]:=\Charts_D\times_{\Diff(\bbR^d,0)}F\longrightarrow D\]
Any bundle of this form will be called \textit{local}. Observe
that this construction is functorial in $D$ on the category of
$d$-dimensional manifolds and local diffeomorphisms. As an example
the bundle $J^r(D,N)$ of $r$-jets of maps $D\rightarrow N$ is a
local bundle as
\begin{equation}\label{tra_jet_bundle}
D[J^r_0(\bbR^d,N)]=
\Charts_D\times_{\Diff(\bbR^d,0)}J^r_0(\bbR^d,N)\cong J^r(D,N)
\end{equation} where $J^r_0(\bbR^d,N)$ is the subspace of
$J^r(\bbR^d,N)$ of $r$-jets with source 0. The bijection is
provided by the map
\[[u,\alpha]\mapsto\alpha\circ j^r_{u(0)}(u^{-1})\]

Having a $\Diff(\bbR^d,0)$-invariant submanifold $B\subset
J^r_0(\bbR^d,N)$ we get an associated subbundle $D[B]\subset
J^r(D,N)$ for any $d$-dimensional manifold $D$. This allows us to
talk about jet transversality conditions on a map $D\rightarrow N$
without specifying what $D$ (and hence also $J^r(D,N)$) is. Let us
fix such a $B\subset J^r_0(\bbR^d,N)$.

Using (\ref{tra_jet_bundle}) one can easily see that $j^r_x(h)\in
J^r_x(D,N)$ lies in $D[B]$ iff in some (and hence any) chart $u\in
G_{0,\mathrm{diff}}(\bbR^d,D)_x$ the local expression $j^r_0(hu)$
of $j^r_x(h)$ in the chart $u$ lies in $B$.
\end{construction}

\begin{lemma} \label{tra_equivalence2} For $h:D\rightarrow N$ the following conditions
are equivalent
\begin{enumerate}
\item[(i)]{$h_*:J^r_{0,\mathrm{imm}}(\bbR^d,D)
\rightarrow J^r_0(\bbR^d,N)$ is transverse to $B$}
\item[(ii)]{$j^r(h):D\rightarrow J^r(D,N)$ is transverse to $D[B]$}
\end{enumerate}
where the \emph{``$\mathrm{imm}$''} index means we take the subspace
of jets of immersions.
\end{lemma}

\begin{proof} Taking the associated bundles (i) is obviously
equivalent to the transversality of
\[h_*:J^r_\mathrm{imm}(D,D)\rightarrow J^r(D,N)\]
to $D[B]$. Let $j^r_x(k)\in J^r_{x,\mathrm{imm}}(D,D)_y$ be an
$r$-jet of a diffeomorphism $k:V\rightarrow W$ between open
subsets of $D$. Then we have a diagram
\[\xymatrix{
J^r_{\mathrm{imm}}(V,D) \ar[r]^{h_*} \ar[d]_{\cong}^{k_*} &
J^r(V,N) \ar[d]_{\cong}^{k_*} & V[B] \ar@{d->}[l] \ar[d]_{\cong}^{k_*}\\
J^r_{\mathrm{imm}}(W,D) \ar[r]^{h_*} & J^r(W,N) & W[B]
\ar@{d->}[l] }\] Now $j^r_x(k)$ in the top left corner is mapped
by $k_*$ down to $j^r_y(\id)$. Hence we see that it is enough
(equivalent) to check the transversality only at $j^r_y(\id)$'s
for all $y\in D$ for which $h_*(j^r_y(\id))=j^r_y(h)\in D[B]$. For
such $y$ the same diagram shows that every $j^r_x(k)$ with target
$y$ is mapped by $h_*$ to $D[B]$. Thus the whole fibre over $y$ of
the target map
\[J^r_\mathrm{imm}(D,D)\xlra{\tau} D\]
is mapped to $D[B]$. The target map $\tau$ has a section
\[j^r(\id):D\rightarrow J^r_\mathrm{imm}(D,D)\]
and so (i) is finally equivalent to the composite
\[D\xlra{j^r(\id)}J^r_\mathrm{imm}(D,D)\xlra{h_*}J^r(D,N)\]
being transverse to $D[B]$. This is (ii).
\end{proof}

We say that a map $g:P\rightarrow N$ is transverse to $B$, denoted
$g\pitchfork B$, if
\[g_*:J^r_{0,\mathrm{imm}}(\bbR^d,P)\rightarrow J^r_0(\bbR^d,N)\]
is transverse to $B$. This is the case for example if $g$ is a
submersion. When $r=0$ this is equivalent to the usual
transversality of a map to a submanifold. Let $f\pitchfork A$
where $f:P\rightarrow M$ and $i:A\subset M$ is a submanifold. Then
we have the following diagram
\[\xymatrix{
J^r_{0,\mathrm{imm}}(\bbR^d,f^{-1}(A)) \ar@{c->}[r] \ar[d]^{j_*} &
J^r_0(\bbR^d,f^{-1}(A)) \ar[r] \ar[d]^{j_*} & J^r_0(\bbR^d,A) \ar[d]^{i_*} \\
J^r_{0,\mathrm{imm}}(\bbR^d,P) \ar@{c->}[r] & J^r_0(\bbR^d,P)
\ar[r]^-{f_*} & J^r_0(\bbR^d,M)}\] where both squares are
pullbacks. This can be easily seen in local coordinates. Also
$i_*$ is a submanifold inclusion and $f_*\pitchfork i_*$.
Combining Lemma \ref{tra_equivalence1} with Lemma
\ref{tra_equivalence2} we get:

\begin{lemma} \label{tra_equivalence3}
Let
\[\xymatrix{
f^{-1}(A) \ar[r] \ar@{c->}[d]^-j & A \ar@{c->}[d] \\
P \ar[r]^{f} \ar[d]^{g} & M \\
N}\] be a diagram of smooth manifolds and smooth maps between them
where the maps denoted by $\xymatrix@1{{}\ar@{c->}[r]&{}}$ are
inclusions of  submanifolds. Let us assume that $f\pitchfork A$
and $g\pitchfork B$ where $B\subset J^r_0(\bbR^d,N)$ is a
$\Diff(\bbR^d,0)$-invariant submanifold with $d=\dim P+\dim A-\dim
M$. Then the following conditions are equivalent:
\begin{enumerate}
\item[(i)]{$j^r(gj)\pitchfork f^{-1}(A)[B]$, where
\[j^r(gj):f^{-1}(A)\rightarrow J^r(f^{-1}(A),N)\] is the jet
prolongation}
\item[(ii)]{$f_*|_Y\pitchfork J^r_0(\bbR^d,A)$}, where
$Y=(g_*)^{-1}(B)$ is defined by a pullback diagram
\[\xymatrix{
Y \ar@{c->}[r] \ar[d] & J^r_{0,\mathrm{imm}}(\bbR^d,P) \ar[d]^{g_*} \\
B \ar@{c->}[r] & J^r_0(\bbR^d,N)}\]
\end{enumerate}
\end{lemma}

\begin{proof} Applying Lemma \ref{tra_equivalence1} to the diagram
\[\xymatrix{
& J^r_{0,\mathrm{imm}}(\bbR^d,f^{-1}(A)) \ar[r] \ar[d]^{j_*} & J^r_0(\bbR^d,A) \ar[d]^{i_*} \\
Y \ar@{c->}[r] \ar[d] & J^r_{0,\mathrm{imm}}(\bbR^d,P) \ar[r]^-{f_*} \ar[d]^{g_*} & J^r_0(\bbR^d,M)\\
B \ar@{c->}[r] & J^r_0(\bbR^d,N)}\] gives an equivalence of (ii)
with transversality of
\[(gj)_*:J^r_{0,\mathrm{imm}}(\bbR^d,f^{-1}(A))\longrightarrow J^r_0(\bbR^d,N)\]
to $B$. By Lemma \ref{tra_equivalence2} this is equivalent to (i).
\end{proof}

Now that we know how $f$ controls the transversality of a map
defined on the preimage $f^{-1}(A)$ of some submanifold, we would
like to see that this transversality condition (any of the two
equivalent conditions in Lemma \ref{tra_equivalence3}) is generic.
This is indeed the case. We first prove a more general result
which at the same time happens to generalize the Thom
Transversality Theorem.

\begin{lemma} \label{tra_transversality_theorem} Let $D$, $M$, $N$ be manifolds, $Y\subset
J^r_\mathrm{imm}(D,M)$ and $Z\subset J^r(D,N)$ submanifolds. Let
us further assume that $\sigma_Y\pitchfork\sigma_Z$, where
\[\sigma_Y=\sigma|_Y:Y\subset J^r(D,M)\rightarrow D\quad\textrm{and}\quad
\sigma_Z=\sigma|_Z:Z\subset J^r(D,N)\rightarrow D\] are the
restrictions of the source maps. For a smooth map $f:M\rightarrow
N$ let $f_*|_Y$ denote the map
\[Y\subset J^r_\mathrm{imm}(D,M)\xlra{f_*} J^r(D,N)\] Then the set
\[\mfX:=\bigl\{f\in C^\infty(M,N)\ \bigl|\ f_*|_Y\pitchfork Z\bigl\}\]
is residual in $C^\infty(M,N)$ with the strong topology, and open
if $Z$ is closed (as a subset) and $\tau_Y:Y\rightarrow M$ proper.
\end{lemma}

\begin{proof} This is an application of Lemma \ref{ttt_basic_lemma}.
We have a map
\[\alpha:C^\infty(M,N)\rightarrow C^\infty(Y,J^r(D,N))\]
sending $f$ to $f_*|_Y$. This map is continuous in the weak
topologies
. Clearly $\mfX=\{f\in C^\infty(M,N)\ |\ \alpha(f)\pitchfork Z\}$.
We have to verify the conditions of Lemma \ref{ttt_basic_lemma}.

Let $f_0\in C^\infty(M,N)$ and $K\subseteq Y$, $L\subseteq Z$
compact disks. We can assume that $\tau(K)$ lies in a coordinate
chart $\bbR^m\cong U\subseteq M$ and that $\tau(L)$ lies in a
coordinate chart $\bbR^n\cong V\subseteq N$. We use these charts
to identify $U$ with $\bbR^m$ and $V$ with $\bbR^n$ when needed.
Let $\lambda:U\rightarrow\bbR$ be a compactly supported function
such that $\lambda=1$ on a neighborhood of $\tau(K)\cap
f_0^{-1}(\tau(L))$ and such that $\lambda=0$ on $U-f_0^{-1}(V)$.

We set $P:=J^r_*(\bbR^m,\bbR^n)\cong J^r_*(U,V)$ and identify it
with the space of polynomial mappings $U\ra V$. Then we get a map
\[\beta:P\rightarrow C^\infty(M,N)\]
sending $g$ to the function $f_0+\lambda g$ where the operations
are interpreted inside $V$ via the chart. It is continuous in the
strong topology
and the map
\[\gamma=(\alpha\beta)^\sharp:P\times Y\rightarrow J^r(D,N)\]
is smooth. Thus it is enough to show that (after a suitable
restriction) $\gamma\pitchfork Z$. Clearly $\gamma$ sends
$(g,j^r_x(h))$ to $j^r_x((f_0+\lambda g)h)$. Suppose now that
$h(x)\in W:=\textrm{int}\lambda^{-1}(1)$ so that this equals to
$j^r_x(f_0h+gh)$. By restriction we get a map
\[\delta:P\cong P\times\{j^r_x(h)\}\xlra{\gamma}J^r_x(D,V)\]
In the affine structure on $J^r_x(D,V)$ inherited from the chart
$\delta$ is clearly affine. Identifying $P$ with $J^r_{h(x)}(U,V)$
the linear part of $\delta$ is just a precomposition with $h$
\[h^*:J^r_{h(x)}(U,V)\rightarrow J^r_{x}(D,V)\]
The map $h$, being an immersion, has (locally - near $x$) a left
inverse $\pi$ which then gives a right inverse $\pi^*$ of $h^*$
and so the linear part of $\delta$ is surjective and hence it is a
submersion.

In the horizontal direction our transversality condition
$\sigma_Y\pitchfork\sigma_Z$ applies and so $\gamma\pitchfork Z$
on $P\times \tau_Y^{-1}(W)$. If $f:M\rightarrow N$ is close enough
to $f_0$ then
\[\tau(K)\cap f^{-1}(\tau(L))\subseteq W\]
and in particular there is a neighbourhood $P'$ of $0$ in $P$ such
that $\beta(P')$ consists only of such maps. Therefore the
restriction of $\gamma$ to
\[P'\times K\lra J^r(D,N)\]
is transverse to $Z$.

If $\tau_Y$ happens to be proper then $\alpha$ is continuous even
in strong topologies and $\mfX$ is a preimage of the open subset
of maps $f:Y\rightarrow J^r(D,N)$ transverse to $Z$.
\end{proof}

\begin{corollary}[Thom Transversality Theorem] Let $M$, $N$ be
manifolds, $Z\subseteq J^r(M,N)$ a submanifold. Then the set
\[\mfX:=\{f\in C^\infty(M,N)\ |\ j^rf\pitchfork Z\}\]
is residual in $C^\infty(M,N)$. If $Z$ is closed (as a subset)
then it is also open.
\end{corollary}

\begin{proof} We apply Lemma
\ref{tra_transversality_theorem} to $D=M$ and
\[\xymatrix{
Y=M \ar@{c->}[rr]^{j^r(\id)} \ar[rd]_{\id} & & J^r(M,M)
\ar[ld]^{\sigma}
\\
& M}\] As $\sigma_Y=\id=\tau_Y$ it is both proper and transverse
to $\sigma_Z$ for any $Z$.
\end{proof}

\begin{corollary} Let $M$, $N$ be manifolds, $Y\subset
J^r_{0,\mathrm{imm}}(\bbR^d,M)$ and $Z\subset J^r_0(\bbR^d,N)$
submanifolds. Then the set
\[\bigl\{f\in C^\infty(M,N)\ \bigl|\ f_*|_Y\pitchfork Z\bigl\}\]
is residual in $C^\infty(M,N)$ with the strong topology, and open
if $Z$ is closed (as a subset) and $\tau_Y:Y\ra M$ proper.
\end{corollary}

\begin{proof} Under the natural identifications
\[\begin{split}
\bbR^d\times J^r_{0,\mathrm{imm}}(\bbR^d,M)\cong & J^r_{\mathrm{imm}}(\bbR^d,M)\\
\bbR^d\times J^r_0(\bbR^d,N)\cong & J^r(\bbR^d,N)
\end{split}\]
we can apply Lemma \ref{tra_transversality_theorem} to $D=\bbR^d$,
$M$, $N$ and
\[\begin{split}
0\times Y\subseteq & J^r_{\mathrm{imm}}(\bbR^d,M)\\
\bbR^d\times Z\subseteq & J^r(\bbR^d,N)
\end{split}\]
\end{proof}

%





\backmatter


\end{document}